\documentclass[12pt]{preprint} 

\usepackage[a4paper,innermargin=1.2in,outermargin=1.2in,
bottom=1.5in,marginparwidth=1in,marginparsep
=3mm]{geometry}

\usepackage[T1]{fontenc} 
\usepackage{fourier} 
\usepackage[english]{babel} 
\usepackage{amsmath,amsfonts,amsthm} 
\usepackage{tikz-cd}

\usepackage{shuffle}
\usepackage{ amssymb }
\usepackage{fancyhdr} 
\usepackage{comment}
\usepackage{stmaryrd}

\usepackage{cprotect}
\usepackage{hyperref}
\usepackage{mathrsfs}
\usepackage{mhequ}
\usepackage{microtype}
\usepackage{esint}

\usepackage{xstring}

\usepackage{orcidlink}
\usepackage{longtable}
\usepackage{tabularx}
\usepackage{booktabs}

\colorlet{darkblue}{blue!90!black}
\colorlet{darkred}{red!90!black}

\newcommand{\fraks}{\mathfrak{s}}

\newcommand{\cC}{\mathcal{C}}

\newcommand{\R}{\mathbb{R}}

\newcommand{\mcL}{\mathcal{L}}

\renewcommand{\bf}{\mathbf{f}}

\newcommand{\eps}{\varepsilon}

\newcommand{\E}{\mathbb{E}}

\newcommand{\bfK}{\pmb{\mathcal{K}}}

\def\L{\operatorname{L}}

\newcommand{\vertiii}[1]{{\left\vert\kern-0.25ex\left\vert\kern-0.25ex\left\vert #1 
		\right\vert\kern-0.25ex\right\vert\kern-0.25ex\right\vert}}

\newtheorem{theorem}{Theorem}[section]

\newtheorem{corollary}{Corollary}[theorem]

\newtheorem{lemma}[theorem]{Lemma}

\newtheorem{prop}[theorem]{Proposition}
\newtheorem{assumption}[theorem]{Assumption}
\theoremstyle{remark}
\newtheorem{remark}[theorem]{Remark}

\DeclareMathOperator{\supp}{supp}

\numberwithin{equation}{section} 
\numberwithin{figure}{section} 
\numberwithin{table}{section} 

\newcommand{\overbar}[1]{\mkern 1.5mu\overline{\mkern-1.5mu#1\mkern-1.5mu}\mkern 1.5mu}
\def\bsquare{\overbar{\square}}

\colorlet{symbols}{blue!90!black}
\colorlet{testcolor}{green!60!black}

\usetikzlibrary{calc}

\tikzset{
	eps/.style={circle,fill=white,draw=symbols,inner sep=0pt,minimum size=0.8mm},
	}
\makeatletter
\def\DeclareSymbol#1#2#3{\expandafter\gdef\csname MH@symb@#1\endcsname{\tikz[baseline=#2,scale=0.15,draw=symbols]{#3}}}
\def\<#1>{\csname MH@symb@#1\endcsname}
\makeatother

\DeclareSymbol{1}{0}{\draw[white] (-.5,0) -- (.5,0); \draw (0,0)  -- (0,1.5) node[eps] {};}
\DeclareSymbol{2}{0}{\draw (-0.5,1.5) node[eps] {} -- (0,0) -- (0.5,1.5) node[eps] {};}
\DeclareSymbol{3}{0}{\draw (0,0) -- (0,1.5) node[eps] {}; \draw (-.7,1.3) node[eps] {} -- (0,0) -- (.7,1.3) node[eps] {};}

\makeatletter

\DeclareRobustCommand{\TitleEquation}[2]{\texorpdfstring{\StrLeft{\f@series}{1}[\@firstchar]$\if%
		b\@firstchar\boldsymbol{#1}\else#1\fi$}{#2}}

\makeatother

\newcommand{\horrule}[1]{\rule{\linewidth}{#1}} 

\title{	
\horrule{0.5pt} \\[0.4cm] 
\large Periodic space-time homogenisation of the \TitleEquation{\phi^4_2}{phi^4_2} equation
 \\ 
\horrule{0.5pt} \\[0.5cm] 
}
\author{Martin Hairer$^{1,2}$\orcidlink{0000-0002-2141-6561} and Harprit Singh$^3$\orcidlink{0000-0002-9991-8393}}
\institute{EPFL, Lausanne Switzerland, \email{martin.hairer@epfl.ch} \and 
Imperial College London, UK, \email{m.hairer@imperial.ac.uk} \and
University of Edinburgh, UK, \email{harprit.singh@ed.ac.uk}}

\begin{document}

\maketitle 
\begin{abstract}
We consider the homogenisation problem for the $\phi^4_2$ equation on the torus $\mathbb{T}^2$, namely the behaviour as $\eps \to 0$ of the solutions to the equation 
\textit{suggestively} written as
\begin{equ}
\partial_t u_\eps -  \nabla\cdot  {A}(x/\eps,t/\eps^2) \nabla u_\eps = -u^3_\eps +\xi
\end{equ}
where 
$\xi$ denotes space-time white noise and $A: \mathbb{T}^2\times  \mathbb{R}$ is uniformly elliptic, periodic
and H\"older continuous.  When the noise is regularised at scale $\delta \ll 1$ we show that any joint limit $\eps,\delta \to 0$ recovers the
classical dynamical $\phi^4_2$ model. In certain regimes or if 
the regularisation is chosen in a specific way adapted to the 
problem, we show that the counterterms can be chosen as explicit local functions of $A$.
\end{abstract}
\tableofcontents

\section{Introduction}\label{sec:intro}

%
For $\bar A\in \mathbb{R}^{2\times 2}$ a strictly positive definite matrix, the constant coefficient $\phi^4_2$ equation
\begin{equ}
\partial_t u -  \nabla\cdot  \bar{A} \nabla u = -u^3  +\xi\ ,
\end{equ}
(up to a change of variables) was considered in \cite{DD03}. It was observed that, just as in the
construction of the $\phi^4_2$ measure \cite{GlimmJaffe,Simon}, in order to obtain a meaningful notion of solution one should instead consider the renormalised equation formally given by
\begin{equation}\label{eq:phi^4_2_const_intro}
\partial_t u -  \nabla\cdot  \bar{A} \nabla u = -u^3 + \infty \cdot  u  +\xi\ ,
\end{equation}
see also \cite{HRW12}. 
%
In \cite{Sin23} it was noted that in the non-translation invariant setting, i.e.\ when {$A\in C^\alpha( \mathbb{R}\times \mathbb{T}^2, \mathbb{R}^{2\times 2})$} for some $\alpha>0$, a meaningful notion of solution corresponds to considering
\begin{equation}\label{eq:phi^4_2_intro}
\partial_t u-  \nabla\cdot  {A}(x,t) \nabla u = -u^3 + \infty \cdot D(x,t) u +\xi,
\end{equation}
where $D= (\det{A}^s)^{-1/2}$ for $A^s$ the symmetric part of $A$, see also \cite[Section~15]{HS23m} for a general discussion in this direction. 

%

In this article, we consider the periodic homogenisation problem for the $\phi^4_2$ equation.
Assuming that $A(x,t)$ is space-time periodic, we consider the equations formally given by
\begin{equ}\label{eq:homogenisation problem}
\partial_t u_\eps -  \nabla\cdot  {A}(x/\eps,t/\eps^2) \nabla u_\eps = -u^3_\eps + \infty \cdot D(x/\eps, t/\eps^2) +\xi
\qquad \text{for} \qquad \eps^{-1}\in \mathbb{N} \ .
\end{equ}
Periodic homogenisation of elliptic and parabolic equations is a well studied subject, see for example \cite{sanchez1980non, BLP78, jikov2012homogenization, bakhvalov2012homogenisation}.
The results therein imply that the oscillatory operator
$\mcL_\eps = \nabla\cdot  {A}(x/\eps,t/\eps^2) \nabla $ converges in the resolvent sense to a homogenised operator $\nabla\cdot  \bar{A} \nabla$, see \eqref{eq:hom_matrix} for an
expression for $\bar{A}$ in our setting.

In particular, these results suggest that as $\eps\to 0$ the solutions $u_\eps$ to \eqref{eq:homogenisation problem}
converge to a solution of \eqref{eq:phi^4_2_const_intro} for $\bar{A}$ the homogenised matrix associated to $A$.
The main results of this article state that instead an \textit{additional} constant renormalisation is required.
More precisely, we show in Theorem~\ref{thm:main1} below that for $\eps,\delta>0$ one can find constants $\alpha_{\eps,\delta}, \bar{\alpha}_{\eps,\delta}\in \mathbb{R}$ 
such that the functions
\begin{equ}
(0,1]^2 \ni (\eps,\delta) \mapsto\alpha_{\eps,\delta}-
\frac{|\log(\delta/\eps)\wedge 0|}{4\pi} 
, \qquad 
(0,1]^2 \ni
(\eps,\delta) \mapsto
\bar{\alpha}_{\eps,\delta}-
  \frac{|\log(\delta)| \wedge |\log(\eps)|}{{4\pi} },
\end{equ}
extend continuously to $\eps,\delta\in [0,1]^2$
and a bounded family of constants ${c_{\eps,\delta}}$ such that the solution 
$u_{\eps, \delta}$ to
\begin{equ}
\partial_t u_{\eps, \delta} -\mathcal{L}_\eps u_{\eps, \delta}= -u_{\eps,\delta}^3+ 3 \alpha_{\eps, \delta}  D(x/\eps,t/\eps^2)u_{\eps,\delta} + 3\bar{\alpha}_{\eps,\delta} \det(\bar {A})^{-1/2}  u_{\eps,\delta} + 3c_{\eps,\delta} u_{\eps,\delta} +\xi_{\eps, \delta}\  , 
\end{equ}
converges to that of \eqref{eq:phi^4_2_intro} as $(\eps,\delta) \to 0$. Let us also point out that the function $(0,1]^2\ni (\eps,\delta) \to c_{\eps,\delta}\in\mathbb{R}$ can not be extended continuously to all of $ [0,1]^2$. Let us also mention that
here, $\xi_{\eps, \delta}$ denotes
a rather specific regularisation of the white noise $\xi$, namely it is regularised by the heat kernel
associated to $\mathcal{L}_\eps$ at time $\delta^2$.
This particular choice is convenient since it allows us to provide an explicit expression for the renormalisation
which covers every possible regime regarding the relative sizes of $\eps$ and $\delta$. 
Our method also applies for usual (translation invariant) mollification of the noise, but one only retains explicit local counter terms in the regime $\eps\lesssim\delta$, see Theorem~\ref{thm:main_translation invariant} and Proposition~\ref{prop:main_flat}. Similarly,
Theorem~\ref{thm:main_non} uses the general class of mollification schemes in \cite[Sec.~3.2]{Sin23} and treats the regime $\delta\lesssim \eps$.
Importantly, in all these results the solutions for $\delta=0$ agree and fall within in the canonical class of solutions constructed in \cite{Sin23}.

As this article was nearing completion, we learned of \cite{Weijun} which considers a very
similar situation. The main\footnote{ Another difference is that we treat space-time homogenisation on the torus while \cite{Weijun} consider spatial homogenisation on regular domains for vanishing boundary condition.}  differences between the two works are our quantitative treatment of renormalisation to ensure that for any $\eps>0$ the notion of solution used is consistent with the solution theory for non-constant coefficient equations in \cite{Sin23} and that we consider joint limits $(\eps,\delta)\to 0$.
While that and the present article are the first results on homogenisation of singular SPDEs, 
 at least to the best of our knowledge, homogenisation of stochastic PDEs that are either more regular or linear 
has been considered for example in \cite{HK12, ichihara2004homogenization, PARDOUX20031, armstrong2022quantitative}. 
For the topic of stochastic homogenisation which is related but quite distinct, we refer to \cite{armstrong2019quantitative, josien2022annealed,armstrong2022elliptic, gloria2020regularity}.

\paragraph{Outline of the article}
In order to establish our results, we follow the ansatz of \cite{DD03} to decompose the solution into (essentially) the solution of the linear equation $\<1>_{\eps,\delta}$ and a remainder. This splits the required arguments into two main steps:
\begin{enumerate}
\item\label{item:intro1} Consider the homogenisation problem for the remainder equation driven by arbitrary distributions in $C^{-\kappa}$ for some $\kappa>0$.
\item\label{item:intro2} Construct (renormalised) polynomials of $\<1>_{\eps,\delta}$.
\end{enumerate}
The article is structured as follows: In Section~\ref{sec:intro} the function spaces used in the article are introduced and classical background on heat kernels for non-constant coefficient differential operators as well as on homogenisation theory is recalled.  Subsection~\ref{subsec:main results}, then, presents the main results of the article, to the proof of which the remainder of the article is dedicated. In Section~\ref{sec:homogenisation estimates} homogenisation estimates on the kernel which we shall need are derived. With these at hand Step~\ref{item:intro1} above is carried out in Section~\ref{sec:fixedpoint}. Section~\ref{sec:stochastic_estimates} contains Step~\ref{item:intro2} for the case of heat kernel regularisations of the noise in full detail. Section~\ref{sec:further stoch est} carries out the same step for the other regularisation, in order to not be too repetitive this section is somewhat streamlined. 
Finally, in Section~\ref{sec:proof of main} the two steps above are combined to conclude the main results. 
Some lemmata and Propositions used in the main text were moved to Appendices~\ref{ap:A},~\ref{ap:B}~\&~\ref{ap:C} in order not to interrupt the flow of the main arguments. Lastly, Appendix~\ref{ap:D} contains an index of the several kernels used throughout this article.

The general idea to use classical homogenisation theory to first establish heat kernel estimates, which can then be used in SPDE arguments is clearly not restricted to the equation considered here. In forthcoming work we shall implement this approach for more singular subcritical SPDEs within the framework of regularity structures. Instead of using estimates on the heat kernel, it might be possible to obtain these results using regularity estimates directly, c.f.~\cite{MR3925533} in the context of quasilinear singular SPDEs.

\paragraph{Acknowledgments}
HS would like to thank Felix Otto and Kirill Cherednichenko for valuable discussions and gratefully acknowledges financial support from the EPSRC via Ilya Chevyrev’s New Investigator Award EP/X015688/1.
MH was partially supported by the Royal Society through the research professorship
RP\textbackslash R1\textbackslash 191065.

%

\subsection{Function Spaces}
We first introduce the function spaces used in this article. We shall always identify $\mathbb{T}^d= \mathbb{R}^d/\mathbb{Z}^d$. 
\subsubsection{Function spaces on \TitleEquation{\mathbb{T}^d}{T^d}}
For $\gamma\in (0,1)$ and $B\subset \mathbb{R}^d$, define for $f:\mathbb{R}^d \to \mathbb{R}$
\begin{equ}
\|f\|_{C^\gamma (B)}= \sup_{x,\zeta\in B}  \frac{|f(x)- f(\zeta)| }{|x-\zeta|^\gamma} \ .
\end{equ}
For $\rho: \mathbb{R}^d \to \mathbb{R}$ we shall write
$\rho^{\lambda}_x= \lambda^{-d}  \rho\left(\frac{\cdot\ - x}{\lambda} \right)$.
Then, for a distribution $F\in \mathcal{D}'(\mathbb{R}^d)$ and $-1<\gamma<0$ let
\begin{equ}
\| F\|_{C^\gamma(B)} = \sup_{x\in B} \sup_{\lambda\in (0,1)}    \sup_{\rho\in \mathfrak{B}^{1}(\mathbb{R}^{d})} \frac{|\langle F, \rho^{\lambda}_x \rangle|}{\lambda^\gamma} \ ,
\end{equ}
where 
\begin{equ}
\mathfrak{B}^{1}(\mathbb{R}^{d})= \{ \rho \in \mathcal{C}_c^\infty(B_1) \ : \| \rho\|_{\L^\infty}+ \|\nabla \rho\|_{\L^\infty}< 1\ \} .
\end{equ}
Next, recall the canonical projection 
\begin{equ}\label{eq:projection_T^d}
\pi_{d}: \mathbb{R}^d\to \mathbb{T}^d\ 
\end{equ}
and the associated pullback $\pi_d^*$ of functions, resp.\ distributions.
We define $C^\gamma(\mathbb{T}^d)$ as the closure\footnote{
We choose this slightly non-standard convention, since it makes the space separable.
} 
of 
$f\in C^\infty(\mathbb{T}^d)$ under the following norms respectively:
\begin{equ}
\|f\|_{C^\gamma(\mathbb{T}^d)}:=
\begin{cases}
\|\pi_d^* f\|_{C^\gamma (\mathbb{R}^d)}+ \|\pi^*_d f\|_{\L^\infty(\mathbb{R}^d)} & \text{if } \gamma\in (0,1), \\
\|\pi_d^* f\|_{\L^\infty(\mathbb{R}^d)}& \text{if } \gamma=0,\\
\|\pi_d^* f\|_{C^\gamma (\mathbb{R}^d)} & \text{if } \gamma\in (-1,0)\ ,
\end{cases}
\end{equ}
where $\|\, \cdot\, \|_{\L^\infty([-1,1]^d)}$ denotes the usual (essential) supremum-norm.
We shall often freely identify functions, resp.\ distributions on $\mathbb{T}^d$ with periodic functions, resp. distributions by pullback under
\eqref{eq:projection_T^d}.

%

\subsubsection{Function spaces on \TitleEquation{\mathbb{T}^d\times \mathbb{R}}{T^dxR}}
We equip $\mathbb{R}^{d+1}$ with the parabolic scaling $(1,\ldots,1,2)$ in the sense of \cite{Hai14} and write $|z|_\fraks= |x|+ |t|^{1/2} $ for $z=(x,t)\in \mathbb{R}^{d+1}$. 
Similarly to above, we shall often identify functions and distributions on $\mathbb{T}^d\times \mathbb{R}$ with their counterpart on $\mathbb{R}^{d+1}$ by pullback under the projection
\begin{equ}
\pi_{d,1}: \mathbb{R}^{d+1}\to \mathbb{T}^d\times \mathbb{R}\, \qquad
(x,t) \mapsto (\pi_d x, t)\ .
\end{equ}
For $\gamma\in (0,1)$ we define the H\"older semi-norm for subsets $B\subset \mathbb{R}^d\times \mathbb{R}$ as
\begin{equ}
\|f\|_{C^\gamma_\fraks(B)}= \sup_{(x,t), (\zeta, \tau)\in B}  \frac{|f(x,t)- f(\zeta, \tau)| }{(|x-\zeta|+ \sqrt{|t-\tau|})^\gamma} = \sup_{z, z'\in B}  \frac{|f(z)- f(z')| }{|z-z'|_\fraks^\gamma} \ .
\end{equ}
Then, define for $\gamma \geq 0$ the space $C^\gamma(\mathbb{T}^d\times[0,T])$ as the closure of 
$ C^\infty(\mathbb{T}^d\times [0,T])$ under the norms
\begin{equ}
\|f\|_{C_{\fraks}^\gamma(\mathbb{R}^d\times [0,T])}:=
\begin{cases}
\|\pi_{d,1}^* f\|_{\L^\infty(\mathbb{R}^d\times [0,T])}& \text{if } \gamma=0,\\
\|\pi_{d,1}^* f\|_{C_{\fraks}^\gamma (\mathbb{R}^d\times [0,T])}+ \|\pi_{d,1}^*  f\|_{\L^\infty(\mathbb{R}^d\times [0,T])} & \text{if } \gamma\in (0,1),
\end{cases}
\end{equ}
In order to define the corresponding distribution spaces, denote for $\lambda>0$ by
\begin{equ}
\mathcal{S}_{\fraks}^\lambda: \mathbb{R}^{d+1}\to \mathbb{R}^{d+1}, \qquad z=(x,t) \mapsto \mathcal{S}_{\fraks}^\lambda(z)= (\lambda^{-1}x, \lambda^{-2}t)\  
\end{equ}
the parabolic scaling map.
This time, writing for $\rho: \mathbb{R}^{d+1}\to \mathbb{R}$
\begin{equ}
\rho^{\lambda}_z= \lambda^{-(d+2)}  \rho\left( \mathcal{S}_{\fraks}^\lambda(\cdot\ - z) \right) \ ,
\end{equ}
define for $\gamma\in (-1,0)$ and bounded $B\subset \mathbb{R}^{d+1}$ the following semi-norms
\begin{equ}
\| F\|_{C_{\fraks}^\gamma(B)} = \sup_{z\in B} \sup_{\lambda\in (0,1)}    \sup_{\rho\in \mathfrak{B}^{1}(\mathbb{R}^{d+1})} \frac{|\langle F, \rho^{\lambda}_x \rangle|}{\lambda^\gamma} \ .
\end{equ}
Define the space $C_{\fraks}^\gamma(\mathbb{R}^{d+1})\subset \mathcal{D}' (\mathbb{R}^{d+1})$ as the completion of smooth compactly supported functions $C^\infty_c(\mathbb{R}^{d+1})$ with respect to the Fréchet topology induced by the seminorms 
\begin{equ}
\| f\|_{C_{\fraks}^\gamma ((-N,N+1)^{d+1})} \qquad \text{ for } N\in \mathbb{N} 
\end{equ}
and identify $C_{\fraks}^\gamma(\mathbb{T}^{d}\times \mathbb{R})$ as the pushforward under 
$\pi_{d,1}$ of the subspace of $C_{\fraks}^\gamma(\mathbb{R}^{d+1})$ consisting of the periodic (in the first $d$ components) distributions.
\subsubsection{Singular function spaces}
In order define function spaces allowing for singular behaviour near time $0$, 
we first introduce the forward parabolic cylinder 
\begin{equ}\label{eq:backwards parabolic cylinder}
\tilde{Q}_{r}(x,t)=B_r(x)\times (t, t+r^2)\subset \mathbb{R}^{d+1}\;,
\end{equ}
at $(x,t)\in \mathbb{R}^d\times \mathbb{R}$ of radius $r>0$.
For $\gamma\in (0,1]$, $\eta\leq \gamma$ we define the weighted H\"older norm on functions $f: \mathbb{T}^d\times (0,T] \to \mathbb{R}^d$
\begin{equs}
\| f\|_{C^{\gamma,\eta}_T}&=  \sup_{(x,t)\in  [-1,1]^d \times (0,T]} \sup_{(\zeta, \tau)\in \tilde{Q}_{\sqrt{t}}(t,x)\cap  ([-1,1]^d \times (0,T])}  \frac{|\pi_{d,1}^*f(x,t)- \pi_{d,1}^*f(\zeta, \tau)| }{\sqrt{t}^{\eta-\gamma} (|x-\zeta|+ \sqrt{|t-\tau|})^\gamma} 
\\ 
&\qquad +  
\sup_{(x,t)\in [-1,1]\times (0,T]} \frac{ |\pi_{d,1}^* f(x,t)|}{|t|^{\eta/2 \wedge 0}}\;.
\end{equs}

%
%

\subsubsection{Functions of two variables}
Furthermore for $\gamma\in (0,1]$ we introduce the following H\"older semi-norms for subsets $B,B'\subset \mathbb{R}^d\times \mathbb{R}$
\begin{equs}
\|F\|_{C^{0,\gamma}_\fraks(B\times B')}&= \sup_{z\in B} \sup_{z',\bar{z}'\in B'}  \frac{|F(z,z') -F({z}, \bar{z}')|}{ |z'-\bar{z}'|_\fraks^{\gamma} } \ ,\\
\|F\|_{C^{\gamma,0}_\fraks(B\times B')}&=\sup_{z,\bar{z}\in B} \sup_{z'\in B'}  \frac{|F(z,z') -F(\bar{z}, {z}') | }{|z-\bar{z}|_\fraks^\gamma} \ ,
\end{equs}
as well as for $\gamma, \gamma'\in (0,1]$,
\begin{equ}
\|F\|_{C^{\gamma,\gamma'}_\fraks(B\times B')}= \sup_{z,\bar{z}\in B} \sup_{z',\bar{z}'\in B'}  \frac{|F(z,z') -F(\bar{z}, {z}') -F({z}, \bar{z}')+ F(\bar{z}, \bar{z}')| }{|z-\bar{z}|_\fraks^\gamma \cdot |z'-\bar{z}'|_\fraks^{\gamma'} } \ .
\end{equ}

\subsection{Non-translation invariant heat kernels}
 Throughout this article we make the following assumption.
\begin{assumption}\label{ass:kernel}
$A: \mathbb{R}^{d+1}\to \mathbb{R}^{d\times d}$ is uniformly elliptic, $\mathbb{Z}^{d+1}$-periodic  and (space-time) $\theta$-H\"older continuous for some $\theta\in (0,1)$, i.e.
\begin{equ}
\sup_{|z-z'|_\fraks\leq 1} \frac{ |A(z)-A(z')|}{|z-z'|_\fraks^\theta}<+\infty\ .
\end{equ}
\end{assumption} 
We shall throughout use the notation $a\lesssim b$ to mean that there exists a constant $C>0$ such that $a\leq Cb$ and whenever applicable further explain what that constant depends on. Throughout this article, constants will depend on $A$ without further mention.

Next we recall some properties of the fundamental solution $\Gamma$ of the differential operator 
\begin{equation}\label{eq:dif_operator}
\partial_t - \nabla \cdot A(z)\nabla \ ,
\end{equation} 
c.f.\ \cite{Fri08}, \cite{Sin23}. We shall write $a^{i,j}$ for the entries of the inverse matrix $A^{-1}$ of $A$ and define
\begin{equ}[e:deftheta]
\vartheta^{z}(\zeta)= \sum_{i,j} a^{i,j}(z) \zeta_i \zeta _j, \qquad
 w^{z}(\zeta, \tau)= \frac{\mathbf{1}_{\{\tau > 0 \}}}{\tau^{d/2}} \exp \left(\frac{\vartheta^{z}(\zeta)}{4\tau}\right)\ .
\end{equ}
The fundamental solution of the differential operator with coefficients ``frozen'' at $z=(\zeta, \tau)$ is given by 
\begin{equation}\label{eq:exlicit Z_0}
Z(x,t; \zeta,\tau): =  C(\zeta,\tau) w^{(\zeta,\tau)}(x-\zeta,t-\tau)\ ,
\end{equation}
where $C(z)= (4\pi)^{-d/2} \det(A^s(z))^{-1/2}$.
We also define 
\begin{equ}\label{eq:explicit bar Z}
\bar{Z}(x,t; \zeta,\tau): =  C(x,t) w^{(x,t)}(x-\zeta,t-\tau)\ ,
\end{equ}
and note by direct computation that
\begin{equ}\label{eq:second_ineq_prop}
|\bar{Z} (x,t;\zeta, \tau-\delta)-\bar{Z}(x,t;\zeta, \tau)| \lesssim 
\frac{\delta^{\theta/2}
 }{ (t-\tau)^{d/2+\theta/2}}  \exp\left( -\frac{\mu |x-\zeta|^2}{t-\tau} \right)
\end{equ}
uniformly over $t>\tau>0$, $x,\zeta\in \mathbb{R}^d$, $\delta\in [0,1]$.

The following is essentially known and summarises some properties of the fundamental solution.
\begin{prop}\label{prop:heat_kernel}
Under Assumption~\ref{ass:kernel} there exists $\mu>0$ such that the (unique) fundamental solution $\Gamma$ of the differential operator \eqref{eq:dif_operator}
satisfies
\begin{equ}
|\Gamma(x,t; \zeta, \tau)-Z(x,t;\zeta, \tau)|+|\Gamma(x,t; \zeta, \tau)-\bar Z (x,t;\zeta, \tau)| \lesssim (t-\tau)^{(\theta-d)/2}  \exp\left( -\frac{\mu |x-\zeta|^2}{t-\tau} \right)\ ,
\end{equ}
and 
\begin{align*}
&|\Gamma(x,t; \zeta, \tau-\delta)-\bar{Z} (x,t;\zeta, \tau-\delta)-\Gamma (x,t; \zeta, \tau)+\bar{Z}(x,t;\zeta, \tau)| \\
&\qquad \lesssim \left(\left(\frac{
\delta^{\theta/2} }{ (t-\tau)^{d/2+\theta/2} } +  \frac{
\delta^{\theta/2} }{ (t-\tau)^{d/2} }   \right) \wedge \frac{
1} { (t-\tau)^{d/2-\theta/2} }  \right)
\exp\left( -\frac{\mu |x-\zeta|^2}{t-\tau} \right)
\end{align*}
uniformly over $t>\tau>0$, $x,\zeta\in \mathbb{R}^d$, $\delta\in [0,1]$.
\end{prop}
\begin{remark}
Throughout this article, $\mu$ will denote some strictly positive constant (depending on $A$), which is allowed to change from line to line.
\end{remark}

\begin{proof}
Recall that by the construction of the fundamental solution by Levi's parametrix method, c.f. \cite[Ch.~1 \& Ch.~9]{Fri08}, we can write $Z_{\geq 1}:= \Gamma-Z$ as
\begin{equ}
Z_{\geq 1} = \sum_{\nu=1}^\infty Z_{\nu} \ ,
\end{equ}
where each summand is explicitly given in \cite[Eq.~2.5]{Sin23}.
Furthermore, it follows from the discussion around \cite[Ch.~1, Sec.~4, Eq.~4.14]{Fri08} that there exist
 $H,\mu>0$ such that
\begin{equ}
|Z_{\nu}(x,t; \zeta,\tau)|
\lesssim
 \frac{H^{\nu}}{\Gamma(\theta \nu)} 
\frac{
1}{ (t-\tau)^{d/2-\theta\nu/2}}  \exp\left( -\frac{\mu |x-\zeta|^2}{t-\tau} \right)
\end{equ}
and
\begin{equ}
| Z_{\nu}(x,t; \zeta,\tau-\delta)-Z_{\nu}(x,t; \zeta,\tau)|
\lesssim
 \frac{H^{\nu}}{\Gamma(\theta \nu)} 
\frac{
\delta^{\theta/2} }{ (t-\tau)^{d/2-\theta/2(\nu-1)}}  \exp\left( -\frac{\mu |x-\zeta|^2}{t-\tau} \right)\ ,
\end{equ}
where the $\Gamma(\cdot)$ in the denominator denotes the Gamma function.
We conclude that 
\begin{equ}\label{eq:A}
|Z_{\geq 1}(x,t; \zeta,\tau)|
\lesssim
\frac{
1}{ (t-\tau)^{d/2-\theta/2}}  \exp\left( -\frac{\mu |x-\zeta|^2}{t-\tau} \right) 
\end{equ}
and
\begin{equ}\label{eq:A'}
|Z_{\geq 1}(x,t; \zeta, \tau-\delta)-Z_{\geq 1}(x,t; \zeta, \tau)|\lesssim \frac{
\delta^{\theta/2} }{ (t-\tau)^{d/2}}  \exp\left( -\frac{\mu |x-\zeta|^2}{t-\tau} \right)\ .
\end{equ}
A direct computation shows that 
\begin{equ}\label{eq:B}
|Z(x,t; \zeta,\tau)-\bar Z (x,t; \zeta,\tau)|\lesssim (t-\tau)^{(\theta-d)/2}  \exp\left( -\frac{\mu |x-\zeta|^2}{t-\tau} \right) \ .
\end{equ}
Thus \eqref{eq:A} together with \eqref{eq:B} imply the first inequality of the proposition.
 Next note that \eqref{eq:B} together with \eqref{eq:second_ineq_prop} and a similar bound on $Z$ imply that
 $\hat{Z} = \bar{Z} - {Z}$ satisfies
\begin{align*}
&| \hat {Z} (x,t; \zeta,\tau-\delta) - \hat {Z} (x,t; \zeta,\tau) | \\
&\qquad \lesssim \left(\left(\frac{
\delta^{\theta/2} }{ (t-\tau)^{d/2+\theta/2} } +  \frac{
\delta^{\theta/2} }{ (t-\tau)^{d/2} }   \right) \wedge \frac{
1} { (t-\tau)^{d/2-\theta/2} }  \right)
\exp\left( -\frac{\mu |x-\zeta|^2}{t-\tau} \right)\ .
\end{align*}
which combined with \eqref{eq:A'} implies the last inequality of the proposition.
\end{proof}

We shall denote by $\Gamma_\eps$ the heat kernel of the differential operator
\begin{equ}
\partial_t +\mathcal{L}_\eps\ \qquad \text{with }\qquad  \mathcal{L}_\eps= - \nabla \cdot (A(x/\eps, t/\eps^2)  \nabla)  \ .
\end{equ}
We observe that for $\eps>0$ and $f\in \mathcal{C}_c^\infty (\mathbb{R}^{d+1})$, if $u$ is a solution to
\begin{equ}
(\partial_t +\mathcal{L}_1)u = f\circ \mathcal{S}_\fraks^{1/\eps}  \qquad \text{ on } \qquad \mathbb{R}^{d+1} \ ,
\end{equ}
then $u_\eps= \eps^2 u \circ \mathcal{S}^\eps_\fraks$ satisfies 
\begin{equ}\label{eq:scaling_property}
(\partial_t +\mathcal{L}_\eps) u_\eps = f \ .
\end{equ} 
We can conclude the following scaling property for the heat kernel.
%

\begin{lemma}\label{lem:formula for Gamma_epsilon}
It holds that
$
\Gamma_\eps (z,\bar{z})= \frac{1}{\eps^d} \Gamma_1 (\mathcal{S}^{\eps}_{\fraks} (z), \mathcal{S}_\fraks^\eps\bar{z})$.\qed
\end{lemma}
%
We define $Z^\eps$ and $\bar{Z}^\eps$ exactly as in \eqref{eq:exlicit Z_0} and \eqref{eq:explicit bar Z} but with $A(x,t)$ replaced by 
$A(x/\eps, t/\eps^2)$. One directly checks that 
\begin{equ}\label{eq:bar Z small scale increment}
|\bar{Z}^{\eps} (x,t;\zeta, \tau-\delta^2)-\bar{Z}^{\eps}(x,t;\zeta, \tau)| \lesssim 
\frac{\delta^{\theta}}{ (t-\tau)^{d/2+\theta/2}} \exp\left( -\frac{\mu |x-\zeta|^2}{t-\tau} \right)\ .
\end{equ}
\begin{corollary}\label{cor:small scale bound}
There exists $\mu>0$ such that
\begin{equation}\label{small scale bound}
|\Gamma_\eps(x,t; \zeta,\tau)- \bar{Z}^\eps (x,t; \zeta,\tau) |\lesssim 
\frac{\eps^{-\theta}}{ (t-\tau)^{d/2-\theta/2} } 
 \exp\left( -\frac{\mu |x-\zeta|^2}{t-\tau} \right)\ ,
\end{equation}
and
\begin{align*}
&|\Gamma_\eps(x,t; \zeta, \tau-\delta^2)-\bar{Z}^\eps (x,t;\zeta, \tau-\delta^2)-\Gamma_\eps (x,t; \zeta, \tau)+\bar{Z}^\eps(x,t;\zeta, \tau)| \\
&\qquad \lesssim \mathbf{1}_{\delta\leq \eps}\cdot\left(\left(\frac{
\delta^{\theta}\eps^\theta }{ (t-\tau)^{d/2+\theta/2} } +  \frac{
\delta^{\theta} }{ (t-\tau)^{d/2} }   \right) \wedge \frac{
\eps^{-\theta}} { (t-\tau)^{d/2-\theta/2} }  \right)
\exp\left( -\frac{\mu |x-\zeta|^2}{t-\tau} \right)\\
&\qquad+
\frac{\mathbf{1}_{\delta>\eps}
\delta^{\theta} }{ (t-\tau)^{d/2} }
\exp\left( -\frac{\mu |x-\zeta|^2}{t-\tau} \right)
\end{align*}
uniformly over $t>\tau>0$, $x,\zeta\in \mathbb{R}^d$ and $\delta,\eps\in (0, 1]$. 
\end{corollary}
\begin{proof}
From Lemma~\ref{lem:formula for Gamma_epsilon} and the fact that  $\eps^{-d}\bar{Z}^1(\mathcal{S}_{{\fraks}}^\eps z; \mathcal{S}_{{\fraks}}^\eps \bar z )= \bar{Z}^\eps(z; \bar z )$ one sees that 
\begin{align*}
|\Gamma^\eps(z;\bar{z})- \bar{Z}^\eps(z;\bar{z})|= \eps^{-d} |\Gamma^1(z;\bar{z})- \bar{Z}^1(z;\bar{z})|
\end{align*}
which, combined with
Proposition~\ref{prop:heat_kernel}, shows \eqref{small scale bound}.
For the latter inequality, we treat the regimes $\delta\leq \eps$ and $\delta>\eps$ separately. The former follows by the same scaling argument as above, while for $\delta>\eps$ the bound follows by the triangle inequality from \eqref{eq:bar Z small scale increment} and Proposition~\ref{prop:uniform_hölder_bound} below.
\end{proof}

\subsection{Periodic Homogenisation}
We now recall some known results from periodic homogenisation, always under Assumption~\ref{ass:kernel}.
One defines for $j=1,\ldots,d$ the \textit{corrector} $\Phi_j: \mathbb{T}^d\times \mathbb{R}\to \mathbb{R}$ as the unique (weak) $1$-periodic solution to the cell problem
\begin{equation}\label{eq:correctors}
(\partial_t+\mathcal{L}_1)\Phi_j= \nabla \cdot (A(x,t) e_j)\; , \qquad 
\int_{\mathbb{T}^d\times [0,1]} \Phi_j(x,t) \,dx\, dt =0\;,
\end{equation}
where $e_j$ denotes the $j$th basis vector of $\mathbb{R}^d$, c.f.\ \cite[Sec.~2]{GS17}.
Then, the \textit{homogenised matrix} defined as 
\begin{equation}\label{eq:hom_matrix}
\bar A= \big(\bar a_{i,j}\big)_{i,j=1}^d, \qquad \bar{a}_{i,j}:= \int_{ [0,1]^{d+1}} \left(a_{i,j}(x,t) +\sum_{k=1}^d a_{i,k}(x,t) \partial_k \Phi_j(x,t) \right)   dx\,dt \;,
\end{equation}
is strictly positive definite.
Writing $\tilde{A}:  \mathbb{R}^{d+1}\to \mathbb{R}^{d\times d}$ for the matrix with entries $\tilde{a}_{i,j}(y,s)={a}_{j,i}(y,-s)$, we set 
$\tilde{\mathcal{L}}_\eps= -\nabla \cdot (\tilde A(x/\eps, t/\eps^2) \nabla)$. One finds that the fundamental solution $\tilde{\Gamma}_\eps(x,t,\zeta,\tau)$ of $\partial_t +\tilde{\mathcal{L}}_\eps$ satisfies for $\tau<t$ (and thus $-t<-\tau$)
\begin{equation}\label{eq:adjoint}
\tilde{\Gamma}_\eps (x,t;\zeta,\tau)= \Gamma_\eps (\zeta, -\tau;x, -t)\ .
\end{equation}
We shall denote by $\bar \Gamma$ the fundamental solution of the homogenised operator $\partial_t- \nabla\cdot \bar{A} \nabla$.
For $j=1,\ldots,d$ we then denote by 
\begin{equation}\label{eq:tilde correctors}
\tilde{\Phi}_j: \mathbb{T}^d\times \mathbb{R}\to \mathbb{R}
\end{equation}
the correctors associated to the homogenisation problem associated to $\tilde{A}$.
For $j,i=1,\ldots,d$, let
\begin{equ}
{b}_{i,j}:= a_{i,j} +\sum_{k=1}^d a_{i,k}(x,t) \partial_k \Phi_j(x,t) - \bar{a}_{i,j} \ .
\end{equ}
as well as $b_{d+1,j}= - \Phi_j\ .$
The functions $\{\Psi_{k,i,j}\}_{k,i=1; j=1}^{d+1;d},$ characterised by the next Lemma \cite[Lem.~2.1]{GS17}, are called \textit{dual correctors}.
\begin{lemma}
There exist continuous periodic functions $\Psi_{k,i,j}: \mathbb{R}^{d+1}\to \mathbb{R}$ satisfying
\begin{equ}
b_{i,j}=   \partial_{t} \Psi_{{d+1},i,j} + \sum_{k=1}^{d} \partial_k \Psi_{k,i,j}, \qquad \Psi_{k,i,j}=-\Psi_{k,j,i}\ .
\end{equ}
\end{lemma}
We write 
\begin{equ}\label{eq:parabolic cylinder}
Q_{r}(x,t)= B_r(x) \times (t-r^2, t)\subset \mathbb{R}^{d+1}\;,
\end{equ}
for the parabolic cylinder at $(x,t)\in \mathbb{R}^d\times \mathbb{R}$ of radius $r>0$.
The following is \cite[Thm~1.1]{GS15}.
\begin{theorem}\label{thm:technical}
Let $p>d+2$ and $\alpha=1-\frac{d+2}{p}$ and $R>0$. There exists a constant $C=C(R,d,p)$ such that for all $\eps\in (0,1]$ and $f=(f_1,\ldots,f_d)\in L^p(Q_{2r}(x_0,t_0))$
any solution
 $u_\eps$ to $(\partial_t+\mathcal{L}_\eps)u_\eps= \nabla\cdot f$ on $Q_{2r}(x_0,t_0)$ satsfies the bound
\begin{equ}
\|u_\eps\|_{C_\fraks^{\alpha}(Q_{r}(x_0,t_0))} \leq Cr^{1-\alpha}\left( \frac{1}{r}  \left(\fint_{Q_{2r}(x_0,t_0)} |u_\eps|^2\right)^{1/2}  + \left(\fint_{Q_{2r}(x_0,t_0)} |f|^p\right)^{1/p} \right)\
\end{equ}
uniformly over $r<R$.
\end{theorem}
The next theorem is well known, see for example \cite[Thm~1.1]{HK04} and \cite[Thm~4.1]{GS15}.
\begin{theorem}\label{thm:uniform}
There exists $\mu>0$ such that 
\begin{equ}
| \Gamma_\eps(x,t;y,s)|\lesssim \frac{1}{|t-s|^{\frac{d}{2}}} \exp\left(-\frac{\mu|x-y|^2}{t-s} \right)
\end{equ}
uniformly over $\eps\in (0,1]$, $x,y\in \mathbb{R}^d$ and $-\infty <s<t<\infty$.
\end{theorem}
The following is \cite[Thm~1]{GS20}. 
\begin{prop}\label{propGS20}
There exists $\mu>0$ such that 
\begin{equation}\label{eq:GS20thm1}
|(\Gamma_\eps-\bar \Gamma)  (x,t; y, s)|\lesssim \frac{C\eps}{(t-s)^{\frac{d+1}{2}}} \exp\left( \frac{\mu|x-y|}{t-\tau} \right)
\end{equation}
uniformly over $x,y\in \mathbb{R}^d$ and $-\infty <s<t<\infty$.
\end{prop}

\subsection{Main results}\label{subsec:main results}
%
%
Denote by $\xi$ space time white noise on $\mathbb{T}^d\times \mathbb{R}$. 
Fix $\phi\in C^\infty_c(B_{1}(0))$ even such that $\int_\mathbb{R}\phi(t)dt=1$.
For $\eps, \delta>0$ set
\begin{equation}\label{eq:noise}
\xi_{\eps,\delta} (x,t)= \int_{\mathbb{R}^{d+1}}
\frac{1}{\delta^2} \phi\left((t-s)/\delta^2\right)
 \Gamma_{\eps}(x,t, \zeta, t-\delta^2) \xi(d\zeta, ds)\ ,
\end{equation}
where we implicitly identified $\xi$ with its periodic counterpart by pullback.
The following lemma suggests that this regularisation is particularly convenient for studying homogenisation of singular SPDEs. It is a direct corollary of Proposition~\ref{prop:IXi regularity} in Section~\ref{sec:stochastic_estimates}.
\begin{lemma}\label{lem:regularisation_for_homogenisation}
For every $\alpha>\frac{d+2}{2}$ there exists a modification of \eqref{eq:noise} which extends to a continuous map
\begin{equ}{}
[0,1]^2\to C_{\fraks}^{-\alpha}(\mathbb{R}^{d+1}) \ , \qquad (\eps,\delta)\mapsto\xi_{\eps,\delta}\ .
\end{equ}
It has the property that for any $\eps\in [0,1]$ it holds that $\xi_{\eps,0}= \xi$ and that for any $\delta>0$
\begin{equ}\label{eq:homogenised noise}
\xi_{0,\delta} (x,t)= \int_{\mathbb{R}^{d+1}} \phi^\delta(t-s) \bar{\Gamma}(x,t, \zeta, t-\delta^2) \xi(d\zeta, ds)\  .
\end{equ}
\end{lemma}

We define the set
\begin{equ}
\square:= \{ (\eps, \delta) \in (0,1]^2 \ : \ \eps^{-1}\in \mathbb{N} \} \;, 
\end{equ}
in order to state the main results of this article, see Remark~\ref{rem:square}. We will always view $\square$ as a subspace of $[0,1]^2$, so its closure
$\bsquare$ equals $\bsquare = \{ (\eps, \delta) \in [0,1]^2 \ : \ \eps\in \{0\} \cup \mathbb{N}^{-1} \}$.
\begin{theorem}\label{thm:main1}
Let $\xi$ denote space-time white noise on $\mathbb{T}^2\times \mathbb{R}$ and let $u_0\in \mathcal{C}^{\alpha}(\mathbb{T}^2)$ for $\alpha>-1/10$.  Consider for $\delta>0$ the  regularisation $\xi_{\eps,\delta}(x,t)$ as in \eqref{eq:noise}. There exist constants\footnote{
Here we write $f\sim g$ (for functions $f, g: (0,1]^2 \to \mathbb{R}$) to mean that 
$f-g$ extends continuously to $[0,1]^2$.
} 
\begin{equ}\label{eq:assymptotics renormalsation constants}
\alpha_{\eps,\delta}\sim 
\frac{|\log(\delta/\eps)\wedge 0|}{4\pi} 
, \qquad 
\bar{\alpha}_{\eps,\delta}\sim
  \frac{|\log(\delta)| \wedge |\log(\eps)|}{{4\pi} },
\end{equ}
and a bounded family of constants ${c_{\eps,\delta}}$ such that if 
 we denote by $u_{\eps, \delta}$ the solution to
\begin{align*}
\partial_t u_{\eps, \delta} -\mathcal{L}_\eps u_{\eps, \delta}=& -u_{\eps,\delta}^3+ 3 \alpha_{\eps, \delta}  D(x/\eps,t/\eps^2)u_{\eps,\delta} + 3\bar{\alpha}_{\eps,\delta} \det(\bar {A})^{-1/2}  u_{\eps,\delta} + 3c_{\eps,\delta} u_{\eps,\delta} +\xi_{\eps, \delta}\  , 
\end{align*}
with initial condition $u_{\eps, \delta}(0)=u_0,$ then, for any $T>0$ the solution map 
\begin{equ}
\square\to \L^{0}\left[C([0,T], \mathcal{D}'(\mathbb{T}^2))\right]\;, \qquad (\eps, \delta) \mapsto u_{\eps, \delta}\;,
\end{equ}
is well defined and has a unique continuous\footnote{Recall that the $\L^0$-topology is characterised by convergence in probability.} extension to $\bsquare$. 
Furthermore, the constants $\alpha_{\eps,\delta}, \bar{\alpha}_{\eps,\delta}, c_{\eps,\delta}$ (depending on $\phi$) can be chosen such that the following hold. 
\begin{enumerate}
\item\label{thm:main_item1} For $\delta>0$ one has $ \lim_{\eps\to 0} \alpha_{\eps,\delta}=0$, $\lim_{\eps\to 0} c_{\eps,\delta}=0$ and the limit $ \bar{\alpha}_{0,\delta}:=\lim_{\eps\to 0} \bar{\alpha}_{\eps,\delta}$ exists. In particular,
$u_{0, \delta}$ agrees with the classical solution to
\begin{equ}
\partial_t u_{0, \delta} -\nabla\cdot \bar{A}\nabla u_{0, \delta}= -u_{\eps,\delta}^3+ 3\bar{\alpha}_{0,\delta}\det(\bar {A})^{-1/2}   u_{0,\delta} +\xi_{0, \delta} \ .
\end{equ}
\item\label{thm:main_item2} For $\eps> 0$, the limits $\bar{\alpha}_{\eps,0}=\lim_{\delta\to 0} \bar{\alpha}_{\eps,\delta}$ and  $c_{\eps, 0}=\lim_{\delta\to 0} c_{\eps,\delta}$ exist. Furthermore there exists a sequence of constants 
$\hat{\alpha}_\eps\sim \frac{1}{4\pi} |\log(\eps)| , $
such that $u_{\eps,0}$ agrees with the\footnote{Assuming we choose the same way to affinely parametrise the solution family, e.g.\ by choosing the same cutoff function $\kappa(t)$, see Section~\ref{subsec:modification for spdes}.
In that case,
\begin{equ}
\hat{\alpha}_\eps=\frac{1}{8}  \int_{\mathbb{R}} \frac{\kappa(\tau)^2(1- \kappa(\tau/\eps^2)^2) }{\tau} d\tau \ .
\end{equ}
}
 solution
to 
\begin{equs}
\partial_t u_{\eps,0}  -\mathcal{L}_\eps u_{\eps,0}=& -u_{\eps,0}^3+ \infty\cdot D(x/\eps,t/\eps^2)u_{\eps,0} + 3\left(\bar{\alpha}_{\eps,0} \det(\bar{A})^{-1/2} + c_{\eps,0}\right) \\
&- 3\hat{\alpha}_\eps D(x/\eps,t/\eps^2)  u_{\eps,0}  + \xi
\end{equs}
constructed in \cite{Sin23}. 
\item\label{thm:main_item3} For all $\eps\in[0,1]$, $u_{\eps,0}$ does not depend on the choice of $\phi$. 
\end{enumerate}
%
%
\end{theorem}
\begin{remark}\label{rem:square}
Note that the restriction $(\eps,\delta)\in \square$ instead of $(\eps,\delta)\in (0,1]^2$ is only so that the differential operator  $\nabla \cdot A(x/\eps, t/\eps^2)\nabla$ can be pushed forward to the torus. We could 
just as well have formulated the results on the full plane but with periodic noise instead, in which case the statement 
holds for $\square$ replaced by $(0,1]^2$. The same remark applies to Theorem~\ref{thm:main_translation invariant}, Proposition~\ref{prop:main_flat} and Theorem~\ref{thm:main_non} below.
\end{remark}
\begin{remark}
The constants $c_{\eps,\delta}\in \mathbb{R}$ have the property that 
 in general
\begin{equ}
\lim_{\delta\downarrow 0} \lim_{\eps\downarrow 0}  c^{(1)}_{\eps,\delta}= 0\neq  \lim_{\eps\downarrow 0}\lim_{\delta\downarrow 0}c^{(1)}_{\eps,\delta}\  ,
\end{equ}
see \eqref{eq:constant_finite shift_c} for an expression of the limit on the right hand side.
\end{remark}
\begin{remark}
It can be seen that Theorem~\ref{thm:main1} as well as its proof extend to arbitrary polynomial non-linearities of odd degree where the highest degree part has a negative coefficient. We focus on the cubic case in order to streamline exposition, the generalisation follows along standard arguments, c.f. \cite{TW18}. 
\end{remark}

\subsubsection{Translation invariant regularisation}
In this section we consider homogenisation of the $\phi^4_2$ equation for usual translation invariant regularisations of the noise. Let us define for $C>0$ the following set
\begin{equ}
\triangle^{<}_{C}:=\{ (\eps,\delta)\in (0,1]^2  \ : \ \eps<C\delta \}  \ .
\end{equ}
\begin{theorem}\label{thm:main_translation invariant}
Let $\xi$ denote space-time white noise on $\mathbb{T}^2\times \mathbb{R}$, let $\phi\in C_c^\infty(B_1)$ be even, non-negative with $\int_{\mathbb{R}^{d+1}} \phi=1$, and let $u_0\in \mathcal{C}^{\alpha}(\mathbb{T}^2)$ for $\alpha>-1/10$.  Consider for $\delta>0$ the regularisation $\xi^\flat_{\delta}(z)=\xi(\phi_z^\delta)$.
There exist constants
\begin{equ}\label{eq:assymptotics renormalsation constants_flat}
\bar{\alpha}^\flat_{\eps,\delta}\sim 
  \frac{|\log(\delta)| \wedge |\log(\eps)|}{{4\pi} },
  \end{equ}
and ${c^\flat_{\eps,\delta}}$ bounded on $\triangle^{<}_{C}$ for any $C>0$
such that if 
we denote by $u^\flat_{\eps, \delta}$ the solution to
\begin{align*}
\partial_t u^{\flat}_{\eps, \delta} -\mathcal{L}_\eps u^{\flat}_{\eps, \delta}=& -(u^{\flat}_{\eps,\delta})^3+ 3\bar{\alpha}^\flat_{\eps,\delta} \det(\bar {A})^{-1/2}  u^{\flat}_{\eps,\delta} + 3c^\flat_{\eps,\delta} u^{\flat}_{\eps,\delta} +\xi^{\flat}_{ \delta}\  , 
\end{align*}
with initial condition $u_{\eps, \delta}^\flat (0)=u_0,$ then for every $C>0$ and $T>0$ the solution map
\begin{equ}
\triangle^{<}_{C}\cap \square
\to \L^{0}\left[C([0,T], \mathcal{D}'(\mathbb{T}^2))\right], \qquad (\eps, \delta) \mapsto u_{\eps, \delta}
\end{equ}
is well defined and has a unique continuous extension to $\triangle^{<}_{C}\cap \bsquare$. 
Furthermore, the constants $\bar{\alpha}_{\eps,\delta}^\flat, c^\flat_{\eps,\delta}$ (depending on $\phi$) can be chosen such that the following hold. 
\begin{enumerate}
\item\label{thm:main_flat_item1} For $\delta>0$ it holds that $\lim_{\eps\to 0} c^\flat_{\eps,\delta}=0$ and the limit $ \bar{\alpha}^\flat_{0,\delta}:=\lim_{\eps\to 0} \bar{\alpha}^\flat_{\eps,\delta}$ exists. In particular,
$u^\flat_{0, \delta}$ agrees with the classical solution to
\begin{equ}
\partial_t u^\flat_{0, \delta} -\nabla\cdot \bar{A}\nabla u^\flat_{0, \delta}= -(u^\flat_{0,\delta})^3+ 3\bar{\alpha}_{0,\delta}^\flat\det(\bar {A})^{-1/2}   u^\flat_{0,\delta} +\xi^\flat_{0, \delta} \ .
\end{equ}
\item\label{thm:main_flat_item2} The process $u^\flat_{0,0}$ does not depend on the choice of $\phi$ and agrees with $u_{0,0}$ in Theorem~\ref{thm:main1}.
\end{enumerate}
\end{theorem}
\begin{remark}
See Item~\ref{prop:main_flatflat_item1} of Proposition~\ref{prop:main_flat} and Remark~\ref{rem:function Dflatflat} for a formula for the constant $c^\flat_{\eps,\delta}$.
\end{remark}

The next proposition shows that the above theorem is sharp in the following sense.
\begin{prop}\label{prop:main_flat}
Let $\xi$, $\phi$, $\xi^\flat_{\delta}$, $u_0\in \mathcal{C}^{\alpha}(\mathbb{T}^2)$ and $\bar{\alpha}^{\flat}_{\eps,\delta}$ be as in Theorem~\ref{thm:main_translation invariant}. There exist
constants ${c^{\flat \flat}_{\eps,\delta}}$ bounded on $(0,1]^2$ and functions $D_{\lambda}^{\flat\flat} : \mathbb{R}^{3}\to \mathbb{R}$ which are $\mathbb{Z}^3$ periodic for each $\lambda\in (0,\infty)$ and  uniformly in $\lambda\in (0,1)$ satisfy
\begin{equ}\label{eq:bound on Dflatflat}
 1+ |\log (\lambda)| \lesssim D_{\lambda}^{\flat\flat} (z) \lesssim 1+|\log (\lambda)| \ ,
 \end{equ} 
 such that if 
 we denote by $u^{\flat\flat}_{\eps, \delta}$ the solution to
\begin{align*}
\partial_t u^{\flat\flat}_{\eps, \delta} -\mathcal{L}_\eps u^{\flat\flat}_{\eps, \delta}=& -(u^{\flat\flat}_{\eps,\delta})^3+ 3  D^{\flat\flat}_{\delta/\eps}(x/\eps,t/\eps^2)u^{\flat\flat}_{\eps,\delta} + 3\bar{\alpha}^{\flat}_{\eps,\delta} \det(\bar {A})^{-1/2}  u^{\flat\flat}_{\eps,\delta} + 3c^{\flat\flat}_{\eps,\delta} u^{\flat\flat}_{\eps,\delta} +\xi^{\flat}_{ \delta}\  , 
\end{align*}
with initial condition $u_{\eps, \delta}^{\flat \flat}(0)=u_0,$ then, for any $T>0$ the solution map
\begin{equ}
\square\to \L^{0}\left[C([0,T], \mathcal{D}'(\mathbb{T}^d))\right], \qquad (\eps, \delta) \mapsto u^{\flat\flat}_{\eps, \delta}
\end{equ}
is well defined and has a unique continuous extension to $\bsquare$. 
Furthermore, the constants $c^{\flat\flat}_{\eps,\delta}$ and the functions $D^{\flat \flat}_\lambda$ can be chosen (depending on $\phi$) such that the following hold. 
\begin{enumerate}
\item\label{prop:main_flatflat_item1} The constants $c^{\flat}_{\eps,\delta}$ from Theorem~\ref{thm:main_translation invariant} can be written as $c^{\flat}_{\eps,\delta}=c^{\flat\flat}_{\eps,\delta} + \int_{[0,1]^{3}} D^{\flat \flat}_{\delta/\eps}$.
\item\label{prop:main_flatflat_item2} For all $\delta>0$, $\lim_{\eps\to 0} D^{\flat \flat}_{\delta/\eps}=0$ and $\lim_{\eps\to 0}  c^{\flat \flat}_{\eps,\delta}=0$.
\item\label{prop:main_flatflat_item3} For all $\delta\in [0,1]$ 
the process $u^{\flat\flat}_{0,\delta}$ agrees with $u^{\flat}_{0,\delta}$ in Theorem~\ref{thm:main_translation invariant}.
\item\label{prop:main_flatflat_item4}
For each $\eps>0$ the limit  $\lim_{\delta\to 0} c^{\flat \flat}_{\eps,\delta}$ exists and agrees with the limit $\lim_{\delta\to 0} c_{\eps,\delta}$ from Theorem~\ref{thm:main1} and, furthermore , $u^{\flat \flat}_{\eps,0}$ agrees with the solution $u_{\eps,0}$ in Theorem~\ref{thm:main1}.
\end{enumerate}
\end{prop}
\begin{remark}\label{rem:function Dflatflat}
Equation~\eqref{eq:formula_Dflatflat} provides an expression for $D^{\flat \flat}_{\lambda}$. 
\end{remark}

\begin{remark}
One can check that 
$\hat{D}^{\phi}(z)= \lim_{\lambda \to 0} D^{\flat \flat}_{\lambda}(z)-\frac{|\log(\lambda)|}{4\pi}D(z)$ 
defines a continuous (but $\phi$-dependent) function. This is a rather generic feature of logarithmic
divergences which tend to come with a regularisation independent prefactor and a 
regularisation dependent part of order one. 
Thus, in the setting of Proposition~\ref{prop:main_flat} it follows by a virtually identical proof that
for $\hat{c}^{\flat\flat}_{\epsilon,\delta}={c}^{\flat\flat}_{\eps,\delta}+ \int_{[0,1]^3} D^{\flat \flat}_{\delta/\eps}-\frac{|\log(\delta/\eps)\vee 0|}{4\pi}\int_{[0,1]^3} D(z)$
the solution map $\square\ni(\eps, \delta) \mapsto \hat{u}^{\flat\flat}_{\eps, \delta} \in \L^{0}\big[C([0,T], \mathcal{D}'(\mathbb{T}^d))\big]$ to
\begin{align*}
\partial_t \hat{u}^{\flat\flat}_{\eps, \delta} -\mathcal{L}_\eps \hat{u}^{\flat\flat}_{\eps, \delta}=& -(\hat{u}^{\flat\flat}_{\eps,\delta})^3+ 3  \frac{|\log(\delta/\eps) \vee 0|}{4\pi}D(x/\eps,t/\eps^2) \hat{u}^{\flat\flat}_{\eps,\delta} \\
&+ 3(\bar{\alpha}^{\flat}_{\eps,\delta} \det(\bar {A})^{-1/2} +\hat{c}^{\flat\flat}_{\eps,\delta})  \hat{u}^{\flat\flat}_{\eps,\delta} +\xi^{\flat}_{\eps, \delta}\  , 
\end{align*}
also extends continuously to $\bsquare$ and it holds that $\hat{u}^{\flat\flat}_{0,\delta}={u}^{\flat\flat}_{0,\delta}$ for all $\delta\in [0,1]$.
But for any $\eps>0$, in general, $\hat{u}^{\flat\flat}_{\eps,0}$ \textit{does} depend on $\phi$ and \textit{does not} fall in the class of solutions exhibited in \cite{Sin23}.
\end{remark}

\subsubsection{Non-translation invariant regularisations}\label{sec:main_result_non}
We first recall the general class non-translation invariant regularisations of \cite[Sec.~3.2]{Sin23} and then establish the analogue result to Theorem~\ref{thm:main_translation invariant}. Let $\rho\in \cC_c^\infty(\R_+)$ non-negative with support in $[0,1)$, all odd derivatives vanishing at the origin, and such that $\int_{\mathbb{R}^d} \rho (|x|^2) dx=1$.
Recalling the definition of $\vartheta^z$ in \eqref{e:deftheta}, let
\begin{equ}
\rho^{(z,\eps, \delta)}(x)=\frac{1}{\delta^2 \det(A(\mathcal{S}^\eps z))^{1/2}} \rho \left(\frac{\vartheta^{\mathcal{S}^\eps z} (x)}{\delta^2}\right)\ .
\end{equ}
We set for $\phi\in  \cC_c^\infty((-1,1))$ even such that $\int \phi(t)\,dt=1$
\begin{equ}
\varrho^{(\eps,\delta; z)} (x,t; \zeta, \tau)  = \frac{1}{\delta^2}\phi\left((t-\tau)/\delta^2\right) \rho^{(z,\eps,\delta)} (x-\zeta)  \ ,
\end{equ}
and $\varrho^{\eps,\delta} (x,t; \zeta, \tau)=\varrho^{(\eps,\delta; (x,t))} (x,t; \zeta, \tau)$.
For later use we also introduce $\bar{\varrho}^\delta (x,t; \zeta, \tau)$ to be defined as  $\varrho^{\eps,\delta} (x,t; \zeta, \tau)$ but with $A$ replaced by the identity matrix.
Let 
\begin{equ}\label{eq:gen_regularisation}
\xi_{\eps,\delta}^{\sharp}(z) = \int \varrho^{\eps,\delta} (z; z') \xi(z')\, dz'\ . 
\end{equ}
Define for $C>0$ the set
$ 
\triangle^{>}_{C}:=\{ (\eps,\delta)\in (0,1]^2  \ : \ \eps>C\delta \}
$. We then have the following result.

\begin{theorem}\label{thm:main_non}
Let $\xi$ denote space-time white noise on $\mathbb{R}\times \mathbb{T}^2$ and let $u_0\in \mathcal{C}^{\alpha}(\mathbb{T}^d)$ for $\alpha>-\frac{1}{10}$.  Denote for $\eps, \delta>0$ by $\xi_{\eps,\delta}^{\sharp}$ its regularisation as in \eqref{eq:gen_regularisation}. Then, there exist sequences of constants 
\begin{equ}\label{thm:main_sharp_constant asympt}
\alpha^{\sharp}_{\delta,\eps}\sim
\frac{|\log(\delta/\eps)\wedge 0|}{4\pi} \ ,
\end{equ}
and $c^{\sharp}_{\delta,\eps}$ bounded on $\triangle^{>}_{C}$ for any $C>0$
 such that if
 we denote by $u^{\sharp}_{\eps, \delta}$ the solution to
\begin{equation}
\partial_t u^{\sharp}_{\eps, \delta} -\mathcal{L}_\eps u^{\sharp}_{\eps, \delta}= -(u^{\sharp}_{\eps,\delta})^3+ 3 \alpha^{\sharp}_{\delta, \eps} D(x/\eps,t/\eps^2)u^{\sharp}_{\eps,\delta} + {c}^{\sharp}_{\eps,\delta}  u_{\eps,\delta}+\xi^{\sharp}_{\eps, \delta} \ ,
\end{equation}
with initial condition $u^\sharp_{\eps, \delta} (0)=u_0,$ then for every $C>0$ and $T>0$ the solution map
\begin{equ}
\triangle^{>}_{C}\cap \square
\to \L^{0}\left[C([0,T], \mathcal{D}'(\mathbb{T}^2))\right], \qquad (\eps, \delta) \mapsto u_{\eps, \delta}
\end{equ}
is well defined and has a unique continuous extension to $\triangle^{>}_{C}\cap \bsquare$.

Furthermore, the constants $\bar \alpha^{\sharp}_{\eps,\delta}, c_{\eps,\delta}^\sharp$ can be chosen (depending on $\phi$ and $\rho$) such that the following hold:
\begin{itemize}
\item\label{thm:main_sharp_item1} For $\eps>0$, the limit $\lim_{\delta\to 0} c_{\eps,\delta}^\sharp$ agrees with $\lim_{\delta\to 0} c_{\eps,\delta}$ from Theorem~\ref{thm:main1}.
\item\label{thm:main_sharp_item2} For all $\eps\in \mathbb{N}^{-1}\cup\{0\}$ the process $u^{\sharp}_{\eps,0}$ agrees with $u_{\eps,0}$ in Theorem~\ref{thm:main1} and does not depend on the choice of $\phi$ and $\rho$.
\end{itemize}
%
%
\end{theorem}

\section{Homogenisation Estimates for Kernels}\label{sec:homogenisation estimates}
In this section we establish further kernel estimates. Recall the definition of the backwards parabolic cylinder $Q_{r}(x,t)$ in \eqref{eq:parabolic cylinder} and the forward parabolic cylinder $\tilde Q_{r}(x,t)$ in \eqref{eq:backwards parabolic cylinder}.
We have the following uniform H\"older bound on heat kernels.
\begin{prop}\label{prop:uniform_hölder_bound}
 For $\alpha\in (0,1)$, there exists $\mu>0$ such that
\begin{equ}
\|\Gamma_{\eps}(\ \cdot\ ; y, s )\|_{C_\fraks^{\alpha}(Q_{\sqrt{ |t-s|}/8}(x,t))}+ \|\Gamma_{\eps}(x,t; \ \cdot\  )\|_{C_\fraks^{\alpha}(\tilde{Q}_{\sqrt{ |t-s|}/8}(y,s))} \lesssim \frac{1}{|t-s|^{\frac{d+\alpha}{2}}} \exp\left(-\frac{\mu|x-y|^2}{t-s} \right)\ .
\end{equ}
uniformly over $x,y\in \mathbb{R}^d$, $-\infty <s<t<\infty$ and $\eps\in (0,1]$.
\end{prop}
\begin{proof}
Note that $u_\eps(\ \cdot\ )=\Gamma_{\eps}(\ \cdot\  ; y, s )$ satisfies the assumptions of Theorem~\ref{thm:technical} with $r=\sqrt{ |t-s|}/8$ and $f=0$.
By Theorem~\ref{thm:uniform}
\begin{align*}
\sup_{z\in Q_{2r}(x,t)} \Gamma_\eps (\ \cdot ,y,s )&\lesssim \sup_{x'\in B_{2r}(x),\  t' \in (t-(2r)^2, t) }  \frac{1}{|t'-s|^{\frac{d}{2}}} \exp\left(-\frac{\mu|x'-y|^2}{t'-s} \right)\\
&\lesssim \frac{1}{|t-s|^{\frac{d}{2}}} \exp\left(-\frac{\mu |x-y|^2}{t-s} \right),
\end{align*}
where we recall that the exact value of $\mu$ is allowed to change from line to line.
This proves that
\begin{equ}
\|\Gamma_{\eps}(\ \cdot\ ; y, s )\|_{C_\fraks^{\alpha}(Q_{\sqrt{ |t-s|}/8}(x,t))} \leq \frac{C}{|t-s|^{\frac{d+\alpha}{2}}} \exp\left(-\frac{\mu|x-y|^2}{t-s} \right)\ .
\end{equ}
The bound on $ \|\Gamma_{\eps}(x,t; \ \cdot\  )\|_{C_\fraks^{\alpha}(\tilde{Q}_{\sqrt{ |t-s|}/8}(y,s))} $ follows similarly using \eqref{eq:adjoint}.
\end{proof}

\begin{prop}\label{prop:uniform hölder bound on kernel}
 For $\alpha,\alpha'\in (0,1)$, there exists $\mu>0$ such that
\begin{equ}
\|\Gamma_{\eps}\|_{C_\fraks^{\alpha, \alpha'}\left(Q_{\sqrt{ |t-s|}/8}(x,t)\times\tilde{Q}_{\sqrt{ |t-s|}/8}(y,s)\right)} \lesssim \frac{1}{|t-s|^{\frac{d+(\alpha+ \alpha')}{2}}} \exp\left(-\frac{\mu|x-y|^2}{t-s} \right)\ 
\end{equ}
uniformly over $x,y\in \mathbb{R}^d$, $-\infty <s<t<\infty$ and $\eps\in (0,1]$.
\end{prop}

\begin{proof}
We apply Theorem~\ref{thm:technical} on $Q_{\sqrt{ |t-s|}/8}(x,t)$ with $f=0$, $r= \sqrt{ |t-s|}/8$ to the increment $\Gamma(\cdot, z')- \Gamma(\cdot, \bar{z}')$ for $z', \bar{z}'\in \tilde{Q}_{\sqrt{ |t-s|}/8}(y,s)$. Thus,
\begin{align}\label{eq:local_in_proof}
\|\Gamma(\cdot, z')- \Gamma(\cdot, \bar{z}')\|_{C_\fraks^{\alpha}(Q_{\sqrt{ |t-s|}/8}(x,t)} &\lesssim r^{-\alpha} \sup_{z\in Q_{\sqrt{ |t-s|}/4}(x,t))} |\Gamma(z, z')- \Gamma(z, \bar{z}')|\ .
\end{align}
Therefore,
\begin{align*}
\|\Gamma_{\eps}\|_{C_\fraks^{\alpha, \alpha'}(Q_{\sqrt{ |t-s|}/8}(x,t)\times\tilde{Q}_{\sqrt{ |t-s|}/8}(y,s))} 
&= \sup_{z',\bar{z}'\in\tilde{Q}_{\sqrt{ |t-s|}/8}(y,s)} \frac{\|\Gamma(\cdot, z')- \Gamma(\cdot, \bar{z}')\|_{C_\fraks^{\alpha}(Q_{\sqrt{ |t-s|}/8}(x,t))} }{|z'-\bar{z'}|^{\alpha'}} \\
 &\lesssim r^{-\alpha} \sup_{z\in Q_{\sqrt{ |t-s|}/4}(x,t)}\sup_{z',\bar{z}'\in\tilde{Q}_{\sqrt{ |t-s|}/8}(y,s)} \frac{ |\Gamma(z, z')- \Gamma(z, \bar{z}')|\ }{|z'-\bar{z'}|^{\alpha'}} \\
&\lesssim 		r^{-\alpha} \sup_{(\zeta, \tau)\in Q_{\sqrt{ |t-s|}/4}(x,t)}		 \frac{1}{|\tau-s|^{\frac{d+\alpha}{2}}} \exp\left(-\frac{\mu |\zeta-y|^2}{\tau-s} \right) \\
&\lesssim	 \frac{1}{|t-s|^{\frac{d+\alpha+\alpha'}{2}}} \exp\left(-\frac{\mu |x-y|^2}{t-s} \right) \ ,
\end{align*}
where we used \eqref{eq:local_in_proof} in the first inequality and Proposition~\ref{prop:uniform_hölder_bound} in the second inequality.
\end{proof}

We define for $I,J\in \{0,1\}$
\begin{equ}\label{eq:GammaIJ}
\Gamma^{I,J}_\eps(x,t;y,s)= \left(1+\eps  \mathbf{1}_{\{I>0\}} \sum_{i=1}^d \Phi_i^\eps (x,t) \partial_{x_i}\right)
\left(1+ \eps \mathbf{1}_{\{J>0\}}\sum_{j=1}^d \tilde{\Phi}_j^\eps (y,-s) \partial_y\right)
 \bar{\Gamma} (x,t; y, s) \ ,
\end{equ}
where ${\Phi}_j^\eps$ and $\tilde{\Phi}_j^\eps$ were defined in \eqref{eq:correctors} and \eqref{eq:tilde correctors} respectively.

\begin{prop}\label{prop:homgenize alpha norm}
Let $R>0$. If
\begin{equ}
(\partial_t- \mathcal{L}_\eps) u_\eps= (\partial_t- \bar{\mathcal{L}}) \bar u
\end{equ}
on $Q_{2r}(x_0,t_0)$, then 
\begin{align*}
\|& u_\eps -\bar{u}- \eps \sum_{i=1}^d \Phi_i^\eps \partial_i \bar{u}  \|_{C_\fraks^{\alpha}(Q_{r}(x_0,t_0))}\\
&\lesssim  \frac{1}{r^\alpha}  \left(\fint_{Q_{2r}(x_0,t_0)} |u_\eps-\bar{u}|^2\right)^{1/2}   +r^{1-\alpha}\eps^2 \sup_{Q_{2r}(x_0,t_0)}\left(   |\nabla^3 \bar{u}| + |\nabla \partial_t \bar{u}|\right)  +\frac{\eps}{r^\alpha} \sup_{z\in Q_{2r}(x_0,t_0)} |\nabla \bar{u}(z)| \\
&+\left( \eps^{2-\alpha} + \frac{r\eps + \eps^2}{r^\alpha}\right) \sup_{z\in Q_{2r}(x_0,t_0)} |\nabla^2 \bar{u}(z)| +\eps^2 \| \nabla^2 \bar{u}\|_{C_\fraks^{\alpha}(Q_{r}(x_0,t_0))} \ ,
\end{align*}
uniformly over $r<R$.
\end{prop}

\begin{proof}
Set
\begin{equ}
w_\eps= u_\eps-\bar{u}-\eps\sum_{i=1}^d \Phi_i^\eps \partial_i \bar{u} - \eps^2 \sum_{i,j=1^d} \Psi^\eps_{d+1,i,j} \partial_i \partial_j \bar{u} \ , 
\end{equ}
then one finds that 
$\partial_{t}+\mathcal{L}_\eps w_\eps= \eps \nabla\cdot F_\eps\ ,$
where
\begin{align*}
F_{\eps, i}(x,t)& = ( a_{i,j}^\eps \Phi_k^\eps +\Psi_{i,k,j}^\eps )\partial_j \partial_k \bar{u}
+ \eps\Psi_{i,d+1,j}^\eps \partial_t \partial_j \bar{u}\\
&\quad +a_{i,j}^\eps (\eps\partial_j \Psi_{d+1,l,k}^\eps) \partial_l \partial_k \bar{u}
+\eps a_{i,j}^\eps \Psi_{d+1,l,k}^\eps \partial_j \partial_l \partial_k \bar{u} \ ,
\end{align*}
see \cite[Thm~2.2]{GS17}.
Therefore by Theorem~\ref{thm:technical} we find that 
\begin{align*}
\|w_\eps\|_{C_\fraks^{\alpha}(Q_{r}(x_0,t_0))} &\leq Cr^{1-\alpha}\left( \frac{1}{r}  \left(\fint_{Q_{2r}(x_0,t_0)} |w_\eps|^2\right)^{1/2}  + \eps\left(\fint_{Q_{2r}(x_0,t_0)} |F_\eps|^p\right)^{1/p} \right)\\
\lesssim& r^{1-\alpha}\Bigg( \frac{1}{r}  \left(\fint_{Q_{2r}(x_0,t_0)} |w_\eps|^2\right)^{1/2}  + \eps\left(\fint_{Q_{2r}(x_0,t_0)} |\nabla^2 \bar{u}|^p\right)^{1/p}\\
&\qquad +\eps^2\left(\fint_{Q_{2r}(x_0,t_0)} |\nabla^3 \bar{u}|^p\right)^{1/p} +\eps^2\left(\fint_{Q_{2r}(x_0,t_0)} |\nabla \partial_t \bar{u}|^p\right)^{1/p}
 \Bigg)\ .
\end{align*}
Since one easily checks that
\begin{equ}
\| \Psi^\eps_{d+1,i,j} \partial_i \partial_j \bar{u}\|_{C_\fraks^{\alpha}(Q_{r}(x_0,t_0))} \lesssim \| \partial_i \partial_j \bar{u}\|_{C_\fraks^{\alpha}(Q_{r}(x_0,t_0))} + \eps^{-\alpha}  \sup_{z\in Q_{r}(x_0,t_0)} |\partial_i \partial_j \bar{u}(z)|
\end{equ}
we conclude that 
\begin{align*}
\|& u_\eps -\bar{u}- \eps \sum_{i=1}^d \Phi_i^\eps \partial_i \bar{u}  \|_{C_\fraks^{\alpha}(Q_{r}(x_0,t_0))} \\
&\lesssim r^{1-\alpha}\Bigg( \frac{1}{r}  \left(\fint_{Q_{2r}(x_0,t_0)} |w_\eps|^2\right)^{1/2}  + \eps\left(\fint_{Q_{2r}(x_0,t_0)} |\nabla^2 \bar{u}|^p\right)^{1/p}\\
&\quad +\eps^2\left(\fint_{Q_{2r}(x_0,t_0)} |\nabla^3 \bar{u}|^p\right)^{1/p} +\eps^2\left(\fint_{Q_{2r}(x_0,t_0)} |\nabla \partial_t \bar{u}|^p\right)^{1/p} \Bigg) \\
&\quad+\eps^2 \| \nabla^2 \bar{u}\|_{C_\fraks^{\alpha}(Q_{r}(x_0,t_0))} + \eps^{2-\alpha}  \sup_{z\in Q_{r}(x_0,t_0)} |\nabla^2 \bar{u}(z)| \ .
\end{align*}
Finally, note that
\begin{equs}
\left(\fint_{Q_{2r}(x_0,t_0)} |w_\eps|^2\right)^{1/2} &\leq \left(\fint_{Q_{2r}(x_0,t_0)} |u_\eps-\bar{u}|^2\right)^{1/2} +\eps \sup_{z\in Q_{2r}(x_0,t_0)} |\nabla \bar{u}(z)|\\
&\qquad+\eps^2 \sup_{z\in Q_{2r}(x_0,t_0)} |\nabla^2 \bar{u}(z)|
\end{equs}
and thus
\begin{align*}
\| &u_\eps -\bar{u}- \eps \sum_{i=1}^d \Phi_i^\eps \partial_i \bar{u}  \|_{C_\fraks^{\alpha}(Q_{r}(x_0,t_0))}\\
&\lesssim r^{1-\alpha}\Bigg( \frac{1}{r}  \left(\fint_{Q_{2r}(x_0,t_0)} |u_\eps-\bar{u}|^2\right)^{1/2}  + \eps\left(\fint_{Q_{2r}(x_0,t_0)} |\nabla^2 \bar{u}|^p\right)^{1/p}\\
& +\eps^2\left(\fint_{Q_{2r}(x_0,t_0)} |\nabla^3 \bar{u}|^p\right)^{1/p} +\eps^2\left(\fint_{Q_{2r}(x_0,t_0)} |\nabla \partial_t \bar{u}|^p\right)^{1/p} \Bigg) \\
&+\eps^2 \| \nabla^2 \bar{u}\|_{C_\fraks^{\alpha}(Q_{r}(x_0,t_0))} + \eps^{2-\alpha}  \sup_{z\in Q_{r}(x_0,t_0)} |\nabla^2 \bar{u}(z)|\\
&+\frac{\eps}{r^\alpha} \sup_{z\in Q_{2r}(x_0,t_0)} |\nabla \bar{u}(z)| + \frac{\eps^2}{r^\alpha} \sup_{z\in Q_{2r}(x_0,t_0)} |\nabla^2 \bar{u}(z)| \ .
\end{align*}
The proof is completed after applying the inequality $$\left(\fint_{Q_{2r}(x_0,t_0)} |g|^p\right)^{1/p} \leq \|g\|_{\L^\infty(Q_{2r}(x_0,t_0))}$$ and collecting terms.
\end{proof}

\begin{prop}\label{prop:alpha-convergence for one sided increment}
 For $\alpha\in (0,1)$ there exists $\mu>0$ such that
\begin{align*}
&\|\Gamma_{\eps}(\ \cdot\ ; y, s )- {\Gamma}^{1,0}_\eps (\ \cdot\ ; y, s )  \|_{C_\fraks^{\alpha}(Q_{\sqrt{ |t-s|}/8}(x,t))}+ \|\Gamma_{\eps}(x,t; \ \cdot\  )- \Gamma^{0,1}_{\eps}(x,t; \ \cdot\  )\|_{C_\fraks^{\alpha}(\tilde{Q}_{\sqrt{ |t-s|}/8}(y,s))} \\
&\qquad \lesssim \frac{\eps}{|t-s|^{\frac{d+1+\alpha}{2}}} \exp\left(-\frac{\mu|x-y|^2}{t-s} \right)\ ,
\end{align*} 
uniformly over $x,y\in \mathbb{R}^d$, $-\infty <s<t<\infty$ and $\eps\in (0,1]$.
\end{prop}

\begin{proof}
We only show the bound for the first term on the left-hand side since
the other one follows analogously. 
First consider the case $\eps \leq r:= \sqrt{ |t-s|}/8$.
Let $u_\eps(\cdot)= \Gamma(\cdot \ ;y,s )$ and $\bar{u}(\cdot)= \bar{\Gamma}(\cdot \ ;y,s )$ on
 $Q_{\sqrt{ |t-s|}/4}(x,t)$. Then it follows from Proposition~\ref{prop:homgenize alpha norm} that 
\begin{align*}
 \| u_\eps -\bar{u}- \eps \sum_{i=1}^d \Phi_i^\eps \partial_i \bar{u}  \|_{C_\fraks^{\alpha}(Q_{r}(x,t))}
\lesssim & \frac{1}{r^\alpha}  \left(\fint_{Q_{2r}(x,t)} |u_\eps-\bar{u}|^2\right)^{1/2}  \\
& +r^{2-\alpha}\eps \sup_{Q_{2r}(x,t)}\left(   |\nabla^3 \bar{u}| + |\nabla \partial_t \bar{u}|\right)  \\
&+\frac{\eps}{r^\alpha} \sup_{z\in Q_{2r}(x,t)} |\nabla \bar{u}(z)| + r^{1-\alpha}\eps \sup_{z\in Q_{2r}(x,t)} |\nabla^2 \bar{u}(z)| \\
&+\eps^2 \| \nabla^2 \bar{u}\|_{C_\fraks^{\alpha}(Q_{r}(x,t))} ,
\end{align*}
where we used $\eps<r$. 
Thus, we conclude similarly to the proof of Theorem~\ref{thm:uniform}, but using \eqref{eq:GS20thm1} to bound the first term.
For $\eps>r=\sqrt{ |t-s|}/8$ note that
\begin{align*}
\|\Gamma_{\eps}(\ \cdot\ ; y, s )- {\Gamma}^{1,0}_\eps (\ \cdot\ ; y, s )  \|_{C_\fraks^{\alpha}(Q_{\sqrt{ |t-s|}/8}(x,t))} 
\leq &
\|\Gamma_{\eps}(\ \cdot\ ; y, s )\|_{C_\fraks^{\alpha}(Q_{\sqrt{ |t-s|}/8}(x,t))}
+ \| \bar{\Gamma}\|_{C_\fraks^{\alpha}(Q_{\sqrt{ |t-s|}/8}(x,t))}\\
&
+\eps \sum_{i=1}^d \|  \Phi_i^\eps (x,t) \partial_{i} \bar{\Gamma}(\ \cdot\ ; y, s )\|_{C_\fraks^{\alpha}(Q_{\sqrt{ |t-s|}/8}(x,t))}\ .
\end{align*}
The claim follows by estimating the first summand using Proposition~\ref{prop:uniform_hölder_bound} and the remaining ones directly.
\end{proof}

The following is a direct consequence of Proposition~\ref{prop:alpha-convergence for one sided increment}.
\begin{corollary}\label{cor:for stochastic estimates}
For $\alpha\in (0,1)$ there exists $\mu>0$ such that 
\begin{align*}
&\|\Gamma_{\eps}(\ \cdot\ ; y, s )- \bar{\Gamma}(\ \cdot\ ; y, s )  \|_{C_\fraks^{\alpha}(Q_{\sqrt{ |t-s|}/8}(x,t))} +\|\Gamma_{\eps}(x,t; \ \cdot\  )- \bar{\Gamma}(x,t; \ \cdot\  )\|_{C_\fraks^{\alpha}(\tilde{Q}_{\sqrt{ |t-s|}/8}(y,s))}\\
& \qquad
 \lesssim  \left(  \frac{\eps}{|t-s|^{\frac{d+1+\alpha}{2}}} \vee \frac{\eps^{1-\alpha}}{|t-s|^{\frac{d+1}{2}}} \right)  \exp\left(-\frac{\mu|x-y|^2}{t-s} \right)\
\end{align*}
uniformly over $x,y\in \mathbb{R}^d$, $-\infty <s<t<\infty$ and $\eps\in (0,1]$.
\end{corollary}

For the purpose of this article the next proposition is the main input from periodic homogenisation theory.
\begin{prop}\label{prop:Holder convergence kernel}
For $\alpha, \alpha'\in (0,1)$, there exists $\mu>0$ such that 
\begin{align*}
\|\Gamma_{\eps}-\Gamma^{1,1}_\eps\|_{C_\fraks^{\alpha, \alpha'}(Q_{\sqrt{ |t-s|}/8}(x,t)\times\tilde{Q}_{\sqrt{ |t-s|}/8}(y,s))}&
 \lesssim \frac{\eps}{|t-s|^{\frac{d+1+\alpha+ \alpha'}{2}}} \exp\left(-\frac{\mu|x-y|^2}{t-s} \right)\\
 &+ \frac{\eps^2}{|t-s|^{\frac{d+2+\alpha+ \alpha'}{2}}} \exp\left(-\frac{\mu|x-y|^2}{t-s} \right)\ ,
\end{align*}
uniformly over $x,y\in \mathbb{R}^d$, $-\infty <s<t<\infty$ and $\eps\in (0,1]$.
\end{prop}

\begin{proof}
We proceed similarly to the proof of Proposition~\ref{prop:uniform hölder bound on kernel}. Consider first the case
$\eps \leq r:= \sqrt{ |t-s|}/8$. For $z', \bar{z}'\in \tilde{Q}_{r}(y,s)$, let
\begin{equ}
u_\eps^{z', \bar{z}'} = \Gamma_\eps(\cdot, z')- \Gamma_\eps(\cdot, \bar{z}')\qquad 
\text{and}\qquad
 u_{0;\eps}^{z', \bar{z}'} = \Gamma^{0,1}_\eps(\cdot, z')- \Gamma^{0,1}_\eps(\cdot, \bar{z}')\ .
\end{equ}
Thus, by Proposition~\ref{prop:homgenize alpha norm}
\begin{equs}
&\|\Gamma_{\eps}-\Gamma_\eps^{1,1}\|_{C_\fraks^{\alpha, \alpha'}(Q_{r}(x,t)\times\tilde{Q}_{r}(y,s))} \\
&= \sup_{z',\bar{z}'\in\tilde{Q}_{r}(y,s)} \frac{\|\Gamma_\eps(\cdot, z')- \Gamma_\eps(\cdot, \bar{z}') -\left(\Gamma^{1,1}_\eps(\cdot, z')- \Gamma^{1,1}_\eps(\cdot, \bar{z}') \right)\|_{C_\fraks^{\alpha}(Q_{r}(x,t))} }{|z'-\bar{z'}|^{\alpha'}} \\
&= \sup_{z',\bar{z}'\in\tilde{Q}_{r}(y,s)} \frac{\|u_\eps^{z', \bar{z}'} -u_{0;\eps}^{z', \bar{z}'}- \eps \sum_{i=1}^d \Phi_i^\eps \partial_i u_{0;\eps}^{z', \bar{z}'} \|_{C_\fraks^{\alpha}(Q_{r}(x,t))} }{|z'-\bar{z'}|^{\alpha'}} \\
&\lesssim  \frac{1}{r^\alpha}  \sup_{z\in Q_{2r}(x,t)}\sup_{z',\bar{z}'\in\tilde{Q}_{r}(y,s)}   \frac{ |u_\eps^{z', \bar{z}'}(z)-u_{0;\eps}^{z', \bar{z}'}(z)|}{|z'-\bar{z'}|^{\alpha'}}  \\
&\qquad+ r^{2-\alpha}\eps \sup_{z\in Q_{2r}(x,t)} \sup_{z',\bar{z}'\in\tilde{Q}_{r}(y,s)}\frac{|\nabla^3 u_{0;\eps}^{z', \bar{z}'} (z)| + |\nabla \partial_t u_{0;\eps}^{z', \bar{z}'}(z)|}{|z'-\bar{z'}|^{\alpha'}}   \\
&\qquad+
\sup_{z\in Q_{2r}(x,t)}
\sup_{z',\bar{z}'\in\tilde{Q}_{r}(y,s)}\left(\frac{\eps}{r^\alpha} 
 \frac{|\nabla u_{0;\eps}^{z', \bar{z}'}(z)| }{|z'-\bar{z'}|^{\alpha'}} 
+ r^{1-\alpha}\eps 
\frac{ |\nabla^2 u_{0;\eps}^{z', \bar{z}'}(z)| }{|z'-\bar{z'}|^{\alpha'}} 
 \right) \\
&\qquad+\eps^2 \sup_{z',\bar{z}'\in\tilde{Q}_{r}(y,s)} \frac{1}{|z'-\bar{z'}|^{\alpha'}}  \| \nabla^2 u_{0;\eps}^{z', \bar{z}'}\|_{C_\fraks^{\alpha}(Q_{r}(x,t))} \ .
\end{equs}
Bounding each term separately, we find by Proposition~\ref{prop:alpha-convergence for one sided increment}
\begin{equ}
\frac{1}{r^\alpha}  \sup_{z\in Q_{2r}(x,t)}\sup_{z',\bar{z}'\in\tilde{Q}_{r}(y,s)}   \frac{ |u_\eps^{z', \bar{z}'}(z)-\bar{u}^{z', \bar{z}'}(z)|}{|z'-\bar{z'}|^{\alpha'}} \lesssim  \frac{C\eps}{|t-s|^{\frac{d+1+\alpha+\alpha'}{2}}} \exp\left(-\frac{\mu|x-y|^2}{t-s} \right)
\end{equ}
and using that $\eps<r$ on each of the remaining terms the same upper bound.
For
$\eps > r:= \sqrt{ |t-s|}/8$ and setting $A = Q_{r}(x,t)\times\tilde{Q}_{r}(y,s)$, note that the claim follows from 
\begin{align*}
\|\Gamma_{\eps}- \Gamma^{1,1}_{\eps}\|_{C_\fraks^{\alpha, \alpha'}(A)}
&\leq
\|\Gamma_{\eps}\|_{C_\fraks^{\alpha, \alpha'}(A)}
+ \| 
 \bar{\Gamma} \|_{C_\fraks^{\alpha, \alpha'}(A)}
 +\eps \sum_{i=1}^d \|  \Phi_i^\eps  \partial_{i;1} 
 \bar{\Gamma} \|_{C_\fraks^{\alpha, \alpha'}(A)}\\
&\quad +\eps\sum_{j=1}^d \|  \tilde{\Phi}_j^\eps  \partial_{j;2}
 \bar{\Gamma}\|_{C_\fraks^{\alpha, \alpha'}(A)}
 +\eps^2  \sum_{i,j=1}^d  \|  \Phi_i^\eps  \tilde{\Phi}_j^\eps  \partial_{i;1} \partial_{j;2}
 \bar{\Gamma} \|_{C_\fraks^{\alpha, \alpha'}(A)} \ ,
\end{align*}
using Proposition~\ref{prop:uniform hölder bound on kernel} to bound $\|\Gamma_{\eps}\|_{C_\fraks^{\alpha, \alpha'}(A)}$.
\end{proof}

\subsection{Post-processing of kernel estimates}\label{subsec:modification for spdes}

We fix a cutoff function $\kappa: \mathbb{R}\to [0,1]$ such that 
 \begin{itemize}
 \item $\kappa(t)=0$ for $t<0$ and for $t>2$,
 \item $\kappa(t)=1$ for $t\in (0,1)$,
 \item $\kappa|_{\mathbb{R}_+}$ is smooth,
 \end{itemize}
and write
\begin{equation}\label{eq:rescaled truncations}
\kappa^\eps(t)=\kappa(t/\eps^2) \qquad \text{as well as} \qquad \kappa^\eps_c(t)= \mathbf{1}_{\{t>0\}}\left(1-\kappa^\eps(t) \right) \ .
\end{equation}
%
We define $\bar{\Gamma}_\eps (x,t;\zeta, \tau)= \kappa^\eps_c(t-\tau) \bar{\Gamma} (x,t;\zeta, \tau)$ and 
\begin{equation}\label{eq:Gamma_tilde}
\tilde{\Gamma}_\eps (x,t; \zeta, \tau) = (1+\eps \Phi_\eps(x,t) \nabla_x) (1+\eps \tilde{\Phi}_\eps(\zeta,\tau) \nabla_\zeta) \bar{\Gamma}_\eps (x,t; \zeta, \tau) \ .
\end{equation}
We also fix $\chi: \mathbb{R}^d\to [0,1]$ smooth and compactly supported on $[-2/3, 2/3]^{d}\subset \mathbb{R}^d$ such that $\sum_{k\in \mathbb{Z}^d} \chi(x+k)=1$ for all $x\in \mathbb{R}^d$. Finally, for $\Gamma\in \{\bar{\Gamma},{\Gamma}_{\eps}, \tilde{\Gamma}_{\eps}\}$ set
\begin{equ}\label{eq:periodising kernels}
{K} (t,x;s,y)=  \sum_{k\in \mathbb{Z}^d}\kappa(t-s)\chi(x-y)	{\Gamma}(x,t;y+k,s)\ 
\end{equ}
and denote the resulting kernels by $\bar{K}\ ,{K}_{\eps}$ and $\tilde{K}_{\eps}$ respectively.

%

We shall use the space of kernels, c.f. \cite[Def.~1]{Sin23}, $\bfK^\beta_{L,R}= \{ K\in \bfK \ : \|K\|_{\beta;L,R}<+\infty\}$ for $L,R\in (0,1)$ with the norm
\begin{equ}
\|K\|_{\beta;L,R} =\inf_{\{K_n\}_{n\geq 0}} \left(\sup_{n\in \mathbb{N}} \frac{\|K_n\|_{\L^\infty }}{2^{n(|\fraks|- \beta) }}   
+ \sup_{n\in \mathbb{N}} \frac{\|K_n\|_{C^{0,R}_\fraks}}{2^{n(|\fraks|- \beta +R)} } 
 +\sup_{n\in \mathbb{N}} \frac{\|K_n\|_{C^{L,0}_\fraks}}{2^{n(|\fraks|- \beta +L)} } 
 + \sup_{n\in \mathbb{N}} \frac{\|K_n\|_{C^{L,R}_\fraks}}{2^{n(|\fraks|- \beta +L+R)} } 
\right)\  ,
\end{equ}
where the infimum is taken over over all kernel decomposition $K(z,z')= \sum_{n\geq 0} K_n(z,z')$ such that each $K_n$ for $n\geq 1$ is supported on 
$ \{(z,z')\in (\mathbb{R}^{d+1})^{\times 2}\ : \ |z-z'|_\fraks \leq 2^{-n+1} \} $ and $K_0$ is supported on $ \{(z,z')\in (\mathbb{R}^{d+1})^{\times 2}\ : \ |z-z'|_\fraks  \leq C \} $ for some $C>0$.

\begin{prop}\label{prop:uniform bounds on modified kernels}
For every $R,L\in (0,1)$ there exists $C>0$ such that
\begin{equ}
\|K_\eps \|_{\beta;L,R} + \|\tilde{K}_\eps \|_{\beta;L,R} + \|\bar K\|_{\beta;L,R} <C
\end{equ}
uniformly over $\eps \in (0,1]$, $\beta\in (0,2]$.
\end{prop}
\begin{proof}
Let $\varphi\in \mathcal{C}^{\infty}_{c}(B_2\setminus B_{1/2})$ be such that $\sum_{n=1}^{\infty}\varphi^{2^{-n}}(x)=1$ for every $x\in B_{1/2}\setminus\{0\}$.
For $K\in \{\bar{K}, K_\eps, \tilde{K}_\eps\} $ set for $n\geq 1$
\begin{equ}
{K}_n(z;z')=  \varphi^{2^{-n}}(z-z') 		K(z;z') 
\end{equ}
and $K_0= K- \sum_{n\geq 1}{K}_n$.
%
We shall show the bound for this decomposition of the kernels.
The bound on $\bar K$ is standard, and the bound on $K_\eps$ follows directly from Theorem~\ref{thm:uniform}, Proposition~\ref{prop:uniform_hölder_bound} and Proposition~\ref{prop:uniform hölder bound on kernel}. Finally, we see from \eqref{eq:Gamma_tilde} that the bound on $\tilde{K}_{\eps,n}$ for $2^{-n}< \eps$ follows directly from the bound on $\bar K$, while for $2^{-n}\geq \eps$ we consider the following two cases separately:
\begin{itemize}
\item $t-s<2\eps$: here, the bounds follow again from the bound on $\bar K$.
\item $t-s\geq2\eps$: here, the bounds follow Corollary~\ref{cor:for stochastic estimates}, resp.\  Proposition~\ref{prop:Holder convergence kernel} and the triangle inequality.
\end{itemize}
\end{proof}
The following Proposition is sufficient to treat the equation considered here. 
\begin{prop}\label{prop:kernels for fixed point convergence}
For every $R,L\in [0,1)$ it holds that
\begin{equation}\label{eq:corrrectet kernel convergence for RS}
\|K_\eps- \tilde{K}_\eps \|_{\beta;L,R} \lesssim \eps^{2-\beta}\vee \eps
\end{equation}
uniformly over $\beta\in [0,2]$. If furthermore $0<R+L<1$, it also holds that
\begin{equation}\label{eq:corrrectet kernel convergence for RS_unmodified kernel}
\|K_\eps- \bar{K} \|_{\beta;L,R} \lesssim \eps^{2-\beta}\vee \eps^{1-R-L}
\end{equation}
uniformly over $\beta\in [0,2]$. 
\end{prop}

\begin{proof}
By interpolating between the bounds of Theorem~\ref{thm:uniform} and Proposition~\ref{propGS20} we find that 
$|{K}_{n,\eps}(z;z')-\bar{K}_{n}(z;z')| \lesssim (\eps^{2-\beta}\vee \eps) 2^{n(|\fraks|-\beta)}\ .$ 
By direct computation, we see that
\begin{equ}
|\tilde{K}_{n,\eps}(z;z')-\bar{K}_{n}(z;z')| \lesssim (\eps^{2-\beta} \vee \eps) 2^{n(|\fraks|-\beta)}
\end{equ}
and thus  $|{K}_{n,\eps}(z;z')-\tilde{K}_{n,\eps}(z;z')|\leq |{K}_{n,\eps}(z;z')-\bar{K}_{n}(z;z')|+ |\tilde{K}_{n,\eps}(z;z')-\bar{K}_{n}(z;z')|\lesssim (\eps^{2-\beta}\vee \eps) 2^{n(|\fraks|-\beta)}$.
By interpolating between Proposition~\ref{prop:Holder convergence kernel} and Proposition~\ref{prop:uniform hölder bound on kernel} one finds that 
$$\|{K}_{n,\eps}-\tilde{K}_{n}\|_{C^{L,R}_\fraks(Q_2\times Q_2)} \lesssim 
\begin{cases}
  \eps 2^{n(|\fraks|-2+L+R)}  & \text{if $\eps<2^{-n}$,} \\
  \left(\eps^{2-\beta}\vee \eps\right) 2^{n(|\fraks|-\beta+L+R)} &  \text{otherwise.}
\end{cases}
$$
The bounds on $\|{K}_{n,\eps}-\tilde{K}_{n}\|_{C^{L,0}_\fraks(Q_2\times Q_2)}$, $\|{K}_{n,\eps}-\tilde{K}_{n}\|_{C^{0,R}_\fraks(Q_2\times Q_2)}$ follow similarly but using Proposition~\ref{prop:alpha-convergence for one sided increment} and the triangle inequality,
which concludes the proof of \eqref{eq:corrrectet kernel convergence for RS}.
Finally,
$$\|\bar{K}_{n}-\tilde{K}_{n,\eps}\|_{C^{L,R}_\fraks(Q_2\times Q_2)}\lesssim
\begin{cases}
   \eps^{1-L-R} 2^{|\fraks|-1}  & \text{if $\eps<2^{-n}$,} \\
 \left(\eps^{2-\beta}\vee \eps\right)  2^{|\fraks|-\beta+L+R} &  \text{otherwise.}
\end{cases}
 $$
The remaining bounds follow from Corollary~\ref{cor:for stochastic estimates}, which combined with \eqref{eq:corrrectet kernel convergence for RS} yields \eqref{eq:corrrectet kernel convergence for RS_unmodified kernel}.
\end{proof}

\section{Fixed Point Theorem}\label{sec:fixedpoint}

For $v_0 \in \mathcal{C}^{\eta} (\mathbb{T}^d)$ write
\begin{equ}
(K_\eps(t) v_0) (x)= v_0\left(K_\eps(x,t;\  \cdot\ ,0)\right) = v_0\left(\kappa(t){\Gamma_\eps}(x,t;\ \cdot\ ,0)\right) 
\end{equ}
It follows that $ v= K_\eps (v_0)$ satisfies the initial value-problem $v(0)=v_0$,
\begin{equ}
(\partial_t - \mathcal{L}_\eps) v=0, \qquad \text{on }  \mathbb{T}^d \times (0,1)\ . 
\end{equ}
\begin{lemma}\label{lem:inital_cond}
For $-1<\eta < \gamma<1$, $\eta \neq 0$ it holds 
 that $\|K_\eps (\, \cdot\, )v_0\|_{C^{\gamma, \eta}_T}\lesssim_T \|v_0\|_\eta$ uniformly over $\eps\in \mathbb{N}^{-1}$ and 
 $v_0\in \mathcal{C}^{\eta} (\mathbb{T}^d)$.
Furthermore, for $\eta<0$, $\kappa\in [0,1)$ it holds that 
$\|K_\eps(\, \cdot\, ) v_0- \bar K (\, \cdot\, ) v_0\|_{C^{\gamma, \eta- \kappa}_T}\lesssim_T \left(\eps^{\kappa}\vee \eps^{1+\eta-\gamma}\right) \|v_0\|_\eta$.
\end{lemma}

\begin{proof}

For the first part of the Lemma we use Proposition~\ref{prop:uniform bounds on modified kernels} with $\beta=2$.
We see that for $n_t\in \mathbb{N}$ such that $\sqrt{t}\in [2^{-n_t-1}, 2^{n_t})$
\begin{equ}
|(K_\eps (t) v_0) (x)|\leq \sum_{n=0}^{n_t} |v_0\left(K_{\eps,n}(x,t;\ \cdot\ ,0)\right) |\lesssim \sum_{n=0}^{n_t}  2^{-n\eta} \lesssim |t|^{\eta/2 \wedge 0}\ . 
\end{equ}
Similarly we observe that for $t<\tau$
\begin{equs}
|(K_\eps(t) v_0) (x)- (K_\eps(\tau) v_0)(\zeta) |& \leq \sum_{n=0}^{n_t} \left| v_0\left(K_{\eps,n}(x,t;\ \cdot\ ,0)-K_{\eps,n}(\zeta,\tau;\ \cdot\ ,0)\right) \right|\\
& \lesssim (|x-\zeta| + \sqrt{|t-\tau|})^\gamma
  \sum_{n=0}^{n_t}  2^{-n(\eta- \gamma)}\\
&\lesssim (|x-\zeta| + \sqrt{|t-\tau|})^\gamma |t|^{(\eta- \gamma)/2}\ ,
\end{equs}
where in the last line we used that $\eta<\gamma$.
%
%
%

In order to conclude the second part of the lemma, we use \eqref{eq:corrrectet kernel convergence for RS_unmodified kernel} of Proposition~\ref{prop:kernels for fixed point convergence} with $\beta= 2-\kappa$ and find for $n_t$ as above,
\begin{equs}
|(\bar K(t) v_0 &- K_\eps(t) v_0) (x)| \lesssim \sum_{n=0}^{n_t} \left| v_0\left((\bar K_n - K_{\eps,n})(x,t;\ \cdot\ ,0)\right)\right|  
\lesssim \sum_{n=0}^{n_t}  (\eps^{2-\beta}\vee \eps^{1+\eta} )2^{n(2-\beta-\eta)}\\
&\quad \lesssim  (\eps^{\kappa}\vee \eps^{1+\eta})   |t|^{(\eta-\kappa)/2\wedge 0}
=  (\eps^{\kappa}\vee \eps^{1+\eta})   |t|^{(\eta-\kappa)/2}\ ,
\end{equs}
where we used that $\eta<0$.
Similarly, we find
\begin{equ}
|(K_\eps(t)-\bar K(t)) v_0 (x)- (K_\eps (\tau) - \bar K (\tau))v_0(\zeta) |
\lesssim (\eps^{\kappa}\vee \eps^{1+\eta-\gamma})   (|x-\zeta| + \sqrt{|t-\tau|})^\gamma |t|^{(-\kappa+\eta-\gamma)/2}\ ,
\end{equ}
where we used that $-\kappa +\eta -\gamma<0$.
\end{proof}

Next, we shall prove the main fixed point theorem. 
\begin{prop}\label{prop:remainder equation}
Let $\alpha, \eta \in (-\frac{1}{10},0)$, $\kappa\in (0, \frac{1}{10}\wedge (-\eta))$  and $\gamma \in (-\alpha, \frac{1}{4})$.
%
%
%
%
 Then, for every family $\mathbb{X}^\eps = (X^\eps_1,X^\eps_2,X^\eps_3)\in  C^{\alpha}(\mathbb{T}^2\times \mathbb{R})^{\times 3}$ and $v_0^\eps\in C^{\eta+\kappa}(\mathbb{T}^2\times \mathbb{R})$ indexed by $\eps \in [0,1]$,
there exists a $T^*\in (0,1)$ depending on 
$\sup_{\eps\in [0,1]} \|\mathbb{X}^\eps\|_{C^\alpha(\mathbb{T}^d\times(-1,2))}$ and $\sup_{\eps\in [0,1]} \| v_0^\eps \|_{C^{\eta+\kappa} (\mathbb{T}^d)}  $, 
such that there exist unique solutions 
$\{w_{\eps}\}_{\eps \in (0,1]}\in C^{\gamma, \eta}_{T^*}$ to the equations
\begin{equ}\label{eq:remainder_equation}
\partial_t-\mathcal{L}_\eps w_{\eps} = -( w_{\eps})^3 +3 w_{\eps}^2 X^\eps_1 - 3  w_{\eps} X_2^\eps + X_3^\eps\ , \qquad w_{\eps} (0)=v^\eps_{0}\ . 
\end{equ}
for $\eps\in (0,1]$ 
 as well as a unique solution $\bar w  \in C^{\gamma, \eta}_{T^*}$ to
\begin{equ}
(\partial_t-\nabla \cdot \bar{A}\nabla) \bar{w} = -( \bar{w} )^3 +3\bar{w}^2 X^0_1 - 3  \bar{w} X^0_2 + X^0_3\ , \qquad \bar{w} (0)={v}_{0}^0\ .
\end{equ}
Furthermore, it holds that
\begin{equ}
\|w_{\eps} - \bar{w}\|_{C^{\gamma,\eta}_{T^*}} \lesssim \eps^{\kappa} + \|\mathbb{X}^\eps-\mathbb{X}^0\|_{C^\alpha(\mathbb{T}^d\times(-1,2))}+ \|v^\eps_0-  v^0_0\|_{C^\eta(\mathbb{T}^d)} \ ,
\end{equ}
where the implicit constant depends on $\sup_{\eps\in [0,1]} \|X^\eps\|_{C^\alpha(\mathbb{T}^d\times(-1,2))}$ and \linebreak
 $\sup_{\eps\in [0,1]} \| v_0^\eps \|_{C^{\eta+\kappa} (\mathbb{T}^d)}
$.
\end{prop}
For convenience, we opt to use the theory of regularity structures in the proof. This is only done in order to shorten exposition as it allows us to avoid reproving standard lemmata.
\begin{proof}
Consider the regularity structure $\mathcal{T}= (T,G)$ where $T_{\alpha}$ is spanned by three elements $\{\Xi_i\}_{i=1}^3$, $T_0$ is spanned by $\mathbf{1}$, and $G$ is the trivial group consisting of only one element.
We consider models $(\Pi,\Gamma)$ such that
$\Pi_x \mathbf{1}= 1$ and 
$\Pi_x \Xi_i= X_i$ for all $x\in \mathbb{R}^{d+1}$.

Next, writing $v_\eps(x,t):=\left( K_{\eps}(t)v^\eps_0 \right)(x)$ and $\bar{v}(x,t):=\left( \bar{K}(t) v^0_0 \right)(x)$ we note that $v_\eps\in C^{\gamma,\eta+\kappa}_T$ for every $T>0$ by Lemma~\ref{lem:inital_cond}.
Thus, we can reformulate the fixed point problem in the language of regularity structures 
\begin{equ}
W= (I + J^{K_\eps}(\cdot) + \mathcal{N}^{K_\eps}) \mathbf{R}^+ F(W,\Xi) + v_\eps \cdot \mathbf{1}\ ,
\end{equ}
where for $\Xi=  (\Xi_1,\Xi_2, \Xi_3)$ we write $F(z,\Xi)=  -( z )^3 +3z^2 \Xi_1 - 3 z \Xi_2 + \Xi_3$.
Here this simplifies to
\begin{equ}
W(t,x)= w(t,x)\cdot \mathbf{1}\ , \qquad w= K_\eps \left( \mathbf{R}^+ f(w,\mathbb{X}^\eps) \right)  + v_\eps
\end{equ}
for $f(z,\mathbb{X})=  -( z )^3 +3z^2 X_1 - 3 z X_2 + X_3$.
We consider $K_\eps$ as an element of $\bfK^\beta_{L,R}$ for $\beta=1, L=R=\frac{1}{4}$.
In particular, we then find by \eqref{eq:corrrectet kernel convergence for RS_unmodified kernel}, that
\begin{equ}
\|K_\eps- \bar{K} \|_{\beta;L,R} \lesssim \eps^{2-\beta}\vee \eps^{1-R-L}\leq \eps^{\kappa} \ .
\end{equ}
Note that for our choice of $\beta, L,R$ it holds that $L>\gamma,\ R>-\alpha$, thus we can apply \cite[Thm~2.11\& Rem.~2.12]{Sin23} (which is a modification of \cite[Thm~7.8]{Hai14} which includes continuity of the fixed point with respect to the kernel)
for the choice of models determined by $\mathbb{X}= \mathbb{X}^\eps$ and $\mathbb{X}= \mathbb{X}^0$ to obtain
 that $T^*\in (0,1)$ 
and solutions
$w^\eps, \bar{w}\in C^{\gamma, \eta}_{T^*}$ exist, and furthermore
\begin{equ}
\|w_{\eps} - \bar{w}\|_{C^{\gamma,\eta}_{T^*}} \lesssim \eps^{\kappa} +  \|\mathbb{X}^\eps-\mathbb{X}^0\|_{C^\alpha(\mathbb{T}^d\times(-1,2))}+ \|v_\eps- \bar{v}\|_{C^{\gamma, \eta}_{T^*}} \ .
\end{equ}
Since
\begin{equ}
v_\eps(x,t)- \bar{v}(x,t)= \left( K_{\eps}(t)(v^\eps_0 -v^0_0 )\right)(x) + \left( (K_{\eps}(t)- \bar{K}(t))  v^0_0 \right)(x)
\end{equ}
and thus by Lemma~\ref{lem:inital_cond}
\begin{equ}
\|v_\eps- \bar{v} \|_{C^{\gamma,\eta}_{T^*}}
\lesssim 
 \|v^\eps_0- {v}^0_0\|_{ C^\eta} 
 + \eps^\kappa \|{v}^0_0\|_{C^{ \eta+\kappa}} 
 \ ,
\end{equ}
this completes the proof.
\end{proof}

%
%

\subsection{An a priori bound}
In this section we show that the a priori estimate for the $\phi^4_2$ model of \cite[Thm.~6.4]{MWplane} holds uniformly in the homogenisation length-scale.
\begin{theorem}\label{thm:a-priori}
Fix $p\in 2\mathbb{N}$ large enough and $\alpha<0$ such that $|\alpha|$ is small enough (depending on $p$, see \cite[Eq.~6.9]{MWplane} for the precise condition).
Let $T>0$ and $M>0$, then there exists $C>0$ such that if for some $\eps\in(0,1]$, $\gamma>0$, $T'\in (0,T]$ and $\mathbb{X}\in  C^{\alpha}(\mathbb{T}^2\times \mathbb{R})^{\times 3}$ satisfying
$\|\mathbb{X}\|_{C^\alpha(\mathbb{T}^d\times(-1,T+1))}< M$, if
$w\in C^\gamma(\mathbb{T}^2\times [0,T'])$ is a solution to
\begin{equ}
\partial_t-\mathcal{L}_\eps w = -w^3 +3 w^2 X_1 - 3  w X_2 + X_3\ , 
\end{equ}
then
\begin{equ}
\sup_{0\leq t\leq T'  } \|w(t)\|_{\L^p(\mathbb{T}^2)}\leq \| w(0)\|_{\L^p(\mathbb{T}^2)} + C \ .
\end{equ}
\end{theorem}
\begin{remark}
Importantly the constant $C>0$ does not depend on $\eps$.
\end{remark}
\begin{proof}
The proof adapts essentially mutatis mutandis from \cite[Theorem~6.4]{MWplane} using only the qualitative fact that $A(x/\eps,t/\eps^2)$ is Hölder continuous together with uniform ellipticity. We point out the minor modifications that are required: 
\begin{enumerate}
\item By Schauder estimates one finds that $(0,T'] \ni  t  \mapsto \nabla w(t) \in L^\infty (\mathbb{T}^2)$ is continuous.
\item Using Hölder continuity of $w$ one obtains the analogue of \cite[Prop.~6.7]{MWplane} for $0<\kappa <t\leq T'$
\begin{equ}
\langle w(t), \phi \rangle -\langle w(\kappa), \phi \rangle 
= \int^t_\kappa - \langle A(\cdot/\eps, s/\eps^2) \nabla w(s), \nabla \phi \rangle  + \langle F(w(s), \mathbb{X}_s), \phi\rangle \ , ds \ .
\end{equ}
\item Then, arguing exactly as in the proof of \cite[Prop.~6.8]{MWplane} one finds that 
for $0<\kappa <t\leq T'$
\begin{equs}
 \frac{1}{p}\left( \|w(t)\|_{\L^p(\mathbb{T}^2)} - \|w(\kappa)\|_{\L^p(\mathbb{T}^2)}\right) &=  -\int^t_\kappa(p-1)\langle  A(\cdot/\eps, s/\eps^2) \nabla w(s), w(s)^{p-2} \nabla w(s) \rangle \, ds \\
&\qquad+ \int^t_\kappa \langle F(w(s), \mathbb{X}_s), w(s)^{p-1}\rangle \ , ds \ .
\end{equs}
\item Using that $A$ is uniformly elliptic and bounded one concludes as in the proof of \cite[Thm.~6.4]{MWplane} that
\begin{equ}
\sup_{\kappa \leq t\leq T'  } \|w(t)\|_{\L^p(\mathbb{T}^2)}\leq \| w(\kappa)\|_{\L^p(\mathbb{T}^2)} + C \ .
\end{equ}
\end{enumerate}
Letting $\kappa\to 0$ completes the proof.
\end{proof}

\section{Continuity of the model for Theorem~\ref{thm:main1}}\label{sec:stochastic_estimates}

%

Let $\xi_{\eps, \delta}$ be as in \eqref{eq:noise}, 
we define
\begin{equ}\label{eq:linear tree}
\<1>_{\eps,\delta}(z)= \int_{\mathbb{R}^{d+1}} K_\eps(z, z') \xi_{\eps,\delta}(z') dz' \ . 
\end{equ}

\begin{prop}\label{prop:IXi regularity}
Let $d\geq 2$. For $\alpha>0$, $1>\alpha'>0$, there exists a modification of \eqref{eq:linear tree} such that for all $T>0$ the map
\begin{equ}{}
[0,1]^2 \to C^{\alpha'/2}\left( [-T,T], C^{-d/2+1-\alpha- \alpha'}(\mathbb{R}^d) \right), \qquad (\eps, \delta)\mapsto \<1>_{\eps,\delta} 
\end{equ}
is a.s. continuous.
\end{prop}

\begin{proof}
Using the stochastic estimates of Proposition~\ref{prop:stochastic estimates_linear}, the claim follows by applying a Kolmogorov-type criterion, c.f. \cite[Prop.~4.20]{HairerNotes}, to the Banach space valued random variables 
$\left\{\<1>_{\eps, \delta} (\cdot,  t)\right\}_{(\eps, \delta, t)\in [0,1]^{2}\times [-T,T] } \ . $
\end{proof}
\begin{remark}
Note that in even though in this article we are exclusively interested in the case $d=2$, we have stated Proposition~\ref{prop:IXi regularity} for general $d\geq 2$ as it might be useful in other settings. We shall sometimes similarly do specific calculations in proofs for general $d$ to emphasise where $d=2$ is used in the argument.
\end{remark}

\begin{prop}\label{prop:stochastic estimates_linear}
Let $\alpha, \alpha' ,\kappa>0$ such that $\alpha'+\kappa<1$. 
Using the shorthand notation $\triangle_{s,t}\<1>_{\eps, \delta} := \<1>_{\eps, \delta} (\cdot, t)- \<1>_{\eps, \delta} (\cdot, s)$
it holds that
\begin{enumerate}
\item\label{eq:0} $\E[\langle \<1>_{\eps, \delta}(\ \cdot\ , t) , \psi^\lambda_x \rangle^2]^{\frac{1}{2}} \lesssim \lambda^{-\alpha-d/2+1}$
\item\label{eq:1} $\E [\langle \triangle_{s,t}\<1>_{\eps, \delta} , \psi^\lambda_x \rangle^2]^{\frac{1}{2}} \lesssim |t-s|^{\alpha'/2} \lambda^{-\alpha-\alpha'-d/2+1}$
\item\label{eq:1} $\E[\langle \triangle_{s,t}\<1>_{\eps, \delta}-\triangle_{s,t}\<1>_{0, \delta} , \psi^\lambda_x \rangle^2]^{\frac{1}{2}} \lesssim \eps^{\kappa} |t-s|^{\alpha'/2} \lambda^{-\alpha-\alpha'-\kappa-d/2+1}$
\item\label{eq:1} $\E[\langle \triangle_{s,t}\<1>_{\eps, \delta}-\triangle_{s,t}\<1>_{\eps, 0} , \psi^\lambda_x \rangle^2]^{\frac{1}{2}} \lesssim \delta^{\kappa} |t-s|^{\alpha'/2} \lambda^{-\alpha-\alpha'-\kappa-d/2+1}\ ,$
\end{enumerate}
uniformly over $|t-s| \vee \lambda \le 1$, $\eps,\delta \in (0,1]$, $x \in \R^d$, 
and $\psi \in \mathfrak{B}^{0}(\mathbb{R}^{d})$.
\end{prop}
\begin{proof}
We note that for $z=(x,t)$ and writing $\phi^\delta(t)= \frac{1}{\delta^2} \phi(t/\delta^2)$
 \begin{equs}
\<1>_{\eps, \delta} (x,t)&= \int_{\mathbb{R}^{d+1}} K_\eps(x,t; \zeta, \tau) \xi_{\eps,\delta}(\zeta, \tau) d\zeta d\tau \\
&= \int_{\mathbb{R}^{d+1}} \Gamma_\eps(x,t; \zeta, \tau ) \kappa(t-\tau) \xi_{\eps,\delta}(\zeta,\tau) d\zeta d\tau \\
&= \int_{\mathbb{R}^{d+1}} \Gamma_\eps(x,t; \zeta, \tau ) \kappa(t-\tau) 
 \left(
\int_{\mathbb{R}^{d+1}} \phi^\delta(\tau-s) \Gamma_{\eps}(\zeta,\tau; y, \tau-\delta^2) \xi(dy, ds) \right)
d\zeta d\tau \\
&= \int_{\mathbb{R}^{d+1}} \left(\int_{\mathbb{R}} \phi^\delta(\tau-s) \kappa(t-\tau)  \Gamma_\eps(x,t;  y, \tau-\delta^2)  d\tau \right)
 \xi(dy, ds) \\
  &= \int_{\mathbb{R}^{d+1}} \left(\int_{\mathbb{R}} \phi^\delta(\tau) \kappa(t-s-\tau)  \Gamma_\eps(x,t;  y, s+ \tau-\delta^2)  d\tau \right)
 \xi(dy, ds) 
  \\
 &= \int_{\mathbb{R}^{d+1}} H_{\eps, \delta}(x,t;y,s)
 \xi(dy, ds) \ ,  \label{eq:kernel_representation} 
 \end{equs}
where $H_{\eps, \delta}$ is defined in \eqref{eq: H_epsilon,delta}. Thus,
the inequalities follow from Proposition~\ref{prop:H bound} and Proposition~\ref{prop:convolution estimate}.
\end{proof}

\subsection{Renormalisation}\label{sec:renormalisation}
Recall the functions $\kappa^\eps(t)$ and 
 $\kappa^\eps_c(t)$ from \eqref{eq:rescaled truncations}, we define
\begin{equ}\label{eq:mathcalZ}
\mathcal{Z}_\eps(x,t;\zeta, \tau):= \kappa^\eps(t-\tau)\bar Z^\eps(x,t;\zeta,\tau) + \kappa_c^\eps(t-\tau)  \bar{\Gamma} (x,t;\zeta, \tau)\ .
\end{equ}

It follows from Corollary~\ref{cor:small scale bound} and Proposition~\ref{propGS20} that 
\begin{equ}\label{eq:renormalisation_approx}
|(\Gamma_\eps-\mathcal{Z}_\eps)(x,t;\zeta, \tau)|\lesssim \left( 
\frac{\eps^{-\theta}\kappa^\eps(t-\tau)}{(t-\tau)^{d/2-\theta/2}} +  \frac{\eps\kappa_c^\eps(t-\tau) }{(t-\tau)^{(d+1)/2}} \right)
 \exp\left( -\frac{\mu |x-\zeta|^2}{(t-\tau)} \right)\ .
\end{equ}
Similarly it follows from
Corollary~\ref{cor:small scale bound} and Corollary~\ref{cor:for stochastic estimates} that 
\begin{align}\label{eq:renormalisation_approx_increment}
&|(\Gamma_\eps-\mathcal{Z}_\eps)(x,t;\zeta, \tau-\delta^2)-(\Gamma_\eps-\mathcal{Z}_\eps)(x,t;\zeta, \tau)|  \\
&\lesssim \delta^{\theta} \kappa^\eps(t-\tau) 
\left(\left(\frac{\mathbf{1}_{\delta\leq \eps}
\eps^\theta }{ (t-\tau)^{d/2+\theta/2} } +  \frac{\mathbf{1}_{\delta\leq \eps}
 }{ (t-\tau)^{d/2} }   \right) \wedge \frac{\mathbf{1}_{\delta\leq \eps}
\delta^{-\theta}\eps^{-\theta}} { (t-\tau)^{d/2-\theta/2} }  +
\frac{\mathbf{1}_{\delta>\eps} }{ (t-\tau)^{d/2} } \right)
\exp\left( -\frac{\mu |x-\zeta|^2}{t-\tau} \right)\nonumber \\
&+ \delta^{\theta} \frac{\eps^{1-\theta}\kappa_c^\eps(t-\tau)}{(t-\tau)^{(d+1)/2}} 
\exp\left( -\frac{\mu |x-\zeta|^2}{(t-\tau)} \right) \ . \nonumber 
\end{align}
Furthermore, it holds that 
\begin{equation}\label{eq:mathcalZ_bound}
|\mathcal{Z}_\eps(x,t;\zeta, \tau)|\lesssim\frac{1}{(t-\tau)^{d/2}} \exp\left( -\frac{\mu |x-\zeta|^2}{(t-\tau)} \right)\ .
\end{equation}
as well as 
\begin{equation}\label{eq:mathcalZ_bound_increment}
|\mathcal{Z}_\eps(x,t;\zeta, \tau-\delta)-\mathcal{Z}_\eps(x,t;\zeta, \tau)|\lesssim\delta^{\theta} \left( \frac{ 1}{(t-\tau)^{d/2+\theta/2} } \vee \frac{ \eps^{-\theta}}{(t-\tau)^{d/2} } \right)\exp\left( -\frac{\mu |x-\zeta|^2}{(t-\tau)} \right)\ .
\end{equation}
\begin{remark}
Note that $\mathcal{Z}_\eps(x,t;\zeta, \tau)$ is an explicit approximation to $\Gamma_{\eps}$ at scales $\ll\eps$ and $\gg \eps$. 
In dimension $d=2$ this approximation is sufficiently good to capture the divergent behaviour of the singular integral, as quantified by the estimate in Lemma~\ref{lem:bounded_periodic1} on $f_{\eps,\delta}$. Note that in higher dimensions or for other equations this will not necessarily be the case.
\end{remark}
We define
\begin{align*}
f_{\eps,\delta}(x,t)
= \int_\mathbb{R}\int_\mathbb{R}
\left(\phi^\delta\right)^{*2} ( \tau-\tau' ) \kappa(t-\tau) & \kappa(t-\tau') \Bigg( \int_{\mathbb{R}^2} \Gamma_\eps(x,t;\eta, \tau-\delta^2) 
 \Gamma_\eps(x,t;  \eta, \tau'-\delta^2) \\
 & - 
\mathcal{Z}_\eps (x,t;\eta, \tau-\delta^2) \mathcal{Z}_\eps (x,t;\eta, \tau'-\delta^2)
\,d\eta \Bigg) d\tau  d\tau'\ .
\end{align*}
\begin{lemma}\label{lem:bounded_periodic1}
For each $\eps\in (0,1]$ the function $f_{\eps, \delta}$ is $(\eps\mathbb{Z})^2\times (\eps^2\mathbb{Z})$-periodic
and satisfies
\begin{equ}
\sup_{\eps, \delta \in (0,1]^2}\| f_{\eps, \delta}\|_{\L^\infty}<+\infty\ .
\end{equ}
Furthermore, for $\delta_0>0$, it holds that $\lim_{\eps \to 0} \|f_{\eps,\delta_0}\|_{\L^\infty}=0$
and for $\eps_0>0$, the functions $f_{\eps_0,\delta}$ converge pointwise as $\delta \to 0$.
Lastly, it holds that for every $\kappa\in [0, \theta)$
\begin{equ}
\sup_{\eps, \delta \in (0,1]^2} \frac{|f_{\eps, \delta}- f_{\eps,0}|}{\delta^{\kappa} \eps^{-\kappa}}<+ \infty \ .
\end{equ}
\end{lemma}

\begin{proof}
Writing
\begin{align*}
&\Gamma_\eps(x,t;\eta, \tau-\delta^2) 
 \Gamma_\eps(x,t;  \eta, \tau'-\delta^2) - 
\mathcal{Z}_\eps (x,t;\eta, \tau-\delta^2) \mathcal{Z}_\eps (x,t;\eta, \tau'-\delta^2) \\
&=\big(\Gamma_\eps(x,t;\eta, \tau-\delta^2) -\mathcal{Z}_\eps (x,t;\eta, \tau-\delta^2)\big) \Gamma_\eps(x,t;  \eta, \tau'-\delta^2)\\
&\qquad -\mathcal{Z}_\eps (x,t;\eta, \tau-\delta^2)\big( \Gamma_\eps(x,t;  \eta, \tau'-\delta^2)-\mathcal{Z}_\eps (x,t;\eta, \tau'-\delta^2)\big)
\end{align*}
we see by \eqref{eq:renormalisation_approx} and \eqref{eq:mathcalZ_bound} that
\begin{align*}
&|\Gamma_\eps(x,t;\eta, \tau-\delta^2) 
 \Gamma_\eps(x,t;  \eta, \tau'-\delta^2) - 
\mathcal{Z}_\eps (x,t;\eta, \tau-\delta^2) \mathcal{Z}_\eps (x,t;\eta, \tau'-\delta^2)| \\
&\lesssim \left( 
\frac{\eps^{-\theta}\kappa^\eps(t-\tau+\delta^2)}{(t-\tau+\delta^2)^{d/2-\theta/2}} +\frac{\eps\kappa_c^\eps(t-\tau+\delta^2)}{(t-\tau+\delta^2)^{(d+1)/2}} \right)\\
&\qquad \times
\frac{1}{(t-\tau'+\delta^2)^{d/2}} \exp\left( -\frac{\mu |x-\zeta|^2}{(t-\tau+\delta^2)} -\frac{\mu |x-\zeta|^2}{(t-\tau'+\delta^2)} \right)\\
&+ \left( 
\frac{\eps^{-\theta}\kappa^\eps(t-\tau'+\delta^2) }{(t-\tau'+\delta^2)^{d/2-\theta/2}} + 
\frac{\eps\kappa^\eps_c(t-\tau'+\delta^2)}{(t-\tau'+\delta^2)^{(d+1)/2}} \right)\\
&\qquad \times
\frac{1}{(t-\tau+\delta^2)^{d/2}} \exp\left( -\frac{\mu |x-\zeta|^2}{(t-\tau+\delta^2)} -\frac{\mu |x-\zeta|^2}{(t-\tau'+\delta^2)} \right) \ .
\end{align*}
Therefore, by a direct calculation
\begin{align*}
&\int_{\mathbb{R}^d} \Big|\Gamma_\eps(x,t;\eta, \tau-\delta^2) 
 \Gamma_\eps(x,t;  \eta, \tau'-\delta^2) - 
\mathcal{Z}_\eps (x,t;\eta, \tau-\delta^2) \mathcal{Z}_\eps (x,t;\eta, \tau'-\delta^2) \Big|
\,d\eta\\
&\lesssim \left( 
\frac{\eps^{-\theta}  \kappa^\eps(t-\tau+\delta^2) }{(t-\tau+\delta^2)^{-\theta/2}} + \frac{\eps\kappa^\eps_c(t-\tau+\delta^2) }{(t-\tau+\delta^2)^{1/2}} +
\frac{\eps^{-\theta}\kappa^\eps(t-\tau'+\delta^2) }{(t-\tau'+\delta^2)^{-\theta/2}} + \frac{\eps\kappa_c^\eps(t-\tau'+\delta^2) }{(t-\tau'+\delta^2)^{1/2}} \right)\\
&\qquad \times 
\frac{1}{(2t-\tau-\tau'+2\delta^2)^{d/2}}
\end{align*}
and thus
\begin{align*}
&|f_{\eps,\delta}(x,t)|\\
&\lesssim\int_\mathbb{R}\int_\mathbb{R}
\left(\phi^\delta\right)^{*2} ( \tau-\tau' ) 
\left(\frac{\eps^{-\theta}\kappa^\eps(t-\tau+\delta^2) }{(t-\tau+\delta^2)^{-\theta/2}} + \frac{\eps\kappa_c^\eps(t-\tau+\delta^2) }{(t-\tau+\delta^2)^{1/2}} \right)
\frac{\kappa(t-\tau)  \kappa(t-\tau')}{(2t-\tau-\tau'+2\delta^2)^{d/2}} d\tau d\tau'\\
&\lesssim\int_\mathbb{R}\int_\mathbb{R}
\left(\phi^\delta\right)^{*2} ( \tau-\tau' ) 
\left(\frac{\eps^{-\theta}\kappa^\eps(\tau+\delta^2) }{(\tau+\delta^2)^{-\theta/2}} + \frac{\eps\kappa_c^\eps(\tau+\delta^2) }{(\tau+\delta^2)^{1/2}} \right)
\frac{\kappa(\tau)  \kappa(\tau')}{(\tau+\tau'+2\delta^2)^{d/2}} d\tau d\tau'\ . 
\end{align*}
Uniform boundedness in $\eps, \delta$ can thus be read off directly. 
We note that $\lim_{\eps \to 0} \|f_{\eps,\delta_0}\|_{\L^\infty}=0$, since 
 $\kappa(\tau) \kappa^\eps(\tau+\delta^2)= 0$ for all $\tau\in \mathbb{R}$ whenever $\delta>2\eps$.  
Note that $f_{\eps_0,\delta}$ converges pointwise as $\delta\to 0$ by dominated convergence.
Lastly, similarly to above but additionally to \eqref{eq:renormalisation_approx} and \eqref{eq:mathcalZ_bound} using \eqref{eq:renormalisation_approx_increment} and \eqref{eq:mathcalZ_bound_increment} 
it follows that 
$|f_{\eps, \delta}(x,t)- f_{\eps,0}(x,t)|\lesssim \delta^{\kappa} \eps^{-\kappa}$.
\end{proof}

\subsubsection{Counterterms}
We define 
\begin{equ}\label{eq:first_counterterm}
{\alpha}_{\eps,\delta}:=\frac{1}{4\pi}\int_\mathbb{R}\int_\mathbb{R}
\left(\phi^\delta\right)^{*2} ( \tau-\tau' ) \kappa(\tau)  \kappa(\tau') 
\frac{
\kappa^\eps(\tau+\delta^2)\kappa^\eps(\tau'+\delta^2) }{
 (\tau +\tau'+2\delta^2)^{d/2}
 }
   d\tau  d\tau'
\end{equ}
 and
\begin{equ}\label{eq:second_counterterm}
\bar{\alpha}^{(1)}_{\eps,\delta}=  \frac{1}{4\pi } \int_\mathbb{R}\int_\mathbb{R}
\left(\phi^\delta\right)^{*2} ( \tau-\tau' ) \kappa(\tau)  \kappa(\tau') 
\frac{
\kappa_c^\eps(\tau+\delta^2)\kappa_c^\eps(\tau'+\delta^2)}{
 (\tau +\tau'+2\delta^2)^{d/2}
 }
   d\tau  d\tau' \ .
   \end{equ}
These constants are chosen such that the function
\begin{align*}
g_{\eps,\delta}(x,t):&=\int_\mathbb{R}\int_\mathbb{R}
\left(\phi^\delta\right)^{*2} ( \tau-\tau' ) \kappa(t-\tau)  \kappa(t-\tau')\\
&\qquad\qquad \times \Bigg(\int_{\mathbb{R}^2} 
\mathcal{Z}_\eps (x,t;\eta, \tau-\delta^2) \mathcal{Z}_\eps (x,t;\eta, \tau'-\delta^2) 
\,d\eta\Bigg) d\tau  d\tau' \\
&\quad - \alpha_{\eps,\delta} D(x/\eps, t/\eps^2) - \bar{\alpha}^{(1)}_{\eps,\delta} \det (\bar A)^{-1/2}
\end{align*}
satisfies the following bound.
\begin{lemma}\label{lem:bounded_periodic2}
For each $\eps>0$, the function $g_{\eps, \delta}$ is $(\eps\mathbb{Z})^2\times (\eps^2\mathbb{Z})$-periodic and satisfies
\begin{equ}
\sup_{\eps, \delta \in (0,1]^2}\| g_{\eps, \delta}\|_{\L^\infty}<+\infty\ .
\end{equ}
Furthermore, for every $\eps>0$, $\lim_{\delta\to 0} g_{\eps, \delta}= g_{\eps, 0}$ exists.
 For $\delta>0$, 
$\lim_{\eps\to 0} \| g_{\eps, \delta}\|_{\L^\infty}=0\ .$
Lastly, it holds that for $\kappa\in [0,\theta)$
\begin{equ}
\sup_{\eps, \delta \in (0,1]^2} \frac{\| g_{\eps, \delta}- g_{\eps, 0}\|_{\L^\infty}}{\eps^{-\kappa}\delta^{\kappa}}<+\infty\ .
\end{equ}
\end{lemma}


\begin{proof}
A straightforward calculation along the lines of \cite[Sec.~3.1.2]{Sin23} shows that 
\begin{align*}
{\alpha}_{\eps,\delta} D(x/\eps, t/\eps^2)=& \int_\mathbb{R}\int_\mathbb{R}
\left(\phi^\delta\right)^{*2} ( \tau-\tau' ) \kappa(t-\tau)  \kappa(t-\tau') 
\kappa^\eps(t-\tau +\delta^2)
\kappa^\eps(t-\tau' +\delta^2) \\
&\qquad \times\left(\int_{\mathbb{R}^d} 
\bar{Z}_\eps (x,t;\eta, \tau-\delta^2) \bar{Z}_\eps (x,t;\eta, \tau'-\delta^2) 
\,d\eta\right) d\tau  d\tau'
\end{align*}
%
and
\begin{align*}
\bar{\alpha}^{(1)}_{\eps,\delta}= & \int_\mathbb{R}\int_\mathbb{R}
\left(\phi^\delta\right)^{*2} ( \tau-\tau' ) \kappa(t-\tau)  \kappa(t-\tau') 
\kappa_c^\eps(t-\tau +\delta^2)
\kappa_c^\eps(t-\tau' +\delta^2)\\
&\qquad \times\left(\int_{\mathbb{R}^d} 
\bar{\Gamma} (x,t;\eta, \tau-\delta^2) \bar{\Gamma} (x,t;\eta, \tau'-\delta^2) 
\,d\eta\right) d\tau  d\tau'\ .
\end{align*}
Thus 
\begin{align}
g_{\eps, \delta}(x, t)= & 2  \int_\mathbb{R}\int_\mathbb{R}
\left(\phi^\delta\right)^{*2} ( \tau-\tau' ) \kappa(t-\tau)  \kappa(t-\tau') 
\kappa^\eps(t-\tau +\delta^2)
\kappa_c^\eps(t-\tau' +\delta^2) \nonumber
\\
&\qquad \times\left(\int_{\mathbb{R}^d} 
\bar{\Gamma} (x,t;\eta, \tau-\delta^2) \bar{Z}^\eps (x,t;\eta, \tau'-\delta^2) 
\,d\eta\right) d\tau  d\tau' \label{eq:formula for g}
\end{align}
and using the following bound which can be directly read off \eqref{eq:explicit bar Z}
\begin{equ}
|\bar{\Gamma} (x,t;\eta, \tau)| + |\bar{Z}^\eps (x,t;\eta, \tau) |\lesssim (t-\tau)^{-d/2} \exp\left(-\frac{\mu |x-\eta|^2}{t-\tau}\right) 
\end{equ}
a straightforward computation shows that
\begin{equs}
|g_{\eps, \delta}(x, t)| &\lesssim \int_\mathbb{R}\int_\mathbb{R}
\left(\phi^\delta\right)^{*2} ( \tau-\tau' ) \kappa(t-\tau)  \kappa(t-\tau') 
\frac{
\kappa^\eps(t-\tau +\delta^2)
\kappa_c^\eps(t-\tau' +\delta^2)}{
 2(t+\delta^2) -\tau -\tau'
 }
   d\tau  d\tau'\\
   &\lesssim \int_\mathbb{R}\int_\mathbb{R}
\left(\phi^\delta\right)^{*2} ( \tau-\tau' ) \kappa(\tau)  \kappa(\tau') 
\frac{
\kappa^\eps(\tau +\delta^2)
\kappa_c^\eps(\tau' +\delta^2) }{
 2\delta^2 +\tau +\tau'
 }
   d\tau  d\tau'
\end{equs}
which by substitution $(\tau, \tau')\mapsto (\delta^2\tau, \delta^2\tau') $ is seen to be bounded.
For $\eps>0$ we see that $g_{\eps, \delta}$ converges pointwise as $\delta\to 0$ from \eqref{eq:formula for g} and dominated convergence.
For $\delta>0$,  since 
$\kappa(\tau) \kappa^\eps(\tau+\delta^2)=0$ for all $\tau\in \mathbb{R}$ we see that
$g_{\eps, \delta}=0$ whenever $\delta>2\eps$.
Lastly, the bound $|g_{\eps, \delta}(x, t)-g_{\eps, 0}(x, t)|\lesssim \eps^{-\kappa} \delta^{\kappa}$ follows similarly to the first estimate of the lemma from \eqref{eq:formula for g} using the bound
\eqref{eq:bar Z small scale increment} on $\bar{Z}^\eps$ as well as the same bound on $\bar{\Gamma}$.
\end{proof}
We set 
\begin{align}
h_{\eps, \delta}(x,t):=& 
 \int_\mathbb{R}\int_\mathbb{R}
\left(\phi^\delta\right)^{*2} ( \tau-\tau' ) \kappa(t-\tau)  \kappa(t-\tau')\nonumber \\
&\quad\times\left( \sum_{k\in \mathbb{Z}^2\setminus\{0\}} \int_{\mathbb{R}^2} \Gamma_\eps(x,t;\eta+k, \tau-\delta^2) 
 \Gamma_\eps(x,t;  \eta, \tau'-\delta^2) \,d\eta \right) d\tau  d\tau' \ . \label{eq:formula h}
\end{align}

\begin{lemma}\label{lem:bounded_periodic3}
For each $\eps>0$, the function $h_{\eps, \delta}$ is $(\eps\mathbb{Z})^2\times (\eps^2\mathbb{Z})$-periodic and for $\eps, \delta>0$, $z\in \mathbb{R}^d$
the limits $h_{\eps, 0}(z)= \lim_{\nu\to 0 }    h_{\eps, \nu}(z)$, 
$h_{0, \delta}(z)= \lim_{\nu\to 0 }   h_{\nu, \delta} (z)$ exist. It holds that
\begin{equ}\label{eq: h continuity}
\sup_{\eps, \delta \in (0,1]^2} 
\left( \| h_{\eps, \delta}\|_{\L^\infty}+\frac{\| h_{\eps, \delta}- h_{\eps, 0}\|_{\L^\infty}}{\delta^{\theta}}+ \frac{\| h_{\eps, \delta}- h_{ 0, \delta}\|_{\L^\infty}}{\eps^{\theta}}\right)<+\infty\ .
\end{equ}
In particular the map $(0,1]^2\ni (\eps, \delta) \mapsto h_{\eps, \delta}\in C^0$ extends continuously to $[0,1]^2$.
\end{lemma}
\begin{proof}
The existence of the limits $h_{\eps, 0}$ and $h_{0, \delta}$ follows by \eqref{eq:formula h} and dominated convergence.
Uniform bounds on $\| h_{\eps, \delta}\|_{\L^\infty}$ follow from Theorem~\ref{thm:uniform}.  The bound on $\| h_{\eps, \delta}- h_{\eps, 0}\|_{\L^\infty}$ follows from Proposition~\ref{prop:uniform_hölder_bound}. The bound on $\| h_{\eps, \delta}- h_{ 0, \delta}\|_{\L^\infty}$ follows from Proposition~\ref{propGS20}.
The last claim is immediate form \eqref{eq:formula h}.

\end{proof}
Recall the definition of  $\bar{\alpha}^{(1)}_{\eps,\delta}$ in \eqref{eq:second_counterterm} and set
\begin{equ}
\bar{\alpha}^{(2)}_{\eps,\delta}:= \det(\bar{A})^{1/2} \int_{[0,1]^{3}} h_{\eps,\delta}, \qquad \bar{\alpha}_{\eps,\delta}:=\bar{\alpha}^{(1)}_{\eps,\delta}+\bar{\alpha}^{(2)}_{\eps,\delta}\ .
\end{equ}
We also define
\begin{equ}
F_{\eps,\delta}(x,t):= \E [ \<1>_{\eps,\delta}^2 (x,t) ]- \alpha_{\eps,\delta} D(x/\eps, t/\eps^2) - \bar{\alpha}_{\eps,\delta}\det(\bar A)^{-1/2}\;,   
\end{equ}
and observe that
\begin{align*}
F_{\eps,\delta}(x,t)&=f_{\eps, \delta}(x,t)+g_{\eps,\delta}(x,t) +h_{\eps, \delta}(x,t) - \int_{[0,1]^{3}} h_{\eps,\delta}\;.
\end{align*}
We further define for $(\eps,\delta)\in (0,1]^2$
\begin{equ}\label{eq:constant c original}
c_{\eps,\delta}:= \int_{[0,1]^{3}} F_{\eps,\delta} =  \int_{[0,1]^{3}} f_{\eps,\delta}+g_{\eps,\delta}, 
\end{equ}
as well as
\begin{equ}
\hat{F}_{\eps,\delta}= F_{\eps,\delta}- c_{\eps,\delta} \ .
\end{equ}
The next lemma follows directly from Lemmata~\ref{lem:bounded_periodic1},~\ref{lem:bounded_periodic2}.
\begin{lemma}\label{lemma:constant c}
For $\eps>0$, the limit $\lim_{\delta\to 0} c_{\eps, \delta}= \int_{[0,1]^{3}} (f_{\eps,0}+g_{\eps,0})$ exists. 
For $\delta>0$, one has $\lim_{\eps\to 0} c_{\eps, \delta}=0$. 
\end{lemma}

\begin{remark}
A straightforward computation shows that
\begin{align*}
&\left(f_{\eps,0}+g_{\eps,0}\right)(x,t)\\
&= 
\int_{\mathbb{R}^{3}} \kappa(t-\tau) 
\left(\Gamma^2_\eps(x,t;\eta, \tau) - \left(\kappa^\eps (t-\tau) \bar{Z}^{\eps}(x,t;\eta, \tau)\right)^2 
-\left(\kappa^\eps_c (t-\tau) \bar{\Gamma}(x,t;\eta, \tau)\right)^2 \right) d\eta d\tau\ 
\end{align*}
and therefore by substitution 
\begin{align*}
&\left(f_{\eps,0}+g_{\eps,0}\right)(\eps x, \eps^2 t)\\
&=\int_{\mathbb{R}^{3}}  \kappa(\eps^2(t-\tau))   \left(\Gamma^2_1(x, t; \eta, \tau) - \left(\kappa^1 (t-\tau) \bar{Z}^{1}(x, t;\eta, \tau)\right)^2 
-\left(\kappa^1_c ( t-\tau) \bar{\Gamma}( x,  t;\eta, \tau)\right)^2 \right)  d\eta d\tau \ .
\end{align*}
Thus,
$\left(f_{\eps,0}+g_{\eps,0}\right)(\eps x, \eps^2 t)$ converges to
\begin{equ}\label{eq:usefull for continuity when delta=0}
\int   \left(\Gamma^2_1(x, t; \eta, \tau) - \left(\kappa^1 (t-\tau) \bar{Z}^{1}(x, t;\eta, \tau)\right)^2 
-\left(\kappa^1_c ( t-\tau) \bar{\Gamma}( x,  t;\eta, \tau)\right)^2 \right)  d\eta d\tau
\end{equ}
as $\eps \to 0$ and therefore
\begin{equs}\label{eq:constant_finite shift_c}
\lim_{\eps\to 0} c_{\eps,0} = \int_{[0,1]^{3}\times\mathbb{R}^{3}}   \Big( \Gamma^2_1(x, t; \eta, \tau) - &\left(\kappa^1 (t-\tau) \bar{Z}^{1}(x, t;\eta, \tau)\right)^2 \\
&\qquad\qquad -\left(\kappa^1_c ( t-\tau) \bar{\Gamma}( x,  t;\eta, \tau)\right)^2\Big) \,  dx dt d\eta d\tau \ .
\end{equs}
\end{remark}

\begin{lemma}\label{lemma:deterministic part bounded}
For each $\eps>0$, the function $F_{\eps, \delta}$ is $(\eps\mathbb{Z})^2\times (\eps^2\mathbb{Z})$-periodic and satisfies for $\kappa\in [0,\theta)$
\begin{equ}\label{eq:lemma deterministic part bounded}
\sup_{\eps, \delta \in (0,1]^2}\| F_{\eps, \delta}\|_{\L^\infty}
+\sup_{\eps, \delta \in (0,1]^2} \frac{\| F_{\eps, \delta}- F_{\eps, 0}\|_{\L^\infty}}{\eps^{-\kappa}\delta^{\kappa}}<+\infty \ 
\end{equ}
Furthermore, the map $[\nu,1] \ni \eps\mapsto F_{\eps, 0}\in C^0([0,1]^3)$ is $\kappa$-Hölder continuous for each $\nu >0$.
\end{lemma}
\begin{proof}
The estimate \eqref{eq:lemma deterministic part bounded} follows directly from Lemmata~\ref{lem:bounded_periodic1},~\ref{lem:bounded_periodic2} and~\ref{lem:bounded_periodic3}. Continuity of
$(0,1] \ni \eps\mapsto f_{\eps, 0}+g_{\eps, 0} \in C^0$ can be read off of \eqref{eq:usefull for continuity when delta=0} which together with Lemma~\ref{lem:bounded_periodic3} implies the last claim.

\end{proof}

%
%
%
%
%
%
%
%

\begin{lemma}\label{lem:convergence_of_constant_contr}
For every $\kappa\in (0,1)$ and  $\kappa' \in [0,\theta)$, it holds that
\begin{equ}
\left\| \hat{F}_{\eps,\delta}  \right\|_{C^{-\kappa}} \lesssim \eps^\kappa\ , \qquad \left\| \hat{F}_{\eps,\delta} -\hat{F}_{\eps,0} \right\|_{C^{-\kappa}} \lesssim \eps^{\kappa-\kappa'} \delta^{\kappa'} \;,
\end{equ}
uniformly over $(\eps,\delta) \in (0,1]^2$.
Furthermore, for every $\alpha<0$ the map
\begin{equation}\label{eq:continous extension}
(0,1]^{2}\to \mathcal{C}^{-\alpha} , \qquad (\eps,\delta) \mapsto \hat{F}_{\eps,\delta}
\end{equation}
extends continuously to $[0,1]^2$.
\end{lemma}
\begin{proof}
Since $c_{\eps,\delta}=  \int_{[0,1]^{d+1}} F_{\eps,\delta}$, the first part of the lemma is a direct consequence of Lemma~\ref{lem:improved convergence to the mean} in the Appendix. The second part of the lemma follows from the first part together with the continuity of
$(0,1] \ni \eps\mapsto \hat{F}_{\eps, 0}\in C^0$ which is a consequence of Lemma~\ref{lem:convergence_of_constant_contr}.
\end{proof}

\subsection{Renormalised products}
For $\eps\geq 0,\delta>0$ we define
\begin{equs}\label{eq:renormalised products}
\<2>_{\eps,\delta}(z)&:= \<1>_{\eps,\delta}(z)^2- \alpha_{\eps,\delta} D(\mathcal{S}^\eps z) - \hat{\alpha}_{\eps,\delta}\det (\bar A)^{-1/2} -c_{\eps,\delta}\  , \\
 \<3>_{\eps,\delta}(z)&:= \<1>_{\eps,\delta}(z)^{3} - 3\left(\alpha_{\eps,\delta} D(\mathcal{S}^\eps z) - \hat{\alpha}_{\eps, \delta}\det (\bar A)^{-1/2} -c_{\eps,\delta}\right)  \ .
\end{equs}
One has the following.

\begin{prop}\label{prop:renormalised products}
For any $\alpha>0$ there exists a modification of \eqref{eq:renormalised products} such that the map
\begin{equ}{}
[0,1]^2\to  C^{-\alpha}(\mathbb{R}^3)^{\times 2}, \qquad (\eps, \delta)\mapsto (\<2>_{\eps,\delta},\<3>_{\eps,\delta})  
\end{equ}
is a.s.\ continuous.
\end{prop}

\begin{proof}
The proof is standard and follows along the lines of the proof of Proposition~\ref{prop:IXi regularity} applying  Kolmogorov's criterion
to the $C^{-\alpha}((-T,T)^3)$-valued random variables
\begin{equ}
\left\{ \<2>_{\eps, \delta} (\cdot)\right\}_{(\eps, \delta)\in [0,1]^{2}}\ , \ \left\{  \<3>_{\eps, \delta} (\cdot)\right\}_{(\eps, \delta)\in [0,1]^{2}}\ ,
\end{equ}
using equivalence of moments for random variables in a finite Wiener Chaos, together with Proposition~\ref{prop:stochastic estimates}.
\end{proof}

%

\begin{prop}\label{prop:stochastic estimates}
For every $\alpha>0, \kappa\in [0, \theta)$ and $T\in \{\<2>,\<3>\}$, it holds that 
\begin{enumerate}
\item\label{eq:0} $\E[\langle T_{\eps, \delta} , \psi^\lambda_{\star} \rangle^2]^{\frac{1}{2}} \lesssim_{\alpha} \lambda^{-\alpha}$
\item\label{eq:1} $\E[\langle T_{\eps, \delta}-T_{0, \delta} , \psi^\lambda_{\star} \rangle^2]^{\frac{1}{2}} \lesssim_{\alpha,\kappa}\eps^{\kappa}  \lambda^{-\alpha-\kappa}\ $
\end{enumerate} 
and for each $\eps\in [0,1]$ there exist random variables 
$ \langle T_{\eps, 0} , \psi^\lambda_{\star} \rangle $ satisfying the above and
such that furthermore
 \begin{equ}\label{eq:2sec}
 \E[\langle T_{\eps, \delta}-T_{\eps, 0} , \psi^\lambda_{\star} \rangle^2]^{\frac{1}{2}} \lesssim_{\alpha,\kappa} \delta^{\kappa}  \lambda^{-\alpha-\kappa}\ ,
 \end{equ}
uniformly over $\star\in \mathbb{R}^{3}$, $\lambda\in (0,1]$ and $\psi \in \mathfrak{B}^{1}(\mathbb{R}^{3})$.
\end{prop}

\begin{proof}
We write $Q_n: \L^2(\Omega)\to \L^2(\Omega)$ for the projection onto the $n$th Wiener chaos.
Then, for $\delta>0$,
\begin{equ}
\<2>_{\eps, \delta} = Q_2 \<2>_{\eps, \delta} + Q_0 \<2>_{\eps, \delta} = \<1>_{\eps, \delta}^{\diamond 2} + \hat{F}_{\eps,\delta}
\end{equ}
and
\begin{equ}
\<3>_{\eps, \delta} = Q_3 \<3>_{\eps, \delta} + Q_1 \<3>_{\eps, \delta} = \<1>_{\eps, \delta}^{\diamond 3} + 3\hat{F}_{\eps,\delta}  \<1>_{\eps, \delta}\ ,
\end{equ}
where we used $\diamond$ to denote the Wick product. (For $\delta=0$ we define the left-hand side by the right-hand side.)
Thus, setting
\begin{equ}
G_{\eps,\delta}(z,z')= \sum_{k\in \mathbb{Z}^{d}}  \int_{\mathbb{R}^{d+1} } H_{\eps,\delta} (z+k,w) H_{\eps, \delta}(z', w) dw\ , 
\end{equ}
one finds for $k\in \{1,2,3\}$
\begin{equ}
\E[\langle  \<1>_{\eps, \delta}^{\diamond k}, \psi^\lambda_{\star} \rangle^2] = \int\int G^k_{\eps,\delta}(z,z') \psi^\lambda_{\star} (z)\psi^\lambda_{\star}(z') dz dz' .
\end{equ}
The bound of Item~\ref{eq:0} on $\<1>_{\eps, \delta}^{\diamond 2}$ and $\<1>_{\eps, \delta}^{\diamond 3}$ follows directly using Proposition~\ref{prop:H bound} and Proposition~\ref{prop:convolution estimate}. The bounds of 
 Item~\ref{eq:1} and Equation~(\ref{eq:2sec}) on $\<1>_{\eps, \delta}^{\diamond 2}$ and $\<1>_{\eps, \delta}^{\diamond 3}$ follow similarly.
%
The bounds on $Q_0 \<2>_{\eps, \delta}= \hat{F}_{\eps, \delta}$ follow directly from Lemma~\ref{lemma:deterministic part bounded} and Lemma~\ref{lem:convergence_of_constant_contr}.
In order to establish the bound of Item~\ref{eq:0} on 
$Q_1 \<3>_{\eps, \delta}= \hat{F}_{\eps,\delta}  \<1>_{\eps, \delta}$ 
note that 
\begin{equ}\label{eq:bound_local}
 \E[\langle \hat{F}_{\eps,\delta}  \<1>_{\eps, \delta}, \psi^\lambda_{(t,x)} \rangle^2]=\E[\langle \<1>_{\eps, \delta}, \hat{F}_{\eps,\delta} \psi^\lambda_{(t,x)} \rangle^2]\lesssim 
  \E[\langle \<1>_{\eps, \delta}, \psi^\lambda_{(t,x)} \rangle^2].
\end{equ}
In order to establish the remaining bounds, write
\begin{equation}\label{eq:local1}
\hat{F}_{\eps,\delta}  \<1>_{\eps, \delta}= \hat{F}_{\eps,\delta}   \left( \<1>_{\eps, \delta} - \<1>_{0, \delta} \right) + 
\hat{F}_{\eps,\delta}   \<1>_{0, \delta} \ ,
\end{equation}
respectively
\begin{equation}\label{eq:local2}
\hat{F}_{\eps,\delta}  \<1>_{\eps, \delta}-\hat{F}_{\eps,0} \<1>_{\eps, 0}
=\hat{F}_{\eps,\delta} \left( \<1>_{\eps, \delta}-\<1>_{\eps, 0}\right) +  \left(\hat{F}_{\eps,\delta}- \hat{F}_{\eps,0} \right)  \<1>_{\eps, 0} \ .
\end{equation}
Then, the desired bound on the first summand of \eqref{eq:local1} and \eqref{eq:local2} follows as in \eqref{eq:bound_local}. 
Next, note that we can write $G_{\eps,\delta}= \sum_{m} G_{\eps, \delta;m}$ where 
$G_{\eps, \delta;m}(z,z')= \sum_{k\in \mathbb{Z}^d} \left(H_{\eps,\delta} \star H_{\eps, \delta}\right)_m (z+k, z')$
as in Proposition~\ref{prop:convolution estimate}.
Furthermore writing $\hat{\mathbf{F}}_{\eps,\delta}(z,z') =  \hat{F}_{\eps,\delta}(z)\hat{F}_{\eps,\delta}(z')$ and 
$\pmb{\psi}^\lambda(z, z')= \psi^\lambda(z)\psi^\lambda(z')$, it follows that
\begin{align*}
\E[ \langle \hat{F}_{\eps,\delta}   \<1>_{0, \delta},\psi^\lambda  \rangle^2 ]
&=\sum_n \int\int G_{\eps, \delta; n} (z,z') \hat{\mathbf{F}}_{\eps,\delta}(z,z')  \pmb{\psi}^\lambda(z, z')dzdz'
 \end{align*}
 Let $N_\lambda\in \mathbb{N}$ be such that $\lambda\in (2^{-(N_\lambda+1)}, 2^{-N_\lambda}]$,
 then by \eqref{eq:periodic_improved convergence}
\begin{align*}
 \sum_{n=0}^{N_\lambda}  \int\int G_{\eps, \delta; n}(z,z') \hat{\mathbf{F}}_{\eps,\delta}(z,z')  \pmb{\psi}^\lambda(z, z')dzdz'
 &\lesssim \sum_{n=0}^{N_\lambda} \eps^{-\alpha} 2^{n\kappa} \lambda^{-\alpha}\\
 &\lesssim \eps^{-\alpha}\lambda^{-\alpha-\kappa}\ .
  \end{align*}
 On the other hand
\begin{align*} 
 &\sum_{n=N_\lambda}^\infty  \int\int G_{\eps, \delta; n}(z,z') \hat{\mathbf{F}}_{\eps,\delta}(z,z')  \pmb{\psi}^\lambda(z, z') dzdz' \\
 & \leq \|\hat{\mathbf{F}}_{\eps,\delta}\|_{\L^\infty} \sum_{n=N_\lambda}^\infty  \int\int G_{\eps, \delta; n}(z,z') \pmb{\psi}^\lambda(z, z') dzdz' \\
& \leq \|\hat{\mathbf{F}}_{\eps,\delta}\|_{\L^\infty} \sum_{n=N_\lambda}^\infty  \int \psi^\lambda(z)\left(\int G_{\eps, \delta; n}(z,z') \psi^\lambda(z') dz \right) dz'\\
& \lesssim  \|\hat{\mathbf{F}}_{\eps,\delta}\|_{\L^\infty} \sum_{n=N_\lambda}^\infty  \lambda^{-|\fraks|} 2^{-n(-\kappa +|\fraks|) } \, \lesssim \|\hat{\mathbf{F}}_{\eps,\delta}\|_{\L^\infty} \lambda^{-\kappa}
 \end{align*}
Finally, the second term of \eqref{eq:local2} can be bounded similarly.
\end{proof}

Finally, observe that by combining Proposition~\ref{prop:IXi regularity} and Proposition~\ref{prop:renormalised products} we in particular obtain for $\alpha<0$ a model
\begin{equ}
\mathbb{Z}_{\eps,\delta}:= ( \<1>_{\eps,\delta},\  \<2>_{\eps,\delta} , \ \<3>_{\eps,\delta} ) \in C_{\fraks}^\alpha(\mathbb{R}^3)^{\times 3}
\end{equ}
continuous in $(\eps,\delta)\in [0,1]^2$.

\section{Continuity of further models}\label{sec:further stoch est}
In this section we check convergence of the model for the regularisations used in Theorem~\ref{thm:main_translation invariant}, Proposition~\ref{prop:main_flat} and Theorem~\ref{thm:main_non}.

\subsection{Translation invariant regularisation}\label{sec:translation invariant regularisation}
In this section we write 
$\phi^\delta:=\delta^{-d-2}\left(\phi\circ \mathcal{S}^{\delta}_\fraks\right)$
and consider
\begin{equ}
\xi^\flat_{\delta}(z):=\int_{\mathbb{R}^{d+1}} \phi^\delta (z-z') \xi(dz') \qquad \text{and}\qquad \<1>_{\eps,\delta}^\flat (z)= \int_{\mathbb{R}^{d+1}} K_\eps(z, z') \xi^\flat_{\eps,\delta}(z') dz' \ . 
\end{equ}
The proof of the following proposition is an ad verbatim adaptation of the proof of Proposition~\ref{prop:IXi regularity}, replacing $H$ by $H^\flat$, see \eqref{eq: H_epsilon,delta^flat}.
\begin{prop}\label{prop:IXi regularity_flat}
Let $d\geq 2$. For $\alpha>0$, $1>\alpha'>0$, there exists a modification of $\ \<1>_{\eps,\delta}^{\flat}$ such that for all $T>0$ the map
\begin{equ}{}
(0,1]^2 \to C^{\alpha'/2}\left( [-T,T], C^{-d/2+1-\alpha- \alpha'}(\mathbb{R}^d) \right), \qquad (\eps, \delta)\mapsto \<1>^\flat_{\eps,\delta} 
\end{equ}
extends continuously to $[0,1]^2$.
\end{prop}

We return to $d=2$ and define for $(\eps, \delta)\in (0,1]^2$
\begin{align*}
\tilde{D}^{\flat\flat}_{\eps,\delta}(x,t):&= \int_{(\mathbb{R}^{3})^{\times 2}}
\left(\phi^\delta\right)^{*2} (\zeta-\zeta', \tau-\tau' ) \kappa (\tau) \kappa (\tau') \\
&\qquad \qquad \times\Big(
\kappa^\eps(\tau)\kappa^\eps(\tau')\bar Z^\eps(x,t;\zeta,t-\tau)\bar Z^\eps(x,t;\zeta',t-\tau') \Big) d\tau d\zeta d\tau' d\zeta'
\end{align*}
and 
\begin{align*}
\bar{\alpha}^{(1),\flat}_{\eps,\delta}:&= \det(\bar A)^{1/2}  \int_{(\mathbb{R}^{3})^{\times 2}}
\left(\phi^\delta\right)^{*2} (\zeta-\zeta', \tau-\tau' ) \kappa (\tau) \kappa (\tau') \\
&\qquad \qquad \times \Big(\kappa_c^\eps(\tau)\kappa_c^\eps(\tau')   \bar{\Gamma} (x,t;\zeta, t-\tau)  \bar{\Gamma} (x,t;\zeta', t-\tau') \Big)
d\tau d\zeta d\tau' d\zeta' \ .
\end{align*}
Furthermore, writing $z=(x,t)$, $w=(\zeta, \tau)$ and $w'=(\zeta', \tau')$ we define analogue functions to the ones in Section~\ref{sec:renormalisation}:
 \begin{align*}
f^{\flat}_{\eps,\delta}(z):&=\int_{(\mathbb{R}^{3})^{\times 2}} \left(\phi^\delta\right)^{*2} (w-w' )  \kappa (\tau) \kappa (\tau') 
\Big( 
 K_\eps(z;w) 
 K_\eps(z;  w') - 
\mathcal{Z}_\eps (z;w)  \mathcal{Z}_\eps (z;w')  \Big)
dw dw' \ ,\\
g^\flat_{\eps, \delta}(z):&=2 \int_{(\mathbb{R}^{3})^{\times 2}}
\left(\phi^\delta\right)^{*2} (w-w' ) \kappa (\tau) \kappa (\tau') \kappa^\eps(\tau) \kappa_c^\eps(\tau') \bar Z^\eps(z;\zeta,t-\tau)\bar{\Gamma} (z;\zeta', t-\tau') dwdw' \ , \\
h^\flat_{\eps, \delta}(z):&= 
 \int_{(\mathbb{R}^{3})^{\times 2}}
\left(\phi^\delta\right)^{*2} (w-w' )  \kappa (\tau) \kappa (\tau') 
\sum_{k\in \mathbb{Z}^2\setminus\{0\}} \Gamma_\eps(z + k, w) 
 \Gamma_\eps(z; w') dwdw' . 
\end{align*}
Lastly, we set $\bar{\alpha}^{(2),\flat}_{\eps,\delta}:= \det(\bar A)^{1/2}  \int_{[0,1]^3} h^\flat_{\eps, \delta}$ and $\bar{\alpha}^{\flat}_{\eps,\delta}= \bar{\alpha}^{(1),\flat}_{\eps,\delta}+\bar{\alpha}^{(2),\flat}_{\eps,\delta}$ as well as
\begin{equ}
F^{\flat\flat}_{\eps,\delta}(x,t):= \E [ (\<1>_{\eps,\delta}^\flat)^2 (x,t) ]-  \tilde{D}^{\flat\flat}_{\eps,\delta}(z)- \bar{\alpha}^{\flat}_{\eps,\delta} \det(\bar A)^{-1/2} \ , \qquad  c^{\flat\flat}_{\eps,\delta}= \int_{[0,1]^3} F^{\flat\flat}_{\eps, \delta} \ .
\end{equ}
Exactly as in the previous section one finds that
\begin{equ}
F^{\flat\flat}_{\eps,\delta}(z):= f^{\flat}_{\eps,\delta}(z) + g^\flat_{\eps, \delta}(z) + h^\flat_{\eps, \delta}(z)- \int_{[0,1]^3} h^\flat_{\eps, \delta} \ , 
\qquad 
c^{\flat\flat}_{\eps,\delta}= \int_{[0,1]^3} f^{\flat}_{\eps, \delta}+ g^\flat_{\eps, \delta} \ .
\end{equ}
\begin{lemma}\label{lem:change to flat}
The function 
$\hat{F}^{\flat\flat}_{\eps, \delta}(z)= {F}^{\flat\flat}_{\eps, \delta}(z)- c^{\flat\flat}_{\eps,\delta}$ is $(\eps\mathbb{Z})^2\times (\eps^2\mathbb{Z})$-periodic,
has mean $0$, i.e. $\int_{[0,1]^{d+1}} \hat{F}_{\eps, \delta}^{\flat\flat}=0$,
and satisfies for $\kappa \in [0,\theta)$ the bounds 
\begin{equ}\label{bound on F^flat}
\sup_{\eps, \delta \in (0,1]^2}\| \hat{F}^{\flat\flat}_{\eps, \delta}\|_{\L^\infty}
+\sup_{\eps, \delta \in (0,1]^2} \frac{\| \hat{F}^{\flat\flat}_{\eps, \delta}- \hat{F}^{\flat\flat}_{\eps, 0}\|_{\L^\infty}}{\eps^{-\kappa}\delta^{\kappa}}<+\infty \ .
\end{equ}
Furthermore, $\hat{F}_{\eps,0}^{\flat\flat}= \hat{F}_{\eps,0}$ for every $\eps\in [0,1]$. 
\end{lemma}
\begin{proof}
Periodicity and the mean $0$ property are immediate.
The bound \eqref{bound on F^flat} follows by observing that the estimates of Lemma~\ref{lem:bounded_periodic1}, Lemma~\ref{lem:bounded_periodic2} and Lemma~\ref{lem:bounded_periodic3} still hold for the  functions $f^{\flat}_{\eps,\delta}$, $g^\flat_{\eps, \delta}$ and $h^\flat_{\eps, \delta}$, respectively.
The last claim follows by noting that for every $\eps\in (0,1]$, 
$\lim_{\delta\to 0} f^{\flat}_{\eps,\delta}= \lim_{\delta\to 0} f_{\eps,\delta}$, $\lim_{\delta\to 0}g^\flat_{\eps, \delta}=\lim_{\delta\to 0}g_{\eps, \delta}$ and $\lim_{\delta\to 0} h_{\eps, \delta}=\lim_{\delta\to 0} h_{\eps, \delta}$.
\end{proof}
Define
\begin{equs}
D^{\flat\flat}_\lambda(x,t) =& 
\int_{(\mathbb{R}^{3})^{\times 2}}
\left(\phi^{\lambda}\right)^{*2} (\zeta-\zeta', \tau-\tau' ) \kappa (\tau) \kappa (\tau')   \\
&\qquad
\times \Big(
\kappa^1(\tau)\kappa^1(\tau')\bar Z^1(x,t;\zeta,t-\tau)\bar Z^1(x,t;\zeta',t-\tau') \Big) \, d\zeta d\tau d\zeta' d\tau' \ . \label{eq:formula_Dflatflat}
\end{equs}

\begin{lemma}\label{lem:behaviour D flatflat and c flatflat}
The function
$D^{\flat\flat}_{\lambda}$ satisfies \eqref{eq:bound on Dflatflat} and $\lim_{\lambda\to \infty}\|D^{\flat\flat}_{\lambda}\|_{\L^\infty}= 0$. It holds that
$\tilde{D}^{\flat\flat}_{\eps,\delta}= D^{\flat\flat}_{\delta/\eps}\circ \mathcal{S}^\eps_\fraks$.  
Furthermore, the constants $c^{\flat\flat}_{\eps, \delta}$ satisfy  $\lim_{\eps\to 0} c^{\flat\flat}_{\eps, \delta} =0$ for $\delta>0$ and $\lim_{\delta\to 0} c^{\flat\flat}_{\eps, \delta} =\lim_{\delta\to 0} c_{\eps, \delta} $ for $\eps>0$ where $c_{\eps, \delta}$ was defined in \eqref{eq:constant c original}.
\end{lemma}
\begin{proof}
All claims follow by a simple direct computation.
\end{proof}
We set for $\eps, \delta>0$
\begin{equs}
\<2>_{\eps,\delta}^{\flat\flat}(z)&:= (\<1>_{\eps,\delta}^{\flat}(z))^2- D^{\flat\flat}_{\delta/\eps}(\mathcal{S}^\eps z) - \bar{\alpha}^{\flat}_{\eps,\delta}\det (\bar A)^{-1/2} -c^{\flat\flat}_{\eps,\delta}\  , \\
 \<3>^{\flat\flat}_{\eps,\delta}(z)&:= (\<1>^{\flat}_{\eps,\delta}(z))^{3} - 3\left(D^{\flat\flat}_{\delta/\eps}(\mathcal{S}^\eps z) + \bar{\alpha}^{\flat}_{\eps,\delta}\det (\bar A)^{-1/2} +c^{\flat\flat}_{\eps,\delta}\right)\ ,
\end{equs}
and write $\mathbb{Z}^{\flat\flat}_{\eps, \delta}= (\<1>^\flat_{\eps,\delta},\  \<2>^{\flat\flat}_{\eps,\delta}, \ \<3>^{\flat\flat}_{\eps,\delta} )$. 
\begin{prop}\label{prop:model convergence flatflat}
For any $\alpha>0$ there exists a modification such that the map
\begin{equ}
(0,1]^2\to  C^{-\alpha}(\mathbb{R}^3)^{\times 3}, \qquad (\eps, \delta)\mapsto \mathbb{Z}^{\flat\flat}_{\eps, \delta} 
\end{equ}
extends continuously to $[0,1]^2$. Furthermore,  $\mathbb{Z}^{\flat\flat}_{\eps, 0}= \mathbb{Z}_{\eps, 0}$ for every $\eps\in [0,1]$.
\end{prop}
\begin{proof}
The claim about $\<1>^\flat_{\eps,\delta}$ follows from Proposition~\ref{prop:IXi regularity_flat}.
It thus suffices to show the analogue of Proposition~\ref{prop:renormalised products} for $\<2>^\flat_{\eps,\delta}$ and $\<3>^\flat_{\eps,\delta} $, the proof of which adapts mutatis mutandis with the main change being the use of Lemma~\ref{lem:change to flat} instead of \ref{lem:convergence_of_constant_contr}.
Lastly, we observe that $\hat{F}_{\eps,0}^{\flat\flat}= \hat{F}_{\eps,0}$ for every $\eps\in [0,1]$, which implies the last claim.
\end{proof}

%
%
%
Next we define for$(\eps, \delta)\in (0,1]^2$ the new constant $c_{\eps,\delta}^{\flat}= c^{\flat \flat}_{\eps, \delta} + \int_{[0,1]^{3}} D^{\flat\flat}_{\delta/\eps}$ and set
\begin{equs}
\<2>_{\eps,\delta}^{\flat}(z)&:= (\<1>_{\eps,\delta}^{\flat})(z)^2- \bar{\alpha}^{\flat}_{\eps,\delta}\det (\bar A)^{-1/2} -c^{\flat}_{\eps,\delta}\  , \\
 \<3>^{\flat}_{\eps,\delta}(z)&:= (\<1>^{\flat}_{\eps,\delta}(z))^{3} - 3\left(\bar{\alpha}^{\flat}_{\eps,\delta}\det (\bar A)^{-1/2} +c^{\flat}_{\eps,\delta}\right)\ ,  \ 
\end{equs}
and $\mathbb{Z}^{\flat}_{\eps, \delta}= (\<1>^\flat_{\eps,\delta},\  \<2>^{\flat}_{\eps,\delta}, \ \<3>^{\flat}_{\eps,\delta} )$. 
\begin{prop}\label{prop:continuity model_flat}
For any $\alpha>0$ and  $C>0$ there exits a modification such that the map
\begin{equ}
\triangle^{<}_C \to  C^{-\alpha}(\mathbb{R}^3)^{\times 3}, \qquad (\eps, \delta)\mapsto \mathbb{Z}^{\flat}_{\eps, \delta} ,
\end{equ}
extends continuously to the closure of $\triangle^{<}_C$ and it holds that  $\mathbb{Z}^{\flat\flat}_{0, \delta}= \mathbb{Z}^{\flat}_{0, \delta}$ for every $\delta\in [0,1]$.
\end{prop}
\begin{proof}
Observe that for any $C>0$ the function $\tilde{D}^{\flat\flat}_{\eps,\delta}$ is bounded on $\triangle^{<}_C$.
The claim for
\begin{equ}
\<2>_{\eps,\delta}^{\flat}= \<2>_{\eps,\delta}^{\flat\flat} + \tilde{D}^{\flat\flat}_{\eps,\delta} -  \int_{[0,1]^{3}} \tilde{D}^{\flat\flat}_{\eps,\delta}
\end{equ}
thus follows by combining Lemma~\ref{lem:improved convergence to the mean} with Prop~\ref{prop:model convergence flatflat}. 
The claim for $\<3>^{\flat}_{\eps,\delta}=\<1>_{\eps,\delta}^{\flat} \diamond \<2>_{\eps,\delta}^{\flat}$ follows along the line of the last part of the proof of Proposition~\ref{prop:stochastic estimates}.
\end{proof}


%
%
%
%
%

\subsection{Generic non-translation invariant regularisations}
Recall the non-translation invariant regularisation $\xi_{\eps,\delta}^\sharp$ of \eqref{eq:gen_regularisation} and let
\begin{equ}
\<1>_{\eps,\delta}^\sharp (z)= \int_{\mathbb{R}^{d+1}} K_\eps(z, z') \xi^\sharp_{\eps,\delta}(z') dz'  = \int_{{(\mathbb{R}^{d+1})}^{\times 2}} K_\eps(z, z') \varrho^{\eps,\delta} (z'; w) \xi(w)\, dw\ dz' \ .   
\end{equ}
The proof of the following proposition is an ad verbatim adaptation of the proof of Proposition~\ref{prop:IXi regularity}, replacing $H$ by $H^\sharp$, see \eqref{eq:H_epsilon,delta^sharp}.
\begin{prop}\label{prop:IXi regularity_sharp}
Let $d\geq 2$. For $\alpha>0$, $1>\alpha'>0$ and $C>0$, there exists a modification of $\ \<1>_{\eps,\delta}^{\flat}$ such that for all $T>0$ the map
\begin{equ}
\triangle^{<}_C \to C^{\alpha'/2}\left( [-T,T], C^{-d/2+1-\alpha- \alpha'}(\mathbb{R}^d) \right), \qquad (\eps, \delta)\mapsto \<1>^\sharp_{\eps,\delta} 
\end{equ}
extends continuously to the closure of $\triangle^{<}_C $ and $ \<1>^\sharp_{\eps,0} =  \<1>_{\eps,0}$ for all $\eps\in [0,1]$. 
\end{prop}

Returning to $d=2$, let $Z^\flat$ denote $\bar{Z}$ but with $A$ replaced by the identity matrix, i.e. $Z^\flat$ is the usual translation invariant heat kernel. Then we set
\begin{equ}\label{eq:sharp counterterm}
 {\alpha}^{\sharp}_{\eps,\delta} = \int \left( \int \bar{\varrho}^\delta(w,v) \bar{\varrho}^\delta(w',v) dv  \right)\kappa (\tau) \kappa (\tau') 
\kappa^\eps(\tau)\kappa^\eps(\tau')  Z^\flat(w') Z^\flat(w) dw dw' 
\end{equ}
and
\begin{equ}
F^{\sharp}_{\eps,\delta}(x,t):= \E [ (\<1>_{\eps,\delta}^\sharp)^2 (x,t) ]  - {\alpha}^{\sharp}_{\eps,\delta}  D(x/\eps,t/\eps^2)
\  , \qquad c^{\sharp}_{\eps,\delta} = \int_{[0,1]^{3}} F^{\sharp}_{\eps,\delta} \ .
\end{equ}
\begin{lemma}\label{lem:change to sharp}
The function 
$\hat{F}^{\sharp}_{\eps, \delta}(z)= {F}^{\sharp}_{\eps, \delta}(z)- c^{\sharp}_{\eps,\delta}$ is $(\eps\mathbb{Z})^2\times (\eps^2\mathbb{Z})$-periodic, has mean $0$
and satisfies for $\kappa \in [0,\theta)$ the bounds 
\begin{equ}\label{bound on F^flat}
\sup_{\eps, \delta \in (0,1]^2}\| \hat{F}^{\sharp}_{\eps, \delta}\|_{\L^\infty}
+\sup_{\eps, \delta \in (0,1]^2} \frac{\| \hat{F}^{\sharp}_{\eps, \delta}- \hat{F}^{\sharp}_{\eps, 0}\|_{\L^\infty}}{\eps^{-\kappa}\delta^{\kappa}}<+\infty \ .
\end{equ}
Furthermore, one has $\lim_{\delta\to 0} c_{\eps,\delta}^\sharp=\lim_{\delta\to 0} c_{\eps,\delta}$ for $\eps\in (0,1]$ and $\hat{F}_{\eps,0}^{\sharp}= \hat{F}_{\eps,0}$ for every $\eps\in [0,1]$. 
\end{lemma}

\begin{proof}
Using the notation of the beginning of Section~\ref{sec:main_result_non} we set
\begin{equ}
\pmb{\varrho}^{\eps,\delta}(w,w')= \int \varrho^{\eps,\delta}(w,v) \varrho^{\eps,\delta}(w',v) dv ,\quad \pmb{\varrho}^{(\eps,\delta;z)}(w,w') =\int \varrho^{(\eps,\delta;z)}(w,v) \varrho^{(\eps,\delta;z)}(w',v) dv
\end{equ}
and
\begin{equ}
\pmb{\bar{\varrho}}^{\delta}(w,w')= \int \bar{\varrho}^\delta(w,v) \bar{\varrho}^\delta(w',v) dv \ .
\end{equ}
This time we set
\begin{equs}
\tilde{D}^\sharp_{\eps,\delta}(z):= \int_{(\mathbb{R}^{3})^{\times 2}} \Big(
\pmb{\varrho}^{\eps,\delta}(w,w') \kappa (\tau) \kappa (\tau')&
\kappa^\eps(\tau)\kappa^\eps(\tau')\\
&\times\bar Z^\eps(x,t;\zeta,t-\tau)\bar Z^\eps(x,t;\zeta',t-\tau') \Big) d\zeta d\tau d\zeta' d\tau' 
\end{equs}
and
\begin{equs}
O^{\sharp}_{\eps,\delta}:= \det(\bar A)^{1/2} \int_{(\mathbb{R}^{3})^{\times 2}} \Big(
\pmb{\varrho}^{\eps,\delta}(w,w') \kappa (\tau) \kappa (\tau') &\kappa_c^\eps(\tau) \kappa_c^\eps(\tau')\\
&\times \bar{\Gamma} (x,t;\zeta, t-\tau)   \bar{\Gamma} (x,t;\zeta', t-\tau') \Big) d\zeta d\tau d\zeta' d\tau' 
\end{equs}
as well as $f^\sharp_{\eps,\delta}, \ g^\sharp_{\eps,\delta}$ and $h^\sharp_{\eps,\delta}$ by replacing in the definition of the analogous functions in Section~\ref{sec:translation invariant regularisation} $\left(\phi^{\lambda}\right)^{*2} (w-w')$ by $\pmb{\varrho}^{\eps,\delta}(w,w')$.
One observes that analogue estimates of Lemma~\ref{lem:bounded_periodic1}, Lemma~\ref{lem:bounded_periodic2} and Lemma~\ref{lem:bounded_periodic3} still hold for the functions $f^{\sharp}_{\eps,\delta}$, $g^\sharp_{\eps, \delta}$ and $h^\sharp_{\eps, \delta}$ if one replaces $(0,1]$ by $\triangle^{>}_C$ for some $C>0$. Furthermore, such an estimate also holds for $O^{\sharp}$ on $\triangle^{>}_C$.
Next we observe that 
\begin{align*}
&\int_{(\mathbb{R}^{3})^{\times 2}}
\pmb{\varrho}^{(\eps,\delta;(x,t))}(w,w') \kappa (\tau) \kappa (\tau')
\kappa^\eps(\tau)\kappa^\eps(\tau')\bar Z^\eps(x,t;\zeta,t-\tau)\bar Z^\eps(x,t;\zeta',t-\tau') \, d\zeta d\tau d\zeta' d\tau' \\
&= \alpha^\sharp_{\eps,\delta} D(x/\eps,t/\eps^2)\ ,
\end{align*} 
and one checks that 
$\tilde{D}^\sharp_{\eps,\delta}(x,t) - \alpha^\sharp_{\eps,\delta} D(x/\eps,t/\eps^2)$ also satisfies the estimate of Lemma~\ref{lem:bounded_periodic1} on $\triangle^{>}_C$.
Thus, we conclude \eqref{bound on F^flat}. 

The final part of the lemma follows by direct inspection of the definitions combined with the observation that for $\eps>0$ one has $\lim_{\delta\to 0} \|\tilde{D}^\sharp_{\eps,\delta} - \alpha^\sharp_{\eps,\delta} D\circ\mathcal{S}_\fraks^\eps\|_{L^\infty}=0$ by the same argument as in \cite[Section~3.2]{Sin23}.
\end{proof}
Define for $(\eps, \delta)\in (0,1]^2$
\begin{equs}
\<2>_{\eps,\delta}^{\sharp}(z)&:= (\<1>_{\eps,\delta}^{\sharp})(z)^2- {\alpha}^{\sharp}_{\eps,\delta}  (D\circ \mathcal{S}^\eps_\fraks)(z) -c^{\sharp}_{\eps,\delta}\  , \\
 \<3>^{\sharp}_{\eps,\delta}(z)&:= (\<1>^{\sharp}_{\eps,\delta}(z))^{3} - 3\left({\alpha}^{\sharp}_{\eps,\delta}  (D\circ \mathcal{S}^\eps_\fraks)(z) +c^{\sharp}_{\eps,\delta}\right)  \ ,
\end{equs}
and $\mathbb{Z}^{\sharp}_{\eps, \delta}= (\<1>^\sharp_{\eps,\delta},\  \<2>^{\sharp}_{\eps,\delta}, \ \<3>^{\sharp}_{\eps,\delta} )$. 

\begin{prop}\label{prop:renormalised products_sharp}
For any $\alpha>0$ and  $C>0$there exits a modification such that the map
\begin{equ}
\triangle^{>}_C \to  C^{-\alpha}(\mathbb{R}^3)^{\times 3}, \qquad (\eps, \delta)\mapsto \mathbb{Z}^{\sharp}_{\eps, \delta},
\end{equ}
extends continuously to the closure of $\triangle^{>}_C$
and $\mathbb{Z}^{\sharp}_{\eps, 0}= \mathbb{Z}_{\eps,0}$ for every $\eps\in [0,1]$.
\end{prop}
\begin{proof}
The claim for $\<1>_{\eps,\delta}^{\sharp}$ follows from Proposition~\ref{prop:IXi regularity_sharp}. The claim for $\<2>_{\eps,\delta}^{\sharp}$ and $\<3>_{\eps,\delta}^{\sharp}$ follows mutatis mutandis as the proof of Proposition~\ref{prop:renormalised products} but using Lemma~\ref{lem:change to sharp} and the bounds on $H^\sharp$ in Proposition~\ref{prop:H^sharp bound}.
\end{proof}

\section{Proof of the main results}\label{sec:proof of main}
The proofs of Theorem~\ref{thm:main1}, Theorem~\ref{thm:main_non}, Proposition~\ref{prop:main_flat} and Theorem~\ref{thm:main_non} share a common generic part, which we provide first.
\begin{proof}[Generic part of the proofs]
Using the classical Da-Prato Debusche trick, \cite{DD03}, we consider the equation for
\begin{equ}
w_{\eps, \delta} = u_{\eps,\delta} - X_{1}^{\eps, \delta}, \qquad  \textit{ for } X_{1}^{\eps, \delta}\in \{\<1>_{\eps, \delta},\ \<1>^\flat_{\eps, \delta},, \ \<1>^\sharp_{\eps, \delta} \}\ .  
\end{equ}
which is exactly of the form \eqref{eq:remainder_equation} with initial condition $v_0=u_0 - X_{1}^{\eps, \delta}( \cdot,0)$ and
\begin{equ}
\mathbb{X} \in \{ \mathbb{Z}, \ \mathbb{Z}^\flat, \mathbb{Z}^{\flat\flat}, \mathbb{Z}^\sharp \} \ .
\end{equ}
Thus, the
existence of a solution up to a random time $T^*>0$ as well as 
 continuity in $(\eps, \delta)$ as a map into $\L^{0}\left[C([0,T^*), \mathcal{D}'(\mathbb{T}^d))\right]$
 follow
by combining Proposition~\ref{prop:remainder equation} with the appropriate continuity results on models:
\begin{itemize}
\item For Theorem~\ref{thm:main1} this is the content of Proposition~\ref{prop:IXi regularity} and Proposition~\ref{prop:renormalised products}.
\item For Theorem~\ref{thm:main_translation invariant} this the content of Proposition~\ref{prop:IXi regularity_flat} and Proposition~\ref{prop:continuity model_flat}.
\item For Proposition~\ref{prop:main_flat} this is content of Proposition~\ref{prop:IXi regularity_flat} and Proposition~\ref{prop:model convergence flatflat}.
\item For Theorem~\ref{thm:main_non} this is the content of Proposition~\ref{prop:IXi regularity_sharp} and Proposition~\ref{prop:renormalised products_sharp}.
\end{itemize}
Note that for $\delta=0$ the solutions do not depend on the choice of regularisation, since $\mathbb{X}_{\eps, 0}$ does not depend on that choice.

Consider $w_{\eps, \delta}(T^*/2)\in L^\infty(\mathbb{T}^2)$ and let $p\in 2\mathbb{N}$ be large enough, using Theorem~\ref{thm:a-priori} there exists a (random) $C>0$ such that the $\L^p$ norm of any possible solution on the full interval
$[T^*/2,T]$ is a priori bounded by $\|w_{\eps, \delta}(T^*/2)\|_{L^p(\mathbb{T}^2)} +C $. Using the elementary consequence of Hölders inequality
\begin{equ}
\|  w  \|_{C^{-d/p}(\mathbb{T}^d) }\leq \|  w  \|_{\L^p(\mathbb{T}^d)} \ ,
\end{equ}
apply Proposition~\ref{prop:remainder equation} repeatedly to construct a solution $u_{\eps,\delta}:[0,T]\to \mathcal{D}'(\mathbb{T}^d)$ as well as to obtain continuity of the solution map as a a function of $(\eps, \delta)$ into
$\L^{0}\left[C([0,T], \mathcal{D}'(\mathbb{T}^d))\right]$.
\end{proof}

Thus, it remains to prove the parts which are specific to the corresponding results.

\begin{proof}[Remainder of the proof of Theorem~\ref{thm:main1}]
One checks \eqref{eq:assymptotics renormalsation constants} directly using the formulas \eqref{eq:first_counterterm} and \eqref{eq:second_counterterm}. 
To check Item~\ref{thm:main_item1}, it suffices to note the that $ \lim_{\eps\to 0} \alpha_{\eps,\delta}=0$ and that 
$\lim_{\eps\to 0} \bar{\alpha}_{\eps,\delta}$ exists both of which can be read off the definition directly, as well as that
$\lim_{\eps\to 0} c_{\eps,\delta}=0$ which is part of Lemma~\ref{lemma:constant c}.
For Item~\ref{thm:main_item2} one checks that
$\lim_{\delta\to 0} \bar{\alpha}_{\eps,\delta}$ and  $\lim_{\delta\to 0} c_{\eps,\delta}$ exist, the former can be read of the definition the latter 
is part of in Lemma~\ref{lemma:constant c}. In order to compare the solution obtained here to the ones of \cite{Sin23} one compares the used counterterms.
Finally, Item~\ref{thm:main_item3} was already explained in the generic part of the proof.
\end{proof}

\begin{proof}[Remainder of the proof of Theorem~\ref{thm:main_translation invariant}]
One checks \eqref{eq:assymptotics renormalsation constants_flat} directly from the definition of $\alpha^\flat_{\eps,\delta}$ Section~\ref{sec:translation invariant regularisation}. To check Item~\ref{thm:main_flat_item1}  it suffices to note that 
for $\delta>0$ it holds that $\lim_{\eps\to 0} c^\flat_{\eps,\delta}=0$ and the limit $ \bar{\alpha}^\flat_{0,\delta}:=\lim_{\eps\to 0} \bar{\alpha}^\flat_{\eps,\delta}$ exists. For the former, note that this follows directly from Lemma~\ref{lem:behaviour D flatflat and c flatflat} together with the definition of 
 $c^\flat_{\eps,\delta}$, while the latter can be read off the definition.
Item~\ref{thm:main_flat_item2} was already observed in the generic part of the proof. 
\end{proof}

\begin{proof}[Remainder of the proof of Proposition~\ref{prop:main_flat}]
The estimate \eqref{eq:bound on Dflatflat}, Item~\ref{prop:main_flatflat_item1} and Item~\ref{prop:main_flatflat_item2} are contained in Lemma~\ref{lem:behaviour D flatflat and c flatflat}. Item~\ref{prop:main_flatflat_item3} follows directly from the fact that $\mathbb{X}^{\flat\flat}_{0,\delta}= \mathbb{X}^{\flat}_{0,\delta}$ in Proposition~\ref{prop:continuity model_flat}. The first claim of Item~\ref{prop:main_flatflat_item4}, namely that $\lim_{\delta\to 0} c^{\flat \flat}_{\eps,\delta}=\lim_{\delta\to 0} c_{\eps,\delta}$ is content of Lemma~\ref{lem:behaviour D flatflat and c flatflat} while the latter claim was already explained in the generic part of the proof.
%
\end{proof}

\begin{proof}[Remainder of the proof of Theorem~\ref{thm:main_non}]
The estimate \eqref{thm:main_sharp_constant asympt} can be read off of \eqref{eq:sharp counterterm}, Item~\ref{thm:main_sharp_item1} is contained in Lemma~\ref{lem:change to sharp} and Item~\ref{thm:main_sharp_item1} was already explained in the generic part of the proof.
\end{proof}

\begin{appendix}
\section{Periodic functions}\label{ap:A}
\begin{lemma}\label{lem:improved convergence to the mean}
Suppose $f\in \L^\infty(\mathbb{R}^{d+1})$ is $(\eps\mathbb{Z})^{d}\times (\eps^2\mathbb{Z})$ periodic. Then, for $\kappa\in (0,1)$
\begin{equ}
\left\|f- \int_{[0,1]^{d+1}} f \right\|_{C_\fraks^{-\kappa}}\leq \eps^{\kappa} \|f\|_{L^\infty}\ .
\end{equ}
\end{lemma}
\begin{proof}
Let $F(z)= f(\mathcal{S}^{\eps}z)$ and $C= \int_{[0,1]^{d+1}} F=  \fint_{\mathcal{S}^{\eps^{-1}}_\fraks( [0,1]^{d+1} )} f$. We shall show that
\begin{equation}\label{eq:periodic_improved convergence}
\left|\int_{\mathbb{R}^{d+1}} (f(z) -C) \phi^{\lambda}_\star (z) \,dz \right|\lesssim \|f\|_{L^\infty} \eps^{\kappa} \lambda^{-\kappa}
\end{equation}
uniformly over $\phi\in C_c(B_1)$ satisfying $\|\phi\|_{C^\kappa_\fraks }<1$.
%
%
By translation, it suffices to consider only the case $\star=0$. 
Note that for $\eps>\lambda$ the bound holds trivially.

In the case $\eps<\lambda$ 
\begin{align*}
&\int_{\mathbb{R}^{d+1}} (f(z) -C) \phi^{\lambda} (z)\,dz = \eps^{|\fraks|} \int_{\mathbb{R}^{d+1}} F(\mathcal{S}^{\lambda^{-1}} z) \phi (\eps z) - C\\
&= \eps^{|\fraks|} \sum_{h\in \mathbb{Z}^{d+1}} \int_{\mathcal{S}^\lambda ([0,1]^{d+1}+ h)} F(\mathcal{S}^{\lambda^{-1}}z) \phi (\eps z) - \eps^{|\fraks|}\sum_{h\in \mathbb{Z}^{d+1}} C \int_{\mathcal{S}^\lambda([0,1]^{d+1}+ h)} \phi (\eps z')\, dz'\\
&= \eps^{|\fraks|} \sum_{h\in \mathbb{Z}^{d+1}} \int_{\mathcal{S}^\lambda ([0,1]^{d+1}+ h)} F(\mathcal{S}^{\lambda^{-1}}z) \phi (\eps z) \\
&\qquad- \sum_{h\in \mathbb{Z}^{d+1}}\lambda^{|\fraks|} \int_{\mathcal{S}^\lambda ([0,1]^{d+1}+ h)} F(\mathcal{S}^{\lambda^{-1}}z)\, dz  \int_{\mathcal{S}^\lambda([0,1]^{d+1}+ h)} \phi (\eps z')\, dz'\\
&= \eps^{|\fraks|} \sum_{h\in \mathbb{Z}^{d+1}} \int_{\mathcal{S}^\lambda ([0,1]^{d+1}+ h)} F(\mathcal{S}^{\lambda^{-1}}z)\left( \phi (\eps z) - \fint_{\mathcal{S}^\lambda([0,1]^{d+1}+ h)} \phi (\eps z')\, dz'\right)\\
&= \eps^{|\fraks|} \sum_{h\in \mathbb{Z}^{d+1}} \int_{\mathcal{S}^\lambda ([0,1]^{d+1}+ h)} F(\mathcal{S}^{\lambda^{-1}}z)\left( \fint_{\mathcal{S}^\lambda([0,1]^{d+1}+ h)}\phi (\eps z) - \phi (\eps z')\, dz'\right) \ .
\end{align*}
Thus
\begin{align*}
&\left|\int_{\mathbb{R}^{d+1}} (f(z) -C) \phi^{\lambda}_x (z) \right| \\
& \qquad\leq \eps^{|\fraks|+\kappa} \lambda^{-\kappa}\sum_{h\in \mathbb{Z}^{d+1}} \int_{\mathcal{S}^\lambda ([0,1]^{d+1}+ h)} F(\mathcal{S}^{\lambda^{-1}}z) \mathbf{1}_{\{ \supp(\phi(\eps \ \cdot \ )) \cap \mathcal{S}^\lambda([0,1]^{d+1}+h)\neq \emptyset   \}} \|\phi\|_{C^\kappa_\fraks }\\
& \qquad\leq \eps^{|\fraks|+\kappa} \lambda^{-(|\fraks|+\kappa)} N_{\lambda, \eps} \|F\|_{L^\infty} \|\phi\|_{C^\kappa_\fraks } \ ,
\end{align*}
where
$N_{\lambda, \eps}:= |\{h\in \mathbb{Z}^d \ : \ \supp(\phi(\eps \ \cdot \ )) \cap \mathcal{S}^\lambda([0,1]^{d+1}+h)\neq \emptyset   \}|\lesssim \eps^{-|\fraks|}\lambda^{|\fraks|}$.
\end{proof}

\section{Regularised Kernels}\label{ap:B}

Fix $\chi: \mathbb{R}^d\to [0,1]$ smooth and compactly supported on $[-2/3, 2/3]^{d}\subset \mathbb{R}^d$ such that $\sum_{k\in \mathbb{Z}^d} \chi(x+k)=1$ for all $x\in \mathbb{R}^d$ and
let 
\begin{equation}\label{eq: H_epsilon,delta}
H_{\eps,\delta}(x,t;y,s)= \sum_{k \in \mathbb{Z}^d} \chi(x-y) \int_{\mathbb{R}} \delta^{-2}\phi(\tau/\delta^2) \kappa(t-s-\tau)  \Gamma_\eps(x,t;  y+k, s+ \tau-\delta^2)  d\tau  
\end{equation}
and 
\begin{equation}\label{eq: H_epsilon,delta^flat}
H_{\eps,\delta}^\flat(x,t;y,s)= \sum_{k \in \mathbb{Z}^d} \chi(x-y) \int_{\mathbb{R}} \delta^{-2-d}\left(\phi\circ \mathcal{S}_{\fraks}^{\delta}\right) (\xi,\tau) \kappa(t-s-\tau)  \Gamma_\eps(x,t;  y+\xi+k, s+ \tau)  d\tau  \ .
\end{equation}
\begin{prop}\label{prop:H bound}
For every $R,L\in (0,1)$ there exists $C>0$ such that 
\begin{equation}\label{eq:uniform_Hbound}
\|H_{\eps, \delta} \|_{\beta;L,R}  <C , \qquad  \|H^\flat_{\eps, \delta} \|_{\beta;L,R} <C
\end{equation}
uniformly over $\eps, \delta \in [0,1]^2$, $\beta\in (1,2]$. 
For $R,L\in (0,1)$, $\beta\in (1,2]$ such that $R+2-\beta<1$
\begin{equation}\label{eq:uniform_Hconv_in_delta}
\|H_{\eps, \delta}- H_{\eps,0} \|_{\beta;L,R} \lesssim \delta^{2-\beta}\  , \qquad \|H^\flat_{\eps, \delta}- H^\flat_{\eps,0} \|_{\beta;L,R} \lesssim \delta^{2-\beta}
\end{equation}
uniformly over $\eps, \delta \in [0,1]^2$.
Finally, for $0<R+L<1$ it holds that
\begin{equation}\label{eq:uniform_Hconv_in_epsilon}
\|H_{\eps, \delta}- H_{0, \delta} \|_{\beta;L,R} \lesssim \eps^{2-\beta}\vee \eps^{1-R-L}  , \qquad  \|H^\flat_{\eps, \delta}- H^\flat_{0, \delta} \|_{\beta;L,R} \lesssim \eps^{2-\beta}\vee \eps^{1-R-L} 
\end{equation}
uniformly over $\beta\in [0,2]$. 
\end{prop}

\begin{proof}
Set for $\varphi^{2^{-n}}$ is as in the proof of Proposition~\ref{prop:uniform bounds on modified kernels} and $n\geq 1$ and
\begin{equation}\label{eq:decompose H_epsilon delta}
H_{\eps,\delta;n}(x,t,y,s)= \varphi^{2^{-n}}(z-z') 	H_{\eps,\delta}(z;z')\ , 	 \qquad H^\flat_{\eps,\delta;n}(x,t,y,s)= \varphi^{2^{-n}}(z-z') 	H^\flat_{\eps,\delta}(z;z')\ , 
\end{equation}
well as $H_{\eps,\delta;0}(x,t,y,s)= \sum_{k \in \mathbb{Z}^d} H_{\eps,\delta}(x,t;y+k,s) - \sum_{n=1}^\infty H_{\eps,\delta;n}(x,t,y,s)$ and $H^\flat_{\eps,\delta;0}(x,t,y,s) \newline = \sum_{k \in \mathbb{Z}^d} H^\flat_{\eps,\delta}(x,t;y+k,s) - \sum_{n=1}^\infty H^\flat_{\eps,\delta;n}(x,t,y,s)$. 
In particular for both $H$ and $H^\flat$ we thus obtain a decomposition as in the definition of the spaces $\bfK^\beta_{L,R}$. 
Inequalities \eqref{eq:uniform_Hbound} follow from Theorem~\ref{thm:uniform}, Proposition~\ref{prop:uniform_hölder_bound} and Proposition~\ref{prop:uniform hölder bound on kernel}.
Similarly, \eqref{eq:uniform_Hconv_in_delta} follows from Proposition~\ref{prop:uniform_hölder_bound} and Proposition~\ref{prop:uniform hölder bound on kernel}.
Finally, \eqref{eq:uniform_Hconv_in_epsilon} follows as \eqref{eq:corrrectet kernel convergence for RS_unmodified kernel} in Proposition~\ref{prop:kernels for fixed point convergence}.
%
%
\end{proof}

\begin{equation}\label{eq:H_epsilon,delta^sharp}
H^\sharp_{\eps,\delta}(x,t;y,s)=\int_{{\mathbb{R}^{d+1}}} \kappa(t-\tau)\Gamma_\eps(x,t, \zeta, \tau) \varrho^{\eps,\delta} (\zeta, \tau; y,s) d\zeta d\tau \ . 
\end{equation}

\begin{prop}\label{prop:H^sharp bound}
For every $R,L\in (0,1)$ 
\begin{equation}\label{eq:uniform_Hsharpbound}
\|H^\sharp_{\eps, \delta} \|_{\beta;L,R}  \lesssim \delta/\eps \vee 1
\end{equation}
uniformly over $\eps, \delta \in [0,1]^2$, $\beta\in (1,2]$. 
For $R,L\in (0,1)$, $\beta\in (1,2]$ such that $R+2-\beta<1$
\begin{equation}\label{eq:uniform_Hsharpconv_in_delta}
\|H^\sharp_{\eps, \delta}- H^\sharp_{\eps,0} \|_{\beta;L,R} \lesssim (\delta/\eps \vee 1) \delta^{2-\beta}\  
\end{equation}
uniformly over $\eps, \delta \in [0,1]^2$.
Finally, for $0<R+L<1$ it holds that 
\begin{equation}\label{eq:uniform_Hsharpconv_in_epsilon}
\|H^\sharp_{\eps, \delta}- H^\sharp_{0, \delta} \|_{\beta;L,R} \lesssim (\delta/\eps \vee 1) \eps^{2-\beta}\vee \eps^{1-R-L}  
\end{equation}
uniformly over $\beta\in [0,2]$. 
\end{prop}
\begin{proof}
We decompose $H^\sharp_{\eps, \delta}$ analogously to \eqref{eq:decompose H_epsilon delta}. Then, \eqref{eq:uniform_Hsharpbound} and\eqref{eq:uniform_Hsharpconv_in_delta} follow by a more tedious but similar computation to \eqref{eq:uniform_Hbound} and \eqref{eq:uniform_Hconv_in_delta}.
To see \eqref{eq:uniform_Hsharpconv_in_epsilon}, we introduce
\begin{equ}
\bar{H}^\sharp_{\eps,\delta}(x,t;y,s)=\int_{{\mathbb{R}^{d+1}}} \kappa(t-\tau)\bar{\Gamma}(x,t, \zeta, \tau) \varrho^{\eps,\delta} (\zeta, \tau; y,s) d\zeta d\tau \
\end{equ}
and estimate separately the two terms
$\|H^\sharp_{\eps, \delta}- \bar{H}^\sharp_{\eps, \delta} \|_{\beta;L,R}$ and
$\|\bar{H}^\sharp_{\eps, \delta}- H^\sharp_{0, \delta} \|_{\beta;L,R}$.
\end{proof}
\section{Convolution}\label{ap:C}

\begin{prop}\label{prop:convolution estimate}
Let $L,L', R, R'<1$ and $\beta, \beta' \in (0, |\fraks|)$. 
Assume that $|\fraks|-\beta-\beta'>0$, $L<\beta$, $L'<\beta'$, then 
for $H\in \bfK^\beta_{L,R}$, $H'\in \bfK^{\beta'}_{L',R'}$ 
  it holds that
\begin{equ}
H\star H':=  \int_{\mathbb{R}^{d+1} } H(z,w) H(z', w)
\end{equ}
belongs to $\bfK^{\beta+\beta'}_{L,L'}$
and furthermore it holds that $\|H\star H'\|_{\beta+\beta'; L, L'}\lesssim \| F\|_{\beta;L,0} \| F\|_{\beta;L,0}$.
\end{prop}

\begin{proof}
We write $H= \sum_{n=0}^\infty H_n$, $H'= \sum_{m=0}^\infty H'_m$, thus
\begin{equ}
H\star H':= \sum_{m\geq 0} (H\star_> H')_m +(H\star_< H')_m + H_m\star H_m\ ,
\end{equ}
where
\begin{equ}
(H\star_> H')_m:= \sum_{n>m}H_n\star H_m, \qquad (H\star_< H')_m:= \sum_{n>m}H_m\star H_n \ .
\end{equ}
Note that $(H\star H')_m$ is supported on $\{(z,z')\in (\mathbb{R}^{d+1})^{\times 2}\ : \ |z-z'|_\fraks<2^{-m+1} \}$.
For $\gamma \in \{0,L\}$, $\gamma' \in \{0,L\}$ with the understanding that $\| \ \cdot\  \|_{C^{0,0}}= \| \cdot \ \|_{L^\infty}$, we find
\begin{itemize}
\item $\|H_m\star H'_m\|_{C^{\gamma,\gamma'}} \lesssim 2^{-m|\fraks| } 2^{m(|\fraks|-\beta+\gamma) } 2^{m(|\fraks|-\beta'+\gamma') }  = 2^{-m(|\fraks|-\beta-\beta'+\gamma +\gamma') }$
\item $\|(H\star_> H')_m\|_{C^{\gamma,\gamma'}}\lesssim \sum_{n>m}2^{-n|\fraks| } 2^{n(|\fraks|-\beta+ \gamma) } 2^{m(|\fraks|-\beta') } \lesssim 2^{-m(|\fraks|-\beta-\beta'+\gamma+\gamma')}$.
\item Similarly, $\|(H\star_< H')_m\|_{C^{\gamma,\gamma'}} \lesssim 2^{-m(|\fraks|-\beta-\beta'+\gamma+\gamma')}$. 
\end{itemize}
Thus, writing $H\star H'= \sum_{m\geq 0} (H\star H')_{m}$ where
\begin{equ}
(H\star H')_{0}= \sum_{i=0}^1\ (H\star_> H')_i +(H\star_< H')_i + H_i\star H_i
\end{equ}
and for $m\geq 1$
\begin{equ}
(H\star H')_{m}= (H\star_> H')_{m+1} +(H\star_< H')_{m+1} + H_{m+1}\star H_{m+1} \qquad
\end{equ}
this concludes the proof.
\end{proof}

\section{Kernel Index}\label{ap:D}

\begin{center}
\begin{longtable}{p{.12\textwidth}p{.60\textwidth}p{.10\textwidth}}
\toprule
Kernel & Description & Ref. \\
\midrule
\endhead
\bottomrule
\endfoot
$\Gamma$ & fundamental solution of $\partial_t - \nabla \cdot A(\cdot)\nabla$  & Page~\pageref{eq:dif_operator}\\
$Z$ & `frozen' in second variable kernel  & Eq.~\ref{eq:exlicit Z_0}\\
$\bar{Z} $ & `frozen' in first variable kernel  & Eq.~\ref{eq:explicit bar Z}\\
$\bar{Z}_\eps$ & defined as $\bar{Z}(x,t; \zeta,\tau) $ but for $A(\cdot, \cdot)$ replaces by  $A(\cdot/\eps, \cdot/\eps^2)$  & Page~\pageref{eq:scaling_property} \\
$\Gamma_\eps$  & fundamental solution of $\partial_t - \nabla \cdot A(\cdot/\eps, \cdot/\eps^2)\nabla $ & Page~\pageref{eq:scaling_property}\\
$\bar{\Gamma}$ & fundamental solution of the homogenised operator & Page~\pageref{eq:adjoint} \\
$\Gamma^{I,J}_\eps$ & corrected fundamental solution for $I,J\in \{0,1\}$  & Eq.~\ref{eq:GammaIJ} \\
$ \bar{\Gamma}_\eps $ & truncated variant of $\bar{\Gamma}$ & Page~\pageref{eq:Gamma_tilde} \\
$\tilde{\Gamma}_\eps$ & truncated variant of $\Gamma^{1,1}$ & Eq.~\ref{eq:Gamma_tilde}\\
$\bar{K}$&  periodisation of $\bar{\Gamma}$ & Eq.~\ref{eq:periodising kernels}\\
${K}_{\eps}$&  periodisation of  ${\Gamma}_{\eps}$ & Eq.~\ref{eq:periodising kernels}\\ 
$\tilde{K}_{\eps}$&  periodisation  of $\tilde{\Gamma}_{\eps}$ &  Eq.~\ref{eq:periodising kernels}\\
$\mathcal{Z}_\eps$ & kernel agreeing on small scales with $\bar{Z}_\eps$ and large scales with $\bar{\Gamma}$ & Eq.~\ref{eq:mathcalZ}
\end{longtable}
\end{center}
\end{appendix}

	\bibliographystyle{Martin}
	\bibliography{Periodic.bib}

\end{document}